\documentclass[11pt]{article}

\usepackage[utf8]{inputenc} % allow utf-8 input
\usepackage{url}            % simple URL typesetting
\usepackage{wrapfig}
\usepackage[authoryear]{natbib}
\usepackage{breakurl}
\usepackage{mathtools}

\usepackage{xcolor}
\definecolor{dgreen}{rgb}{0, .4, 0}

\usepackage{url}% for url's in bib
\usepackage[colorlinks,linkcolor=magenta,citecolor=blue]{hyperref}    %, pagebackref=true

\usepackage{graphicx}

\usepackage{amsmath,amssymb,fullpage,graphicx}
\usepackage{comment}
\usepackage{color}
\usepackage{setspace}
\usepackage{pdflscape}

\usepackage{algorithm}
\usepackage{algpseudocode}
\usepackage{amssymb,amsmath,amsthm}
\usepackage[english]{babel}

\usepackage{macros}

\usepackage{setspace}

\usepackage{amsmath}
\makeatletter
\let\save@mathaccent\mathaccent
\newcommand*\if@single[3]{%
  \setbox0\hbox{${\mathaccent"0362{#1}}^H$}%
  \setbox2\hbox{${\mathaccent"0362{\kern0pt#1}}^H$}%
  \ifdim\ht0=\ht2 #3\else #2\fi
  }
\newcommand*\rel@kern[1]{\kern#1\dimexpr\macc@kerna}
\newcommand*\widebar[1]{\@ifnextchar^{{\wide@bar{#1}{0}}}{\wide@bar{#1}{1}}}
\newcommand*\wide@bar[2]{\if@single{#1}{\wide@bar@{#1}{#2}{1}}{\wide@bar@{#1}{#2}{2}}}
\newcommand*\wide@bar@[3]{%
  \begingroup
  \def\mathaccent##1##2{%
    \let\mathaccent\save@mathaccent
    \if#32 \let\macc@nucleus\first@char \fi
    \setbox\z@\hbox{$\macc@style{\macc@nucleus}_{}$}%
    \setbox\tw@\hbox{$\macc@style{\macc@nucleus}{}_{}$}%
    \dimen@\wd\tw@
    \advance\dimen@-\wd\z@
    \divide\dimen@ 3
    \@tempdima\wd\tw@
    \advance\@tempdima-\scriptspace
    \divide\@tempdima 10
    \advance\dimen@-\@tempdima
    \ifdim\dimen@>\z@ \dimen@0pt\fi
    \rel@kern{0.6}\kern-\dimen@
    \if#31
      \overline{\rel@kern{-0.6}\kern\dimen@\macc@nucleus\rel@kern{0.4}\kern\dimen@}%
      \advance\dimen@0.4\dimexpr\macc@kerna
      \let\final@kern#2%
      \ifdim\dimen@<\z@ \let\final@kern1\fi
      \if\final@kern1 \kern-\dimen@\fi
    \else
      \overline{\rel@kern{-0.6}\kern\dimen@#1}%
    \fi
  }%
  \macc@depth\@ne
  \let\math@bgroup\@empty \let\math@egroup\macc@set@skewchar
  \mathsurround\z@ \frozen@everymath{\mathgroup\macc@group\relax}%
  \macc@set@skewchar\relax
  \let\mathaccentV\macc@nested@a
  \if#31
    \macc@nested@a\relax111{#1}%
  \else
    \def\gobble@till@marker##1\endmarker{}%
    \futurelet\first@char\gobble@till@marker#1\endmarker
    \ifcat\noexpand\first@char A\else
      \def\first@char{}%
    \fi
    \macc@nested@a\relax111{\first@char}%
  \fi
  \endgroup
}
\makeatother

\begin{document}

\begin{center} {\LARGE{\bf{Martingale Methods for Sequential  Estimation \\[0.07in] of  Convex  Functionals and Divergences 
}}}

\vspace*{.1in}

{\large{

%\vspace*{.05in}

\begin{tabular}{cccc}
Tudor Manole$^{1}$, Aaditya Ramdas$^{12}$
\end{tabular}

{
\vspace*{.1in}
\begin{tabular}{c}
				$^1$ Department of Statistics and Data Science\\
				$^2$ Machine Learning Department\\  
				Carnegie Mellon University\\ 
				\texttt{\{tmanole,aramdas\}@andrew.cmu.edu }
			\end{tabular}
}

\vspace*{.1in}

}}

\today

\vspace*{.1in}
\begin{abstract}
We present a unified technique for sequential estimation of convex divergences between distributions, including integral probability metrics like the kernel maximum mean discrepancy,  $\varphi$-divergences like the Kullback-Leibler divergence, 
and optimal transport costs, 
such as powers of 
Wasserstein distances. This is achieved by observing that empirical convex divergences are (partially ordered) reverse submartingales with respect to the exchangeable filtration, coupled with
maximal inequalities for such processes. These techniques appear to be complementary and powerful additions to the existing literature on both confidence sequences and convex divergences. We construct an offline-to-sequential device that converts a wide
array of existing offline concentration inequalities into time-uniform confidence sequences that can be continuously monitored, providing valid tests or confidence intervals at arbitrary stopping times. The resulting sequential bounds pay only an iterated logarithmic price over the corresponding fixed-time bounds, retaining the same dependence on problem parameters (like dimension or alphabet size if applicable). 
These results are also  applicable
to more general convex functionals---like the negative differential entropy, suprema of empirical processes,
and V-Statistics---and to more general processes satisfying a key leave-one-out property.
\end{abstract}
\end{center}

%%%%% Table of contents %%%%%%%%
 \begingroup
 \hypersetup{linkcolor=black}
 {\renewcommand{\baselinestretch}{0.75} 
 \scriptsize  \tableofcontents}
 \endgroup
\section{Introduction}
\label{sec:introduction}
Divergences between 
probability distributions 
arise pervasively in 
information theory, statistics,
%probability theory, 
and machine learning \citep{liese2006,sriperumbudur2012}. 
Common examples %of divergences
include $\varphi$-divergences, 
such as the Kullback-Leibler divergence,
integral probability metrics (IPMs), 
such as the kernel maximum mean discrepancy,
and optimal transport costs, 
such as powers of 
Wasserstein distances. 
Increasingly many applications 
employ such quantities as methodological
tools, in which it is of interest
to estimate a given divergence
between unknown probability
distributions on the basis of random
samples thereof. 
Statistical inference 
is well-studied for many divergences  when
the data is available in a fixed batch
ahead of the user's analysis. 
When data is instead
collected sequentially in time, 
such methods are typically invalid
for repeatedly assessing uncertainty
of a divergence estimate, 
and may therefore provide overly 
optimistic confidence.

The aim of this paper is to develop
rigorous {\it sequential} uncertainty
quantification methods for a large
class of divergences between probability distributions. 
Throughout the sequel, we denote
by $D : \calP(\calX) \times \calP(\calX)\to \bbR_+$
a generic convex divergence, where 
$\calP(\calX)$ denotes the set of Borel probability
measures over a set $\calX \subseteq \bbR^d$. Recall that the functional $D$ is said to be convex 
if
for all measures $\mu_1, \mu_2, \nu_1, \nu_2 \in \calP(\calX)$
and all $\lambda \in [0,1]$,
\begin{equation} 
\label{eq:convex_divergence}
D\big(\lambda \mu_1 + (1-\lambda)\mu_2 ~ \|~  \lambda \nu_1 + (1-\lambda) \nu_2\big) \leq 
\lambda D(\mu_1 \| \nu_1) + (1-\lambda) D(\mu_2 \| \nu_2).
\end{equation}
This property can be  verified for IPMs, 
$\varphi$-divergences, 
and optimal transport costs 
 (see Section~\ref{sec:background_divergences}).
% Our results
% apply more generally to convex functionals
% $\Phi:\calP(\calX) \to \bbR$, 
% satisfying $\Phi(\lambda \mu_1 + (1-\lambda) \mu_2)
%  \leq \lambda \Phi(\mu_1) + (1-\lambda) \Phi(\mu_2)$
% for all $\lambda \in [0,1]$ and $\mu_1,\mu_2 \in \calP(\calX)$. 
% Note that a convex divergence with either argument kept fixed is a convex functional.
% [Say $\Phi(P) = D(P\|Q)$ is
% such a functional]

Given two independent sequences $(X_t)_{t=1}^\infty$
and $(Y_s)_{s=1}^\infty$
of i.i.d. observations arising
respectively
from unknown distributions 
$P,Q \in \calP(\calX)$, 
we aim to construct 
a sequence of confidence
intervals $(C_{ts})_{t,s=1}^\infty$
with the uniform coverage property
\begin{equation}
\label{eq:confseq}
\bbP\big(\forall t,s \geq 1: 
 D(P\|Q) \in C_{ts}\big) \geq 1-\delta,
\end{equation}
for some pre-specified level
$\delta \in (0,1)$.
Such a sequence $(C_{ts})_{t,s=1}^\infty$
is called a {\em confidence
sequence}, and differs
from the standard notion of confidence
interval through the uniformity
in times $t,s$ of the probability 
in equation \eqref{eq:confseq}.
As stated precisely in 
Section~\ref{sec:choice_filtrations},
the guarantee \eqref{eq:confseq}
is equivalent to the requirement
that for all stopping times 
$(\tau,\sigma)$, 
\begin{equation}
\label{eq:confseq_stopping_time}
\bbP(D(P\|Q) \in C_{\tau\sigma }) \geq 1-\delta.
\end{equation}
The scope of potential
applications of such confidence sequences
is far-reaching. For instance, 
a confidence sequence
$C_{ts}$ directly
gives rise to a sequential
two-sample 
test for the null hypothesis
$H_0: P=Q$, where the null
is rejected when $0 \not\in C_{ts}$. 
Fixed-time, two-sample testing 
is a classical problem which
continues to receive
a wealth of attention, 
but \emph{sequential, nonparametric} two-sample testing is relatively less explored;
two exceptions include~\cite{balsubramani2016},
\cite{lheritier2018}. 
We also
note that confidence sequences
for divergences can sometimes  
lead to confidence sequences for other
estimands of interest, as illustrated in Section~\ref{sec:smooth}. 
% hinging upon its continuity with
% respect to the 1-Wasserstein
% distance \citep{polyanskiy2016}. 
Though our work is 
motivated by such practical 
applications, 
our results
can also be viewed from a purely
theoretical standpoint
as deriving  
concentration inequalities
for divergences 
which hold uniformly over time.
We are not aware
of analogous 
time-uniform concentration inequalities in
the literature for the majority of 
divergences studied explicitly
in this paper. 

\textbf{Our Contributions.}
The primary contribution
of our work is to provide
a general recipe for deriving
confidence sequences for 
 convex divergences $D$. 
Our key
observation is that 
the process
\begin{equation}
\label{eq:process_Mts}
M_{ts} = D(P_t\|Q_s) - D(P\|Q), \quad t,s \geq 1,
\end{equation}
is a partially-ordered reverse submartingale, with respect
to the so-called exchangeable filtration
introduced below.
Here, 
$P_t = (1/t) \sum_{i=1}^t \delta_{X_i}$
and $Q_s = (1/s) \sum_{i=1}^s \delta_{Y_j}$
denote empirical measures.
A related
property was previously identified 
by \cite{pollard1981a} for suprema of empirical
processes. 
% , with the convention
% $P_0=P$ and $Q_0=Q$. 
This reverse submartingale property
allows us to apply maximal inequalities
to (functions of) $(M_{ts})_{t,s=1}^\infty$,
which will lead to confidence sequences
for $D(P\|Q)$ based on the plug-in estimator
$D(P_t\|Q_s)$. 
We note
that this estimator
is inconsistent for divergences
requiring the absolute continuity 
of the probability measures being compared,
such as $\varphi$-divergences for 
distributions
supported over $\bbR^d$. 
We therefore extend our results
by showing that 
the process in equation \eqref{eq:process_Mts}
continues to be a reverse submartingale
in each of its indices
when $P_t$ and $Q_s$
are replaced by their smoothed 
counterparts, $P_t\star \calK_\sigma$
and $Q_s \star \calK_\sigma$, 
where $\calK_\sigma$ denotes a kernel
with bandwidth $\sigma$.

We illustrate these findings
by deriving explicit confidence
sequences for the kernel 
Maximum Mean Discrepancy (MMD), 
Wasserstein distances, Total Variation distance, 
and Kullback-Leibler divergence, among others, some for distributions over finite alphabets and others for arbitrary distributions. In all cases, we take care to track the effect of dimensionality, matching the best known rates is non-sequential settings.
To the best of our knowledge, 
there are
no other existing confidence sequences
for these quantities, apart
from the (statistically suboptimal but computationally efficient) linear-time kernel MMD 
\citep{balsubramani2016}. 
We also derive a sequential analogue
of the celebrated Dvoretzky-Kiefer-Wolfowitz
inequality~\citep{dvoretzky1956,massart1990} 
quite differently from both \cite{howard2019} and a recent preprint by \cite{odalric-ambrym2020},
and demonstrate how these
results can be used to obtain confidence
sequences for convex functionals which do not necessarily arise from divergences.

\textbf{Outline.}
We organize this manuscript 
as follows. 
In Section~\ref{sec:background}, 
we provide background on martingales
and maximal inequalities, which form 
the main technical tools used
throughout this paper. We also
provide background on sequential analysis
and on several classes
of convex divergences. 
In Section~\ref{sec:main}, we
state our main results and derive a general confidence
sequence for $D(P\|Q)$ on the basis of the process
$(M_{ts})$, which also applies to more general convex functionals. We 
explicitly illustrate
its application to a wide range
of commonly-used divergences in Section~\ref{sec:applications}.
We close with a conclusion 
in Section \ref{sec:conclusion}. Most proofs of our 
results are relegated to Appendices~\ref{app:additional_lemmas}--\ref{sec:forward_submart}
of the Supplementary Material.

\textbf{Notation.} For a vector
$x \in \bbR^d$, and for some integer
$p \geq 1$,
$\norm{x}_p = (\sum_{i=1}^d |x_i|^p)^{1 /p}$
denotes the $\ell_p$ norm of $x$, and $\norm{x}_\infty = \max_{1 \leq i \leq d} |x_i|$
denotes its $\ell_\infty$ norm.
We also denote the $L^\infty$ norm
of a function $f:\bbR^d \to \bbR$
by $\norm{f}_\infty = \esssup_{x \in \bbR^d}|f(x)|$. 
We abbreviate sequences $(A_t)_{t=1}^\infty$
of sets, functions, or real numbers
by $(A_t)$ when doing so causes no confusion.
For any $a,b \in \bbR$, $a\vee b := \max\{a,b\}$
and $a\wedge b := \min\{a,b\}$.
The diameter of a set $A \subseteq \bbR^d$
is given by $\diam(A) = \sup\{\norm{x-y}_2:x,y\in A\}$.
The Dirac measure placing
mass at a point $x \in \bbR^d$ is
denoted $\delta_x$. The Borel $\sigma$-algebra
on $\calX \subseteq \bbR^d$ is denoted $\bbB(\calX)$.
The floor and ceiling
of a real number $x \in \bbR$ are respectively
denoted $\lfloor x\rfloor$ and 
$\lceil x \rceil$. We write
$\widebar\bbR = \bbR\cup \{\infty\}, \bbR_{+} = [0,\infty)$,
$\bbN=\{1, 2, \dots\}$ and 
$\bbN_0 = \{0, 1, \dots\}$. 
% Given an integer $T > 0$, we write
% $[T] = \{1, \dots, T\}$ and $[T]_0 = \{0, 1, \dots, T\}$.
The convolution of functions $f,g:\bbR^d \to \bbR$ 
is denoted
by $(f\star g)(x) = \int f(x-y)g(y)dy$.
Furthermore, the convolution
of two Borel probability
measures $P,Q$ is the measure 
$(P\star Q)(B) = \int I_B(x+y)dP(x)dQ(y)$
for all $B \in \bbB(\bbR^d)$, 
where $I_B(x) = I(x \in B)$ is the indicator function of $B$. 
The convolution of  
$P$  with $f$ is the function 
$(P \star f)(x) = \int f(x-y)dP(y)$.
For a convex function
$\varphi: \bbR \to \widebar\bbR$, 
its Legendre-Fenchel transform
is denoted
$\varphi^*(y) = \sup_{x \in \bbR} \big\{x^\top y - \varphi(x)\big\}$. 
% We shall use
% the same symbol to denote the Cram\'er
% transform over an interval $I\subseteq \bbR$, 
% $\varphi^*(y) = \sup_{x \in I} \big\{x^\top y - \varphi(x)\big\}$, when
% doing so causes no confusion.
The logarithm of base $b > 1$ is denoted
$\log_b(x)=\log x / \log b$, where $\log$ is
the natural logarithm. 
An $\bbR^d$-valued random variable $Y$ 
 is said to
be $\sigma^2$-sub-Gaussian if 
$\bbE \exp(\lambda^\top (Y-\bbE Y)) \leq
\exp(\norm\lambda_2^2 \sigma^2/2)$ for all $\lambda \in \bbR^d$. 
Given  a probability space $(\Omega, \calF, \bbP)$,
the Lebesgue spaces of equivalence
classes 
are denoted $L^p(\bbP)$ for $p \geq 1$.
For any two sub-$\sigma$-algebras $\calG, \calH \subseteq \calF$, 
we use the standard notation for joins and intersections of 
$\sigma$-algebras, respectively given by
$\calG \bigvee \calH := \sigma(\calG \cup \calH)$
and $\calG \bigwedge \calH := \calG \cap \calH$
to emphasize that they both result in $\sigma$-algebras.

\section{Background}
\label{sec:background}

\subsection{IPMs, 
Optimal Transport Costs, and $\varphi$-Divergences}
\label{sec:background_divergences}

Let $\calX \subseteq \bbR^d$ 
and let $\calA \subseteq 
\calP(\calX)\times \calP(\calX)$ 
be a set of pairs of 
probability measures.
Throughout this paper,
we use the term divergence
to refer to any 
map $D:\calA \to \bbR$
which is nonnegative and 
satisfies 
$D(P\|Q) = 0$ if  
$P = Q$, for all $(P,Q) \in \calA$. 
When the divergence $D$ is convex, 
we extend the domain
of $D$ from $\calA$
to the entire set $\calP(\calX) \times \calP(\calX)$ by letting
$D$ take the value~$\infty$ wherever it is not defined---as
such, convex divergences
will always be understood
as maps
$D:\calP(\calX) \times \calP(\calX)
\to\widebar \bbR_+:= \bbR_+\cup \{\infty\}.$
% We will sometimes
% require $D$ to satisfy the triangle inequality,
% but we will only make this assumption
% explicitly when needed. 

We consider sequential estimation 
of generic convex divergences,
but  
the following three classes will be used as recurring examples
throughout our development. 
Let $P,Q \in \calP(\calX)$ in the sequel. 
\begin{itemize}
\item \textbf{Integral Probability Metrics (IPMs).} Let $\calJ$ denote a set of Borel-measurable,
real-valued functions on $\calX$. The IPM  \citep{muller1997} associated with $\calJ$ is given by
\begin{equation} 
\label{eq:ipm}
D_\calJ(P \| Q) = \sup_{f \in \calJ} \int f d(P-Q).
\end{equation}
For instance, when $P$ and $Q$ have supports contained
in $\bbR$, the class of indicator functions
$\calJ = \{I_{(-\infty,x]}: x \in \bbR\}$ 
gives rise to the Kolmogorov-Smirnov distance, 
$D_\calJ(P\|Q) = \norm{F-G}_\infty$. 
When $\calJ$ is the unit ball of a reproducing
kernel Hilbert space, 
$D_\calJ$ is called the (kernel)
Maximum Mean Discrepancy (\cite{gretton2012};
see also Section \ref{sec:mmd}).
When $\calJ$ is the set of Borel-measurable
maps $f:\bbR^d \to \bbR$ satisfying $\norm f_\infty \leq 1$, 
$D_\calJ$ becomes the total variation distance,
\begin{equation}
\label{eq:tv}
\norm{P-Q}_{\mathrm {TV}} = \sup_{A \in \bbB(\calX)} |P(A) - Q(A)|.
\end{equation}
When the function space $\calJ$
is sufficiently small, 
$D_\calJ(P_t\|Q_s)$ 
is a consistent estimator of 
$D_\calJ(P\|Q)$, and will form
the basis of our confidence sequences. We refer to 
\cite{sriperumbudur2012}
for a study of 
convergence rates for 
such plug-in estimators.

\item \textbf{Optimal Transport Costs.} 
Let $\Pi(P,Q)$ denote the set of joint probability distributions $\pi \in \calP(\calX \times \calX)$
with marginals $P,Q$, that is, satisfying $\pi(B \times  \calX) = P(B)$ and $\pi(\calX \times B) = Q(B)$ for all
$B \in \bbB(\calX)$. 
Given a nonnegative cost function
$c:\calX\times \calX \to \bbR_+$,
the optimal transport cost between $P$ and $Q$
is given by
\begin{equation} 
\label{eq:ot_cost} 
\calT_c(P,Q) = \inf_{\pi\in\Pi(P,Q)} \int c(x,y) d\pi(x,y).
\end{equation}
$\calT_c(P,Q)$ admits the natural interpretation
of measuring the work required to 
couple the measures $P$ and $Q$---we refer to 
\cite{panaretos2019a,villani2003} 
for surveys. When $c$ takes the form
$c=d^{p}$ for some $p\geq 1$
and some metric $d:\calX\times\calX \to \bbR_+$,
the quantity $W_p:= \calT_{d^p}^{1/p}$
becomes a metric on $\calP_p(\calX)= \{P\in \calP(\calX): \int d(x,x_0)^p dP(x)<\infty, x_0 \in \calX\}$ known as the $p$-Wasserstein distance.

The minimization problem \eqref{eq:ot_cost} is an infinite-dimensional linear program, and admits a dual formulation
known as the Kantorovich dual problem. When $c$ is lower semi-continuous,
strong duality holds (\cite{villani2003}, Theorem 1.3), 
leading to the following representation,
\begin{equation} 
\label{eq:kantorovich}
\calT_c(P,Q) = \sup_{(f,g) \in \calM_c} \int f dP + \int g d Q,
\end{equation}
where $\calM_c\equiv \calM_c(P,Q)$ denotes the  set of pairs $(f,g)$ of
functions $f\in L^1(P)$ and $g \in L^1(Q)$
such that $f(x) + g(y) \leq c(x,y)$
for $P$-almost all $x \in \calX$ and $Q$-almost all $y \in \calX$. 
Moreover, when the cost function $c$ is bounded,
the supremum
in the above display is achieved and 
the pairs $(f,g) \in \calM_c$ may be restricted to take the
form $0 \leq f \leq \norm c_\infty$, $-\norm c_\infty \leq g \leq 0$.
Finally, when $c=d$ is itself a metric, 
the Kantorovich-Rubinstein formula implies
that $\calM_c$ may be further reduced to the set
of pairs $\{(f,-f):f \in \calJ\}$, with 
$\calJ$ denoting the set of 1-Lipschitz functions
with respect to~$d$. In this case, 
$\calT_d$ is precisely the IPM generated by~$\calJ$.
We refer to 
\cite{fournier2015},
% \cite{weed2019},
% \cite{lei2020},
\cite{sommerfeld2018},
\cite{niles-weed2019}, 
\cite{chizat2020}, 
\cite{manole2021a}, 
\cite{hundrieser2022},
and references therein 
for convergence
rates and fixed-time concentration inequalities
for plug-in estimators
of optimal transport costs and Wasserstein distances.

\item \textbf{$\varphi$-Divergences.} Let $\varphi:\bbR \to \bbR$ be a convex function, 
and let $\nu\in \calP(\calX)$ be a $\sigma$-finite
measure which
dominates both $P$ and $Q$ (for instance, 
$\nu = (P+Q)/2$). Let $p=dP/d\nu$ and $q = dQ/d\nu$ be the respective densities. 
Then, the $\varphi$-divergence \citep{ali1966}
between $P$ and $Q$ is given by  
\begin{equation} 
\label{eq:phi_divergence} 
D_\varphi(P\|Q) = \int_{q>0} \varphi\left(\frac{p}{q}\right) dQ
 + P(q =0) \lim_{x\to\infty}\frac{\varphi(x)}{x},
\end{equation}
with the convention that the second term of the above
display is equal to zero whenever $P(q=0)=0$,
which is in particular the case if $P \ll Q$. 
For instance, assuming the latter
condition holds, the Kullback-Leibler divergence
$\KL(P\|Q) = \int \log\left(\frac{dP}{dQ}\right)dP$ 
corresponds to the map
$\varphi(x)=x\log x$. 
The total variation
distance in equation 
\eqref{eq:tv} is the unique nontrivial
IPM
which is also a $\varphi$-divergence
\citep{sriperumbudur2012}.
$\varphi$-divergences also admit the variational
representation
\begin{equation} 
\label{eq:phi_div_variational}
D_\varphi(P\|Q)
 \geq \sup_{g \in \calJ} \int 
 \big[g dQ - (\varphi^*\circ g)dP],
 \end{equation} 
for any collection $\calJ$
of functions mapping $\calX$ to $\bbR$.
Equality holds in the above display
if and only
if the subdifferential
$\partial \varphi(dQ/dP)$ 
contains an element of $\calJ$ \citep{nguyen2010}.

Unlike most common IPMs and optimal transport costs,
$\varphi$-divergences are typically
uninformative when $P$ is not absolutely continuous
with respect to $Q$, as the expression~\eqref{eq:phi_divergence} 
becomes dominated by its second term.
This fact sometimes prohibits the estimation
of $\varphi$-divergences via
the plug-in estimator $D(P_t\|Q_s)$---for instance, 
$P_t$ is almost surely not absolutely continuous with
respect to $Q_s$ when
$P$ and $Q$ are both absolutely continuous
with respect to the Lebesgue measure.
One exception is the situation
where $P$ and $Q$ are supported
on countable sets, in which case
\cite{berend2013},
\cite{agrawal2020}, 
\cite{guo2020}, 
\cite{cohen2020},
\cite{han2015},
\cite{kamath2015},
study concentration and convergence
rates of the empirical measure
$P_t$ under the Kullback-Leibler and
Total Variation  divergences.
We develop time-uniform bounds which 
build upon these results
in Section \ref{sec:phi_divergence}.
For distributions $P$ and $Q$ which
are not countably-supported, 
distinct estimators
have been developed by 
\cite{nguyen2010},
\cite{poczos2012},
\cite{sricharan2010},
\cite{krishnamurthy2015},
\cite{rubenstein2019},
\cite{singh2014}, 
\cite{wang2005},
\cite{berrett2019},
and references therein.

\end{itemize} 

The following Lemma is standard, 
and stated without proof.
\begin{lemma}
\label{lem:convex}
For any class $\calJ$  of Borel-measurable functions from $\bbR^d$ to $\bbR$,
the IPM generated by $\calJ$ is convex. Furthermore, the $\varphi$-divergence 
generated by any convex function $\varphi:\bbR \to \bbR$ is convex. 
Finally, for any nonnegative cost function $c$, the optimal transport cost $\calT_c$ is convex.
\end{lemma}
Though our main focus is on the above divergences, our results also apply to generic convex
functionals 
$\Phi:\calP(\calX) \to \bbR$, 
which satisfy $\Phi(\lambda \mu_1 + (1-\lambda) \mu_2)
 \leq \lambda \Phi(\mu_1) + (1-\lambda) \Phi(\mu_2)$
for all $\lambda \in [0,1]$ and $\mu_1,\mu_2 \in \calP(\calX)$. 
Notice that for any fixed distribution
$Q \in \calP(\calX)$, the map $\Phi(P) = D(P\|Q)$ is a convex 
functional---additional examples include the negative differential entropy
(cf. Section~\ref{sec:smooth})
and certain expectation functionals~(cf. Section~\ref{sec:mmd}).  

Maximal martingales inequalities
form a key tool 
in the development of confidence
sequences, 
thus we provide an overview
in what follows.
Let $(\Omega, \calF, \bbP)$ be a probability space, 
over which all processes hereafter will be taken. Before discussing time-reversed concepts, it is useful to first  briefly overview standard martingales and filtrations.

\subsection{Forward filtrations, martingales and maximal inequalities}
\label{sec:martingales}

A forward filtration is a sequence
of $\sigma$-algebras $(\calF_t)_{t=1}^\infty$ 
contained in $\calF$ 
which is nondecreasing:
\[
\calF_t \subseteq \calF_{t+1}, \quad
\text{ for all $t\geq 1$. }
\]
Any process $(S_t)_{t=1}^\infty$
on $\Omega$ is adapted to its canonical (forward) filtration
$(\calC_t)_{t=1}^\infty$ defined by $\calC_t = \sigma(S_1, \dots, S_t)$, for all $t \geq 1$. A (forward) martingale
with respect to a (forward) filtration $(\calF_t)_{t=1}^\infty$
is a stochastic process   $(S_t)_{t=1}^\infty$ 
such that for all $t \geq 1$, $S_t$
is $\bbP$-integrable,
$\calF_t$-measurable, and satisfies 
$$\bbE[S_{t+1} | \calF_t] = S_t, \quad \text{ for all }  t\geq 1.$$
Supermartingales and submartingales
are respectively defined by replacing the equality in the above
display by $\leq$ and $\geq$. 
When it causes no confusion, 
we frequently abbreviate 
$(S_t)_{t=1}^\infty$ and 
$(\calF_t)_{t=1}^\infty$ by $(S_t)$ and $(\calF_t)$.
The construction of confidence sequences
typically relies on maximal inequalities. 
A prominent example is Ville's inequality
\citep{ville1939}, 
which states that any nonnegative supermartingale
$(S_t)$ satisfies
\begin{equation} 
\label{eq:forward_ville}
\bbP(\exists t \geq t_0: S_t \geq u) \leq \frac{\bbE[S_{t_0}]}{u}, \quad \text{ for all $u>0$ and all integers
$t_0 \geq 1$.}
\end{equation}
 Inequality \eqref{eq:forward_ville}
is a time-uniform extension of Markov's inequality
for nonnegative random variables. 
The unbounded range in Ville's inequality is made possible
by nonnegativity, and the fact that supermartingales have nonincreasing
expectations---see for instance~\citet{howard2020} for a formal proof. In contrast, submartingales
admit nondecreasing expectations, 
and do not generally satisfy an infinite-horizon 
inequality such
as \eqref{eq:forward_ville}. 
Nonnegative submartingales $(S_t)_{t=1}^\infty$
instead satisfy Doob's submartingale inequality
(e.g. \cite{durrett2019}, Theorem 4.4.2): 
\begin{equation}
\label{eq:forward_doob}
\bbP(\exists t \leq T: S_t \geq u) \leq \frac{\bbE[S_T]}{u} \quad \text{ for all $u > 0$ and any integer $T \geq 1$.}
\end{equation}
The prototypical example of a martingale
is a sum $S_t = \sum_{i=1}^t X_i$ of i.i.d. random
variables $(X_t)_{t=1}^\infty
\subseteq \calX \subseteq \bbR^d$, with respect
to its canonical filtration 
$\calC_t=\sigma(X_1, \dots, X_t)$. 
As described in
Section~\ref{sec:background_confseq}, 
a wealth of existing works have developed
confidence sequences for the expected
value of a sequence of i.i.d.
random variables, on the basis of inequalities
such as $\eqref{eq:forward_ville}$ and $ \eqref{eq:forward_doob}$.

\subsection{Reverse filtrations, martingales and maximal inequalities}
\label{sec:martingales2}
An alternate approach is to apply maximal inequalities
to the sample mean itself, namely to the process
$R_t = (1/t)\sum_{i=1}^t X_i$.
 This approach is relatively underexplored but turns out to be well-suited to our goals.
Unlike $(S_t)$, the process
$(R_t)$ is a
\textit{reverse} martingale.
To elaborate, a reverse filtration is 
a \textit{nonincreasing} sequence of $\sigma$-algebras 
$(\calF_t)_{t=1}^\infty$ contained in $\calF$:
\[
\calF_{t} \supseteq \calF_{t+1} \quad \text{ for all $t \geq 1$.}
\]
A $\bbP$-integrable process $(R_t)_{t=1}^\infty$ is then
said to be a reverse martingale with respect
to $(\calF_t)_{t=1}^\infty$
if for all $ t \geq 1$,
$R_t$ is $\calF_t$-measurable and
$$\bbE[R_t|\calF_{t+1}] = R_{t+1}, \quad \text{ for all }t \geq 1.$$
Reverse submartingales and supermartingales are
defined analogously, with the above equality
respectively replaced by $\geq$ and $\leq$. 
The sample average $R_t = ( 1 /t )\sum_{i=1}^t X_i$
defines a reverse martingale with respect
to its canonical reverse
filtration 
$\big(\sigma(R_{t}, R_{t+1}, \dots)\big)_{t=1}^\infty$.
$(R_t)$ is also a reverse
martingale with respect to
a richer filtration known as the
{\em symmetric} or 
{\em exchangeable filtration}
\citep{pollard2002, durrett2019}.
\begin{definition}[Exchangeable Filtration]
\label{def:exchangeable_filtration}
Given a sequence of random variables
$(X_t)_{t=1}^\infty$, the exchangeable
filtration is the reverse filtration
$(\calE_t)_{t=1}^\infty$, 
%where $\calE_0 = \sigma(X_1, X_2, \dots)$, 
%and 
% for all $t \geq 1$,
where $\calE_t$ denotes the $\sigma$-algebra
generated by all 
real-valued Borel-measurable functions
of $X_1, X_2, \dots$ which are permutation-symmetric in their 
first $t$ arguments.
% Keep this comment...
% Equivalently,
% $$\calE_t = \sigma\left(\left\{ 
% \sum_{i=1}^t I(X_i\in B)
% : B \in \bbB(\calX)\right\}
% \cup\{ X_{t+1}, X_{t+2}, \dots\}\right),\quad t \geq 1.$$
\end{definition} 
Informally, one may think of $\calE_t$ as $\sigma(P_t,X_{t+1},X_{t+2},\dots)$ or $\sigma(P_t,P_{t+1},P_{t+2},\dots)$, where $P_t = (1/t) \sum_{i=1}^t \delta_{X_i}$ is the empirical measure---while we do not formally define such $\sigma$-algebras, 
we present them here because we have found this intuition to be useful.
Forward filtrations are reflective of 
the worlds that humans inhabit, in which 
information increases with time, but the exchangeable filtration can be an odd object to initially gain intuition about, so we provide another informal analogy for the unfamiliar reader. One can imagine $\calE_t$ to be the information accessible to an ``amnesic oracle'', who at any time $t$ knows the entire future, but is confused about the ordering of events in the past. The amnesia causes the oracle to know less and less as time passes, because events in the initially clear future occur, and get transformed to a muddled~past.

Given a sequence $(X_t)_{t=1}^\infty$ of 
exchangeable random variables, recall that $P_t = (1/t) \sum_{i=1}^t \delta_{X_i}$ 
denotes their empirical measure.
It can be seen that the process $(P_t(B))_{t\geq 1}$ is a reverse martingale
with respect to the exchangeable filtration $(\calE_t)_{t=1}^\infty$,
for any fixed Borel set $B \subseteq \bbR^d$.
This property makes $(P_t)_{t=1}^\infty$ into a 
so-called 
measure-valued reverse martingale \citep{kallenberg2006}. 
% role in this paper because the empirical
% measure $P_t = (1/t) \sum_{i=1}^t \delta_{X_i}$
% forms a measure-valued reverse martingale
% with respect to $(\calE_t)$, in the sense
In fact, the converse nearly holds true:
if $(P_t)$ is a measure-valued reverse martingale,
and if $(X_t)$ is stationary, then $(X_t)$ 
is exchangeable \citep{bladt2019}.
% a more general statement is true, and we highlight this pivotal fact below.
% \begin{fact}[\cite{kallenberg2006}, Theorem 2.4]
% The following two statements are equivalent:
% \begin{enumerate} 
%     \item $X_1,X_2,\dots$ is an exchangeable sequence of random variables.
%     \item $(P_t)$ is a measure-valued reverse martingale with respect to  $(\calE_t)$.
% \end{enumerate}
% \end{fact}
 Note that the exchangeability condition 
 is weaker  than the i.i.d.\  assumption
 which we shall assume for the majority of this paper. We will later see that, with some care, this measure-valued reverse martingale gets effectively translated into a (real-valued) reverse submartingale property for convex functionals.

A key technical tool for handling
real-valued, 
 reverse submartingales will be the following analogue
of Ville's inequality \eqref{eq:forward_ville}, 
first proved by~\citeauthor{doob1940} (\citeyear{doob1940}; Theorem 1.1, p. 458)
for reverse martingales;
see also~\citeauthor{lee1990} (\citeyear{lee1990}; Theorem 3, p.~112). 
\begin{theorem}[Ville's Inequality for Nonnegative
Reverse Submartingales]
\label{thm:reverse_ville}
Let $(R_t)_{t=1}^\infty$ be a nonnegative reverse submartingale
with respect to a reverse
filtration $(\calF_t)_{t=1}^\infty$.
Then, for any integer $t_0 \geq 1$ and real number $u > 0$,
$$\bbP\big( \exists t \geq t_0: R_t \geq u\big) \leq \frac{\bbE[R_{t_0}]}{u}.$$
 \end{theorem}
Since Theorem~\ref{thm:reverse_ville}
plays a central role in our development, we  provide two
self-contained proofs in Appendix 
\ref{app:proofs_sec2} for completeness. As an example, 
it can be deduced from Theorem~\ref{thm:reverse_ville}
that for any sequence of exchangeable random variables $(X_t)_{t=1}^\infty$,
\begin{equation}
\bbP\left( \sup_{t \geq 1} \left| \frac 1 t \sum_{i=1}^t X_i\right| \geq 
u\right) \leq \frac{\bbE |X_1|}{u},
\end{equation}
for all $u > 0$.
This can be seen as a strengthening of Markov's inequality, whose left-hand side
is $\bbP(|X_1| \geq u)$. 

\subsection{\smash{Partially ordered martingales}}
\label{sec:background_partial_order}
In order to handle the two-sample process $(M_{ts})$ in equation 
\eqref{eq:process_Mts}, 
we also employ (reverse) martingales indexed
by $\bbN^2$. We endow $\bbN^2$ with the standard partial
ordering, that is we write 
$(t,s) \leq (t',s')$ for all $(t,s),(t',s')\in \bbN^2$ 
such that $t \leq t'$ and $s \leq s'$. 
The notation $(t,s) < (t',s')$ indicates that at least
one of the componentwise inequalities holds strictly.
The definitions of filtration and martingale extend to this
setting in the natural way: a forward (or reverse) filtration
is a sequence $(\calF_{ts})_{t,s \geq 1}$ of $\sigma$-algebras
which is nondecreasing (or nonincreasing) with respect
to the partial ordering, and a forward (or reverse) 
martingale $(S_{ts})_{t,s \geq 1}$ (or $(M_{ts})_{t,s \geq 1}$)
is an $L^1(\bbP)$ process adapted to $(\calF_{ts})$ which satisfies
$\bbE[S_{ts} | \calF_{t's'}] = S_{t's'}$ for all 
$(t,s) > (t',s')$ (or $\bbE[R_{ts}|\calF_{t's'}] = R_{t's'}$
for all $(t,s) < (t',s')$). When the latter equalities
are replaced by the inequalities $\leq$ and $\geq$, 
$(M_{ts})$ is respectively called a (reverse) 
supermartingale and submartingale.
We refer to \cite{ivanoff1999} for a survey.

 \cite{cairoli1970} showed that a direct analogue of Ville's inequality
cannot hold for partially ordered nonnegative martingales---the inequality $\bbP(\exists t,s \geq 1: S_{ts}\geq u) \leq \bbE[S_{11}]/u$
does not generally hold for all $u > 0$. Distinct maximal inequalities
for partially ordered forward
and reverse martingales have nevertheless been established 
by \cite{christofides1990} 
under suitable
moment assumptions, and under the so-called
conditional independence (CI) assumption introduced by 
\cite{cairoli1975}. A reverse filtration $(\calF_{ts})_{t,s\geq 1}$
is said to satisfy the CI property 
if for all $t,t',s,s' \geq 1$, 
\begin{equation} 
\label{eq:conditional_independence}
\bbE \big\{ \bbE[ ~\cdot ~|~ \calF_{ts'}] ~\big|~\calF_{t's}\big\} = 
\bbE\big\{~\cdot ~|~\calF_{(t\vee t')(s\vee s')}\big\},
\end{equation}
or equivalently, that $\calF_{ts'}$ and $\calF_{t's}$ are conditionally
independent given $\calF_{(t\vee t')(s\vee s')}$
\citep{merzbach2003}. 
The following Ville-type inequality for 
partially ordered reverse submartingales
can be obtained by employing 
Corollary~2.9 of \cite{christofides1990}.
\begin{proposition}
\label{prop:maximal_inequality_moments}
Let $(R_{ts})_{t,s \geq 1}$ be a nonnegative
reverse submartingale 
with respect to a reverse filtration $(\calF_{ts})_{t,s\geq 1}$
satisfying the conditional independence assumption. 
Assume that for some $\alpha > 1$,
$R_{ts} \in L^\alpha(\bbP)$ for 
all $t,s \geq 1$. Then, for all $u > 0$, 
$$\bbP(\exists t \geq t_0, s \geq s_0: 
R_{ts} \geq u) \leq 
\left(\frac \alpha{\alpha-1}\right)^\alpha 
\frac{\bbE[R_{t_0s_0}^\alpha]}{u^\alpha}.$$
\end{proposition}
The special case $\alpha=2$ was stated in Corollary~2.10
of \cite{christofides1990}. A proof
and further discussion of this result is given in Appendix~\ref{app:proof_prop_max_inequality_moments}, 
and forms the basis for our two-sample results.

\subsection{Time-uniform confidence sequences} 
\label{sec:background_confseq}
Confidence sequences
are defined similarly as in~\eqref{eq:confseq} for estimands
other than divergences. 
Given a functional $\theta\equiv\theta(P)$ of interest, 
and an error level $\delta \in (0,1)$,
a $(1-\delta)$-confidence sequence 
$(C_t)_{t=1}^\infty$ based on an i.i.d.
sequence of random
variables $(X_t)_{t=1}^\infty$ from $P$
is a sequence of sets $C_t \in \sigma(X_1, \dots, X_t)$
satisfying
$\bbP(\exists t \geq 1 : \theta \not\in C_t) \leq \delta.$
When $\theta$
is real-valued, we say two sequences
$(\ell_t)$ and $(u_t)$ 
are  lower and upper
confidence sequences
if 
$\bbP(\exists t \geq 1: \theta \leq \ell_t)
\leq \delta$ and $
\bbP(\exists t \geq 1 : \theta \geq u_t)\leq \delta$
respectively.
Confidence sequences
% are a central methodology for 
% statistical inference
% in settings where
% the sample size is not fixed in advance. 
% Early work was 
were pioneered by Robbins, Darling, Siegmund and Lai
\citep{darling1967, 
robbins1969, lai1976}, and  
new techniques have been recently developed that enable their 
extensions to
new, nonparametric settings 
\citep{kaufmann2018,howard2021}. 
This resurgence of interest in sequential
analysis has been driven in part
by its applications
to best-arm identification algorithms
for multi-armed
bandits \citep{jamieson2014,kaufmann2016,shin2019a} and reinforcement learning \citep{karampatziakis2021}, to name a few. 

The mean of a distribution is perhaps
the target of inference which has received the most attention
in prior work on confidence sequences. As described in Section 
\ref{sec:martingales}, 
the process $S_t = \sum_{i=1}^t (X_i-\mu)$
forms a canonical example
of a forward martingale, where $\mu$ denotes
the common (finite) mean of the i.i.d. random variables $X_i$. 
To obtain a confidence sequence for
$\mu$, a maximal inequality
such as that of Ville (see equation \eqref{eq:forward_ville})
cannot be directly be applied to $(S_t)$, however, 
since it is not a nonnegative process. 
Inspired by the Cram\'er-Chernoff
method for deriving concentration inequalities, 
it is instead natural
to consider the nonnegative process
\[
U_t(\lambda) = \exp(\lambda S_t),\quad t\geq 1,
\] for $\lambda > 0$.
Here, we assume a tail assumption is placed 
on $X_1$, such that it admits a finite cumulant generating
function satisfying 
$\log\{\bbE[\exp(\lambda (X_1-\mu))]\} \leq 
\phi(\lambda)$, for some (say, known) 
map $\phi:[0,\mlambda) \to \bbR$, where $\mlambda > 0$.
The exponential process  $(U_t(\lambda))$
is a nonnegative submartingale by Jensen's inequality, 
and forms the basis of several confidence
sequences described below. 
A distinct line of 
work constructs confidence sequences for 
$\mu$ by downweighting this process to recover
a (super)martingale. 
For instance, it can be verified
that
\[
L_t(\lambda)=\exp\{\lambda S_t - t \phi(\lambda)\},
\quad t  \geq 1,
\]
is a nonnegative supermartingale with respect
to the canonical filtration. 
Variants of the process $(L_t(\lambda))$ appear in a long
line of work aimed at deriving
sequential concentration inequalities for means---see 
\cite{howard2020} for a comprehensive review of such approaches.

Applying a maximal martingale inequality
to either $(U_t(\lambda))$
or $(L_t(\lambda))$ does not, on its own, lead  
to satisfactory confidence sequences for $\mu$. 
Infinite-horizon maximal inequalities are not
 available for submartingales $(U_t(\lambda))$, 
and even for the supermartingale $(L_t(\lambda))$, 
a direct application of Ville's inequality leads to a confidence
sequence for $\mu$ with nonvanishing length.
To obtain
confidence sequences with lengths scaling at
the optimal rate 
$O(\sqrt{\log\log t/t})$,
implied by the law of the iterated logarithm, 
it is instead common to use a variant of the ``method of mixtures'', for example by repeatedly applying
a maximal inequality over
geometrically-spaced epochs in time---such methods
are often known as peeling, chaining, or stitching.
\cite{jamieson2014,zhao2016} use stitching arguments based
on $(U_t(\lambda))$ and Doob's submartingale inequality, 
while \cite{garivier2013,kaufmann2016, howard2021},
use $(L_t(\lambda))$ and 
Ville's inequality. The resulting
confidence sequences decay at similar rates, though
with varying constants and tail assumptions---see \cite{howard2021,waudby-smith2021}
for a comparison of such approaches.

In Appendix \ref{sec:forward_submart},
we describe the difficulties in extending 
the above framework to targets
of inference distinct from means. We provide therein
a naive approach to deriving confidence sequences
for some
divergences between probability distributions, using
an analogue of the process $(U_t(\lambda))$.
However, our main results
which follow provide a far-reaching 
improvement, 
which hinges
upon {\em reverse} submartingales---a 
rarely used tool for deriving confidence sequences.

\section{\smash{Confidence sequences
for convex functionals}}
\label{sec:main}

Let $(X_t)_{t=1}^\infty$
and $(Y_s)_{s=1}^\infty$ respectively denote 
independent sequences of i.i.d. 
observations from two distributions
$P,Q \in \calP(\calX)$,
with support contained in a set $\calX \subseteq \bbR^d$. 
% The primary goal of this section is to derive confidence
% sequences for the divergence $D(P\|Q)$. 
% In the interest of generality, we shall formulate our
% results as confidence sequences for general  
Given convex functionals $\Phi: \calP(\calX) \to \widebar\bbR$ and 
$\Psi: \calP(\calX) \times \calP(\calX) \to \widebar\bbR$,
the goal of this section is to derive confidence sequences for
$\Phi(P)$ and $\Psi(P,Q)$. 
Our primary interest is in the special case where $\Psi$ is a convex divergence $D$, 
and $\Phi = D(\cdot\|Q)$ when $Q$ is known, but we formulate our results for arbitrary convex  functionals in the interest of generality. 
We shall make use of the processes
\begin{equation}\label{eq:one-sample-divergence-Nt}
N_t = \Phi(P_t) - \Phi(P), \quad 
M_{ts} = \Psi(P_t, Q_s) - \Psi(P, Q), \quad t,s \geq 1.
\end{equation}
We prove in Section~\ref{sec:lcb},
that $(M_{ts})_{t,s \geq 1}$ and $(N_t)_{t \geq 1}$ are reverse
submartingales with respect
to suitable filtrations, whose
choice is further discussed
in Section \ref{sec:choice_filtrations}.
We then derive maximal inequalities
for these processes, from which lower confidence sequences
will follow using an epoch-based analysis. 
We follow a distinct strategy to obtain
upper confidence sequences in Section~\ref{sec:ucb}.

\subsection{\smash{Lower confidence sequences via 
reverse submartingales}}
\label{sec:lcb}
Let $(\calE_t^X)_{t=1}^\infty$ and
$(\calE_s^Y)_{s=1}^\infty$ denote 
the exchangeable
filtrations associated with the sequences
$(X_t)_{t=1}^\infty$ and $(Y_s)_{s=1}^\infty$
respectively, and define
$\calE_{ts} = \calE_t^X \bigvee \calE_s^Y$ 
for all $t,s \geq 1$. 
As stated
in Section \ref{sec:martingales},
the sequences of empirical measures
$(P_t)$ and $(Q_s)$ form measure-valued reverse martingales
with respect to $(\calE_t^X)$ and $(\calE_s^Y)$
respectively.
This fact suggests that any convex functional evaluated at 
$(P_t,Q_s)$ is a reverse submartingale, as we now show. 
\begin{theorem} 
\label{thm:submart} 
Let  $\Phi,\Psi$
be convex
functionals such that $\Phi(P_t), \Psi(P_t,Q_s) \in L^1(\bbP)$
for all $t,s \geq 1$. Then, 
\begin{enumerate}
\item[(i)]
% (One-Sample.)  
 $(\Phi(P_t))_{t\geq 1}$ is a reverse submartingale with respect to $(\calE_t^X)$.
\item[(ii)]
% (Two-Sample.) 
 $(\Psi(P_t,Q_s))_{t,s \geq 1}$ is a 
partially ordered reverse submartingale with respect to $(\calE_{ts})$.
\end{enumerate}
% In particular, $(N_t)$ in~\eqref{eq:one-sample-divergence-Nt} is a 
% reverse submartingale
% with respect to $(\calE_t^X)$, while $(M_{ts})$ in~\eqref{eq:process_Mts} is a  reverse submartingale with respect to $(\calE_{ts})$, whenever these processes are in $L^1(\bbP)$.
\end{theorem} 
Theorem 
\ref{thm:submart}(i) is known 
in the special case where $\Phi$ is the supremum
of an empirical processes---see for instance \cite{pollard1981a}
% , that is, for the
% one-sample process $(N_t)$
% when $D$ is taken to be an IPM and $P=Q$
and Lemma 2.4.5 of \cite{vandervaart1996}.
Our proof below is inspired by these works, and in particular extends them to general
convex functionals and to partially ordered filtrations. 

Notice that the processes in Theorem~\ref{thm:submart}
lie in $L^1(\bbP)$, and are therefore assumed to be measurable. 
% There are well-known instances when this assumption may be violated, such as when $D$ is an IPM over an uncountable collection of functions. 
When this assumption is removed, 
it can be shown that there exist
measurable \emph{covers} of the processes in Theorem~\ref{thm:submart} 
for which the result continues to hold---this is the approach
taken by \cite{vandervaart1996} for suprema
of empirical processes; see also \cite{pollard1981a} and \cite{strobl1995}. With this modification, we believe
that all subsequent claims in this paper 
can be made to hold in outer probability, but we prefer
to retain the measurability condition to avoid being overwhelmed
by technicalities. In particular, $(N_t)$ and $(M_{ts})$
are tacitly assumed to be measurable throughout the sequel. 

\textbf{Proof of Theorem~\ref{thm:submart}.}
We will prove Theorem~\ref{thm:submart}(ii)
and a similar argument can be used
to prove Theorem~\ref{thm:submart}(i).
For all $t,s \geq 1$ $\Psi(P_t, Q_s)$ is 
invariant to permutations of 
$X_1, \dots, X_t$ and of $Y_1, \dots, Y_s$. 
It follows that
$(\Psi(P_t, Q_s))$ is adapted to $(\calE_{ts})$,
thus it suffices to prove the 
reverse submartingale property.
Fix $t,s \geq 1$. Define the $(t+1)$ different ``leave-one-out'' empirical measures 
$$P_t^i = \frac 1 t\left[ \sum_{j=1}^{i-1} \delta_{X_j} +  \sum_{j=i+1}^{t+1} \delta_{X_j}\right], \quad i=1, 2, \dots, t+1.$$
Then
$P_{t+1} = \frac 1 {t+1} \sum_{i=1}^{t+1} P_t^{i}
%, \quad
%  Q_{t+1} = \frac 1 {t+1} \sum_{i=1}^{t+1} Q_t^i.
$, and the convexity of $\Psi$ implies
$$\Psi(P_{t+1}, Q_{s})\leq \frac 1 {t+1}\sum_{i=1}^{t+1} \Psi(P_t^i ,  Q_s).$$
Since $\Psi(P_{t+1}, Q_s)$ is $\calE_{(t+1)s}$-measurable,
%and the sequences $(X_t)$ and $(Y_s)$ are mutually
%independent, 
we deduce that
\begin{equation} 
\label{eq:pf_reverse_submart_convex}
\Psi(P_{t+1} ,  Q_{s})\leq
\frac 1 {t+1}\sum_{i=1}^{t+1} \bbE[\Psi(P_t^i ,  Q_s)|\calE_{(t+1)s}].
\end{equation}
Since $P_t^{t+1}=P_t$ by definition, the claim will 
follow upon proving the key identity
\begin{equation}
\label{eq:pf_reverse_submart_key_eq}
\bbE[ \Psi(P_t^i ,  Q_s) | \calE_{(t+1)s}] = \bbE[\Psi(P_t, Q_s)|\calE_{(t+1)s}],\quad i=1, \dots, t.
\end{equation}
Notice that 
% since $(X_t)$ and $(Y_s)$ are independent sequences, 
we can write
$\calE_{(t+1)s} = \calE_{t+1}^X \bigvee \calE^Y_s
 = \sigma(\calI)$ where 
 $$\calI =  \{A^X \cap A^Y: A^X \in \calE_{t+1}^X,
 A^Y \in \calE_s^Y\}.$$
$\calI$ is clearly a $\pi$-system. 
To prove
\eqref{eq:pf_reverse_submart_key_eq},
it will thus suffice to prove that for any 
set $A=A^X \cap A^Y$, with $A^X \in \calE_{t+1}^X$
and $A^Y \in \calE_s^Y$, and for all
$1 \leq i \leq t$, we have 
$\bbE[ \Psi(P_t^i ,  Q_s) I_A] = \bbE[\Psi(P_t, Q_s)I_A],$
where $I_A : \Omega \to \{0,1\}$ is the indicator
function of $A$ and $I_A = I_{A^X} I_{A^Y}$. 
% Here, for all $\omega \in \Omega$,
% $I_A(\omega)$ takes the value 1 if
% $\omega \in A$ and 0 otherwise. 
%%Fix $A \in \calE_{t+1}^X$ and $1 \leq i \leq t$. 

Notice that 
$I_{A^X}$ is  $\calE^X_{t+1}$-measurable, 
thus it is a function $f_{A^X}$ of $X_1, X_2, \dots$
which is permutation symmetric in $X_1, \dots, X_{t+1}$. For convenience, we write $I_{A^X}=f_{A^X}(X_1, X_2, \dots)$, 
whence
we may also write 
$I_A = I_{A^Y} f_{A^X}(X_1, X_2, \dots)$.
Now, let $\tau:\bbN \to \bbN$
be the permutation such that 
$\tau(j) = j$ if $j\not\in\{t+1,i\}$
and $\tau(t+1) = i$, $\tau(i) = t+1$. By exchangeability of $X_1, X_2, \dots$, and by~their independence from 
$Y_1, Y_2, \dots$, 
we have
$(X_1, Y_1, X_2, Y_2, \dots)
\overset{d}{=} (X_{\tau(1)}, Y_1, X_{\tau(2)}, Y_2, \dots)$,~whence,
$$\bbE[\Psi(P_t, Q_s)I_{A}]
= \bbE[\Psi(P^{t+1}_t, Q_s)f_{A^X}(X_1, X_2, \dots)I_{A^Y}]
 = \bbE[\Psi(P_t^i, Q_s)f_{A^X}(X_{\tau(1)}, X_{\tau(2)}, \dots)I_{A^Y}].$$
% where we have again used the fact that $P_t = P_t^{t+1}$ in the first equality. 
Since $f_{A^X}$ is permutation-symmetric in its first $t+1$ arguments, and the permutation $\tau$ fixes
all natural numbers greater or equal to $t+2$, we obtain
$$f_{A^X}(X_{\tau(1)}, X_{\tau(2)}, \dots, ) = 
 f_{A^X}(X_1, X_2, \dots),$$
 implying that 
$$\bbE[\Psi(P_t, Q_s)f_{A^X}(X_1, X_2, \dots) I_{A^Y} ]
 = \bbE[\Psi(P_t^i, Q_s)f_{A^X}(X_1, X_2, \dots)I_{A^Y}].$$
% The definition of $(\calE_t^X)$ together with the fact
% that each of $D(P_t^i\|Q)$ is permutation symmetric
% in $X_1, \dots X_{i-1}, X_{i+1}, \dots, X_{t+1}$, for all $i=1, \dots, t+1$,
% implies that 
% $\bbE[ D(P_t^i \| Q_s) | \calE_{t+1}^X] = \bbE[D(P_t\|Q_s)|\calE_{t+1}^X].$
Equation \eqref{eq:pf_reverse_submart_key_eq}
now follows. 
Returning to equation 
\eqref{eq:pf_reverse_submart_convex} we deduce that
$$\Psi(P_{t+1} ,  Q_s) %\leq \frac 1 {t+1}\sum_{i=1}^{t+1} D_(P_t^i \| Q_t^i)
% = \bbE \big[ D(P_{t+1} \| Q_{s}) \big| \calF_{t+1}^X\big]
% \leq \frac 1{t+1} \sum_{i=1}^{t+1}
%\bbE\left[D(P_t^i \| Q_t^i) | \calF_{t+1}\right]
 \leq   \bbE[\Psi(P_t ,  Q_s) | \calE_{(t+1)s}].$$
A symmetric argument shows that
$\Psi(P_{t} ,  Q_{s+1})
 \leq   \bbE[\Psi(P_t ,  Q_s) | \calE_{t(s+1)}]$, 
implying that $(\Psi(P_t, Q_s))$ is a
partially ordered reverse submartingale with respect to $(\calE_{ts})$.
\qed

It is apparent from the proof that the convexity requirement 
is stronger than necessary. For example, 
the following more general statement can be inferred  from
the proof of Theorem~\ref{thm:submart}. We record it formally
as it may be of independent interest. 
\begin{proposition}
\label{prop:loo}
Suppose $(R_t)_{t=1}^\infty$  is an $L^1(\bbP)$ process of the form
$R_t = f_t(X_1, \dots, X_t)$, 
for some sequence of permutation invariant 
maps $f_t:\calX^t \to \bbR$. Suppose further that for all $t \geq 1$, 
\begin{equation} 
\label{eq:loo} 
R_{t+1} \leq \frac 1 {t+1} \sum_{i=1}^{t+1}
f_{t}(X_1, \dots, X_{i-1}, X_{i+1}, \dots, X_{t+1}).
\end{equation}
Then, $(R_t)$ is a reverse submartingale with respect
to $(\calE_t^X)$. 
Furthermore, if the above display holds with equality, 
$(R_t)$ is a reverse martingale with respect to $(\calE_t^X)$. 
\end{proposition}
We refer to condition~\eqref{eq:loo} as the {\it leave-one-out property}. 
In Sections~\ref{sec:mmd} and~\ref{sec:rademacher}
we will briefly make use of 
 processes which cannot be expressed as evaluations of a convex
 functional at the empirical measure, 
 but which
 satisfy the leave-one-out property.

Theorem~\ref{thm:submart} 
permits the use of 
maximal inequalities discussed in 
Sections~\ref{sec:martingales2} and 
\ref{sec:background_partial_order}
for reverse submartingales. In view of generalizing the Cram\'er-Chernoff
technique 
\citep{boucheron2013} to our sequential setup, we shall state our bounds
under the assumption that $N_t$ and $M_{ts}$ admit finite 
cumulant generating functions
over an interval $[0,\mlambda)$, and are upper bounded 
by known convex functions 
$\psi_{ts},\psi_t:[0,\mlambda) \to \bbR$, 
\begin{equation} 
\label{eq:tail_bound_main}
\log\Big\{\bbE\big[\exp(\lambda M_{ts})\big] \Big\} \leq \psi_{ts}(\lambda), \quad
\log\Big\{\bbE\big[\exp(\lambda N_t )\big] \Big\} \leq \psi_{t}(\lambda),
\quad t,s \geq 1, ~\lambda \in [0,\mlambda). 
\end{equation}
% Notice that $\psi_t = \psi_{t0}$. 
% We also define the Cram\'er 
% transform of $M_{ts}$ by
% $\psi_{ts}^*(u) = \sup_{\lambda \in [0,\mlambda)} \{\lambda u - \psi_{ts}(\lambda)\}$,
% and similarly for $N_t$.
We shall discuss in Section \ref{sec:applications}
how such tail assumptions may be replaced
by tail assumptions on the distributions
$P$ and $Q$ themselves, for various
special cases of functionals.
Under equation \eqref{eq:tail_bound_main}, 
we obtain the following result.
\begin{proposition}
\label{prop:maximal_exponential_bounds}
Let $\Phi,\Psi$ be convex functionals such that 
the processes $(M_{ts})$ and $(N_t)$ satisfy
the bounds of equation \eqref{eq:tail_bound_main}.
Then, for all $u > 0$,
and all integers $t_0,s_0 \geq 1$, the following hold. 
\begin{enumerate}
\item[(i)] (One-Sample) $\bbP(\exists t \geq t_0: N_t \geq u) \leq \exp(-\psi_{t_0}^*(u)).$
\item[(ii)] (Two-Sample) $\bbP(\exists t \geq t_0,s \geq s_0: M_{ts} \geq u) 
\leq e \cdot \exp(-\psi_{t_0s_0}^*(u)).$
\end{enumerate} 
\end{proposition}
A proof of Proposition 
\ref{prop:maximal_exponential_bounds}
appears in Appendix \ref{app:proofs_sec4}.
In view of Theorem~\ref{thm:submart}, 
Proposition~\ref{prop:maximal_exponential_bounds}(i) 
is obtained through
an application of the Cram\'er-Chernoff technique, 
together with Ville's inequality for reverse submartingales (Theorem~\ref{thm:reverse_ville}). It can thus also be seen as an extension of the supermartingale techniques in \citet{howard2019} to the reversed setting.
% , which was itself a generalization of the Chernoff and martingale methods.
In Proposition~\ref{prop:maximal_exponential_bounds}(ii), an analogous bound is obtained 
for the two-sample case, though with
the additional factor $e$ in the probability bound. 
The presence of such a factor greater than 1 is necessary due to the aforementioned counterexample
of \cite{cairoli1970} regarding maximal inequalities for partially ordered martingales, though we do not know
if the value $e$ is sharp. The bound itself
is obtained via Proposition \ref{prop:maximal_inequality_moments}.

Inverting the probability inequalities of Proposition~\ref{prop:maximal_exponential_bounds} leads to one-
and two-sample confidence
sequences for $\Phi(P)$ or $\Psi(P,Q)$, though with lengths which
are constant with respect to $t,s$. To obtain confidence sequences scaling at 
rate-optimal lengths, we employ a stitching construction 
inspired by those described in
Section \ref{sec:background_confseq}, together with
Proposition~\ref{prop:maximal_exponential_bounds}. 
Our result will depend on user-specified functions
$\osfn,\tsfn: [0,\infty) \to [1,\infty)$, known as stitching functions,
which dictate the shape of the resulting confidence sequences below. 
We construct these functions to satisfy
\begin{equation} 
\label{eq:stitching_fn_main}
\sum_{k=1}^\infty \frac 1 {\osfn(k)} \leq 1, \quad
\sum_{j,k=1}^\infty \frac e {\tsfn(k+j)} \leq 1,
\end{equation}
as well as $\ell(j) = \ell(1)$ for all $j\in [0,1]$, and 
$g(k) = g(2)$ for all $k\in [0,2]$. Typical choices include 
$\ell(k) =  (1\vee k)^\alpha \zeta(\alpha)$
and
$g(k) = e(2\vee k)^{\alpha+1} (\zeta(\alpha) - \zeta(\alpha+1))
$
(\cite{borwein1987}, p.\ 305), 
where $\alpha > 1$ and $\zeta(\alpha) = \sum_{k=1}^\infty (1/k^\alpha)$. 
For ease of exposition, we  also insist
that the sequences
 $2^{-u}\log \ell(u)$ and $(2^{-u}+2^{-v})\log g(u+v)$
are chosen to be decreasing  in each of the indices $u,v \geq 1$, implying that $\ell(u), g(u) = o(\exp(\exp(u)))$.  
Our main result is below.
 \begin{theorem}
\label{thm:main}
Let $\Phi,\Psi$ denote convex functionals for which
the processes $(M_{ts})$ and $(N_t)$ satisfy
the bounds of equation \eqref{eq:tail_bound_main}.
For any integer $t \geq 1$, let
$\bar t = \lceil t/2\rceil$, and fix $\delta \in (0,1)$.
\begin{enumerate}
\item[(i)] (One-Sample) Assume $\psi_t^*$ is invertible for all $t\geq 1$, and that the  sequence
$$\gamma_t = 
(\psi_{\bar t}^*)^{-1}\Big(\log \osfn(\log_2  t) + \log(2/\delta)\Big), \quad t \geq 1$$
is nonincreasing.
Then,
\begin{align*}
\bbP&\Big\{\exists t \geq 1: \Phi(P_t) \geq  \Phi(P)+ \gamma_t \Big\} \leq \delta/2.
\end{align*}
\item[(ii)] (Two-Sample) Assume $\psi_{ts}^*$ is invertible for all $t,s\geq 1$, and that the sequence 
$$
 \gamma_{ts}  = 
(\psi_{\bar t \bar s}^*)^{-1}\Big(\log g(\log_2 t + \log_2 s) + \log(2/\delta)\Big),\quad t,s \geq 1$$
is nonincreasing with respect
to the partial order on $\bbN^2$.
Then,
\begin{align*}
\bbP&\Big\{\exists t,s \geq 1: \Psi(P_{t}, Q_s)\geq
  \Psi(P, Q) +
\gamma_{ts}\Big\}
\leq \delta/2.
\end{align*}
\end{enumerate}

\end{theorem}
We begin by noting that 
for a fixed sample size $n$, if we denote $\bar \gamma_{n} = (\psi_{n}^*)^{-1}\left( \log(2/\delta)\right)$, then
the fixed-time Cram\'er-Chernoff
concentration bound corresponding to
part (i) is given by
\citep{boucheron2013},
\[
\bbP\Big\{ \Phi(P_{n}) \geq  \Phi(P)+
\bar \gamma_{n}\Big\}
\leq \delta/2,
\]
and an analogous statement can also be made for part (ii) above. Thus, our time-uniform bounds are essentially an iterated logarithm factor worse than the usual fixed-time bounds, but now also apply at arbitrary stopping times. 

Before further commenting on the above result,
we instantiate it in the special 
case where $N_t$ and $M_{ts}$ are sub-Gaussian
for all $t,s\geq 1$. 
In Section \ref{sec:applications},
we illustrate how such a condition can be satisfied
under tail assumptions on the distributions
$P$ and $Q$ themselves.
% In what follows, given any $\sigma > 0$,
% a random variable $Y$ is said to be $\sigma^2$-sub-Gaussian
% if for all $\lambda \in \bbR$,
% $\bbE[\exp(\lambda Y)] \leq \exp(\lambda^2 \sigma^2/2).$
\begin{corollary}
\label{cor:sub_gaussian_bound}
Fix $\delta \in (0,1)$, and recall that $\bar t = \lceil t/2\rceil$.
\begin{enumerate}
\item[(i)] (One-Sample) Assume $N_t$ is  $\kappa_t^2$-sub-Gaussian
for some $\kappa_t > 0$, 
and for all $t \geq 1$. Choose $\ell$ so that
$(\kappa_{\bar t}^2\log \ell(\log_2 t))_{t \geq 1}$ is nonincreasing. 
Then,
$$\bbP \left\{ \exists t \geq 1 : N_t \geq \bbE\left(N_{\bar t}\right) + \sqrt{2 \kappa_{\bar t}^2 \Big[ \log \osfn(\log_2 t) + \log(2/\delta) \Big] }\right\} \leq \delta/2,$$
\item[(ii)] (Two-Sample) Assume $M_{ts}$ is $\sigma_{ts}^2$-sub-Gaussian for some $\sigma_{ts} > 0$, and for all $s,t \geq 1$.  Choose $g$ so that
$(\sigma_{\bar t\bar s}^2\log g(\log_2 t + \log_2 s))_{t,s \geq 1}$ is nonincreasing. 
Then, 
$$\bbP\left\{ \exists t,s \geq 1: M_{ts} \geq \bbE\left(M_{\bar t \bar s}\right) + \sqrt{2\sigma_{\bar t \bar s}^2\Big[ \log \tsfn(\log_2 t + \log_2 s) + \log(2/\delta)\Big]}\right\} \leq \delta/2.$$
\end{enumerate}
\end{corollary}
Theorem~\ref{thm:main}(i) 
is proved by dividing
time $t \geq 1$ into geometrically
increasing epochs of the form $[2^j, 2^{j+1}]$,
$j \geq 0$, over each of 
which we construct confidence boundaries at 
the level $\delta_j/2 = \delta/(2\osfn(j+1)) \in (0,1)$
using Proposition~\ref{prop:maximal_exponential_bounds}.
Taking a union bound over these boundaries
leads to a miscoverage probability of at most
$\sum_{j=0}^\infty (\delta_j/2) \leq \delta/2$. 
The two-sample process $(M_{ts})$ is handled similarly, 
by instead forming two sequences of epochs. 
In Appendix \ref{app:proofs_sec4}, 
we also state and prove
a more general version of Theorem~\ref{thm:main} 
in terms of epoch sizes different than 2, 
which will be needed in Section~\ref{sec:lil}.

{\bf Remark.} Given a convex divergence $D$, Theorem~\ref{thm:main}(i)
implies that $(1-\delta/2)$-upper confidence sequences
for  the processes $N_t^X = D(P_t\|P)$
and $N_s^Y = D(Q_s\|Q)$ are respectively given by
$$\gamma_t^X = (\psi_{X,\bar t}^*)^{-1}
\Big(\log \osfn(\log_2  t) + \log(2/\delta)\Big), ~~ 
\gamma_s^Y = (\psi_{Y,\bar s}^*)^{-1}
\Big(\log \osfn(\log_2 s) + \log(2/\delta)\Big),$$
where $\psi_{X,t}$ is an upper bound on the cumulant generating
function of $N_t^X$, and similarly for $\psi_{Y,s}$. 
% It follows that that
% the set $\{\mu \in \calP(\calX): D(P_t,\mu) \leq \gamma_t^X\}$ forms a 
% $(1-\delta/2)$-confidence sequence for the measure $P$, 
% and similarly for $Q$. 
When $D$ satisfies the triangle inequality, 
one may deduce the following two-sided confidence sequence for $D(P\|Q)$,
\begin{align}  
\label{eq:triangle_bound}
\nonumber 
\bbP\big(\forall t,s \geq 1: &|D(P_t\|Q_s)-D(P\|Q)| 
\leq \gamma_t^X + \gamma_s^Y\big) \\
 &\geq 1- 
      \bbP\big(\exists t \geq 1: N_t^X > \gamma_t^X\big) -
      \bbP\big(\exists  s \geq 1: N_s^Y > \gamma_s^Y\big)
 \geq 1-\delta.
 \end{align} 
Equation \eqref{eq:triangle_bound} is significant in that it provides a time-uniform bound for a partially ordered
reverse submartingale on the basis of 
two totally ordered reverse submartingales.
Doing so bypasses the nearly unavoidable factor $e$
in condition~\eqref{eq:stitching_fn_main}, 
but may 
nevertheless be looser in general due to the application of the triangle inequality. 
We also remark that equation~\eqref{eq:triangle_bound} does not require 
$P_t$ to be independent of $Q_s$, unlike Theorem~\ref{thm:main}(ii).

\subsection{Upper confidence sequences
via affine minorants}
\label{sec:ucb}
The lower confidence
sequences derived in Theorem~\ref{thm:main}
hinged upon the reverse submartingale property 
of the processes $(N_t)$ and $(M_{ts})$---an inherently one-sided
condition. We show in this section how a different
approach can be used to derive
upper confidence sequences, motivated
both by technical and statistical considerations. 
\begin{itemize} 
\item 
On the technical side, 
it would seem natural to repeat the 
steps of Theorem~\ref{thm:main}
with respect to 
the process $(-N_t)$ to obtain
an upper confidence sequence.
However, $(-N_t)$ is a reverse
\emph{super}martingale and thus cannot
satisfy
infinite-horizon (Ville-type) maximal
inequalities. Furthermore, 
the exponential process $(\exp(-N_t))$ may
generally be neither a reverse supermartingale nor
a submartingale, thus an analogue
of Proposition~\ref{prop:maximal_exponential_bounds}
cannot be derived. 
\item 
On the statistical side, 
the plug-in estimators
$\Phi(P_t)$ and $\Psi(P_t,Q_s)$ 
are typically upward biased
in estimating $\Phi(P)$ and $\Psi(P,Q)$ respectively. This fact can be deduced from
Corollary~\ref{cor:upward_bias} 
below, but can already be anticipated
from the fact that $\bbE [\Phi(P_t)]$
and $\bbE[ \Psi (P_t, Q_s)]$ are {\it nonincreasing}
sequences, since $(\Phi(P_t))$ and $(\Psi(P_t, Q_s))$
are reverse submartingales (Theorem~\ref{thm:submart}).
This upward bias suggests that
confidence sequences for $D(P\|Q)$ 
of the form $[D(P_t, Q) - \ell_t, D(P_t, Q) + u_t]$
should typically be asymmetric, 
with the sequence $(u_t)$ potentially decaying
at a faster rate than $(\ell_t)$.
% This indeed turns out the case, and motivates
% our distinct approach below.
 \end{itemize}

% We state our first result in terms of 
% two general convex functionals
% $\Phi: \calP(\calX) \to \widebar\bbR$ and 
% $\Psi: \calP(\calX) \times \calP(\calX) \to \widebar\bbR$. 
Our approach is summarized as follows.
The convexity of $\Phi$ guarantees
that it can be minorized by an affine functional 
on $\calP(\calX)$. 
Notice that an affine functional 
evaluated at the empirical measure is 
a sample average, and therefore a
reverse martingale. 
Furthermore, when  
a convex duality result guarantees that 
$\Phi$ is equal to the supremum over a set of minorizing affine functionals, 
it can be shown that the difference  
$\Phi(P_t) - \Phi(P)$
is in fact minorized by a \textit{mean-zero} 
sample average, for which confidence sequences
of length $O(\sqrt{\log\log t/t})$ can be obtained in a standard way
under appropriate tail conditions. Doing so leads
to a time-uniform lower bound on $\Phi(P_t)-\Phi(P)$,
which is easily rephrased as an upper confidence sequence for $\Phi(P)$ scaling at a near-parametric rate.

While this intuition can be made rigorous for a broad collection of convex 
functionals,
we shall
avoid doing so in full generality 
to avoid introducing additional terminology. 
We shall instead assume that $\Phi$ and $\Psi$ 
take the following form, which is sufficiently general
to cover the divergences of primary interest in our 
development,
\begin{equation} 
\label{eq:sup_affine}
\Phi(\mu) = \sup_{
\substack{f \in \calF_\Phi}}
\int f d\mu, \quad
\Psi(\mu,\nu)
 = \sup_{\substack{(f,g) \in \calH_\Psi}}
\int fd\mu + \int gd\nu,
 \end{equation}
for all $\mu,\nu \in \calP(\calX)$. 
Here, $\calF_\Phi$
denotes a set of Borel-measurable functions
on $\calX$, $\calH_\Psi$ 
a set of pairs of such functions,
% and $\calB_\Phi, \calB_\Psi \subseteq \bbR$.
and we also write $\calF_\Psi = \{f: (f,g) \in \calH_\Psi\}$
and $\calG_\Psi = \{g: (f,g) \in \calH_\Psi\}$.
It is clear from equations \eqref{eq:ipm},
\eqref{eq:kantorovich}
and
\eqref{eq:phi_div_variational} 
that if $D$ is an IPM, $\varphi$-divergence, or 
optimal transport cost, then 
the functionals $\Psi = D$ and
$\Phi = D(\cdot\|Q)$ (for a fixed
measure $Q \in \calP(\calX)$) 
admit the representation~\eqref{eq:sup_affine}---for instance, 
in the case of IPMs generated
by a function class $\calJ$, one may take
$\calF_\Phi = \{f - \int fdQ: f \in \calJ\}$ and
$\calH_\Psi = \{(f,-f): f \in \calJ\}$. 
With this notation in place, the following observation
is straightforward.
\begin{proposition} 
\label{prop:affine_rev_martingale}
If $\Phi$ and $\Psi$ admit the representation
\eqref{eq:sup_affine}, and the suprema
therein are achieved for $\mu = P$ and $\nu = Q$,
respectively by $f_\Phi\in \calF_\Phi$
and by $(f_\Psi, g_\Psi)\in \calH_\Psi$, 
then
\begin{enumerate}
\item The process $(N_t)_{t \geq 1}$
is bounded below by
$$R_t = \int f_\Phi d(P_t-P),$$
which is a (mean-zero) reverse martingale with respect to the exchangeable filtration $(\calE_t^X)$. 
\item The process $(M_{ts})_{t,s \geq 1}$ is bounded below by  
$(R_{t}^X + R_{s}^Y)_{t,s \geq 1}$,
where
$$R_{t}^X = \int f_\Psi d(P_t-P),~t \geq 1, \quad
\text{and}\quad
R_{s}^Y = \int g_\Psi d(Q_s-Q),~s \geq 1,$$
are (mean-zero) reverse martingales
with respect to $(\calE_t^X)$ and 
$(\calE_s^Y)$
respectively.
\end{enumerate}
\end{proposition} 
We begin by noting
that Proposition~\ref{prop:affine_rev_martingale}
implies the aforementioned 
upward bias of the plug-in
estimators $\Phi(P_t), \Psi(P_t,Q_s)$, including at arbitrary stopping times.
Indeed, by the optional stopping theorem (see for instance
\cite{durrett2019}, Theorem 4.8.3.), we can 
easily infer the following fact that we record formally for reference.
\begin{corollary}
\label{cor:upward_bias}
Assume the same conditions
as Proposition~\ref{prop:affine_rev_martingale},
and that 
the processes $(R_t), (R_t^X), (R_s^Y)$ therein
are uniformly integrable. 
Then, for any stopping times $\tau$ and $\sigma$ 
with respect to the canonical
forward filtrations
$(\sigma(X_1, \dots, X_t))_{t=1}^\infty$
and
$(\sigma(Y_1, \dots, Y_s))_{s=1}^\infty$ respectively, we have 
\begin{equation} 
\label{eq:upward_bias} 
\bbE[\Phi(P_\tau)] \geq \Phi(P), \quad 
  \bbE[\Psi(P_\tau, Q_\sigma)] \geq \Psi(P,Q).
\end{equation}
\end{corollary}

Furthermore, Proposition~\ref{prop:affine_rev_martingale}
can readily be
used to form upper confidence
sequences for $\Phi(P)$ or $\Psi(P,Q)$
on the basis of the reverse
martingales $(R_t)$,
or $(R_{t}^X)$ and $(R_s^Y)$. 
These processes are
simply sample averages, 
hence they can already be controlled
using the existing literature
on sequential mean estimation
(summarized
in Section~\ref{sec:background_confseq}). 
Nevertheless, for the purpose of being self-contained,
% we prove the following 
% result along similar lines
% as Theorem~\ref{thm:main}, using 
% reverse martingale techniques.
we use reverse martingale techniques to derive
such upper confidence sequences under tail conditions 
on $\calF_\Phi$ and $\calH_\Psi$,
in Proposition~\ref{prop:ucb} of Appendix~\ref{app:proofs_subsec3_2}. 
The following special case of this result 
% We state a special case of this result as follows, which 
will
be used repeatedly in Section~\ref{sec:applications}.  
\begin{corollary}
\label{cor:ucb_bounded}
Assume the same conditions as Proposition~\ref{prop:affine_rev_martingale}.
Assume further 
that for any $h \in  \calF_\Psi \cup \calG_\Psi$, $\diam(h(\calX)) \leq B < \infty$, and define
% there exist constants $-\infty< a < b < \infty$ such that 
% $a \leq f,g \leq b$ for all $f \in \calF_\Psi$ and
%  $g \in \calG_\Psi$. Let $B = b-a$, and define
\begin{equation} 
\label{eq:kappa_ts}
\kappa_{t} = \sqrt{\frac{\log\ell(\log_2 t) + \log(4/\delta)}{t}},\quad
\kappa_{ts} = \kappa_t+\kappa_s.
\end{equation}
 % $B > 0$ such that
% \begin{equation} 
% \label{eq:bounded_classes}
% \sup_{f \in \calF_\Psi} \sup_{x \in \supp(P)} |f(x)| \leq \frac {B}2, \quad
%   \sup_{g \in \calG_\Psi} \sup_{y \in \supp(Q)} |g(x)| \leq \frac{B}{2}.
% \end{equation}
 Then,  for any $\delta \in (0,1)$, we have
% $$\bbP\Big( \exists t \geq 1 : D(P_t\|Q) \leq  D(P\|Q)- B\kappa_{t} \Big) \leq \delta/2,\quad 
$\bbP( \exists t,s \geq 1 : \Psi(P_t,Q_s) \leq  \Psi(P,Q)- B\kappa_{ts}) \leq \delta/2.$
\end{corollary}
To summarize, under the assumptions and notation 
of Theorem~\ref{thm:main}
and Proposition~\ref{cor:ucb_bounded},
we deduce that a two-sided, two-sample,
$(1-\delta)$-confidence
sequence for $\Psi(P,Q)$ is given by
\begin{equation}
    \label{eq:two_sample_interval}
C_{ts} = \Big[\Psi(P_t,Q_s) -  \kappa_{ts},
\Psi(P_t,Q_s) +  \gamma_t^X + \gamma_s^Y \Big],
\end{equation} 
and similarly for the functional $\Phi$.

\subsection{\smash{On the choice of filtrations and stopping times}}
\label{sec:choice_filtrations}
The majority of confidence sequences derived
in past literature, such as those described in Section~\ref{sec:background_confseq},
employ  
martingales with respect to 
the canonical, or ``data-generating'' filtration.
A notable exception is the work of \cite{vovk2020}, which shows that
the power of certain sequential tests can be increased by coarsening 
the canonical filtration.
% , 
% whose constituent $\sigma$-algebras are
% generated by the observed data up to an increasing time.
It was similarly fruitful in our work to distinguish the 
data-generating filtration, with respect to which our processes do not appear to admit any martingale-type property, 
from a different filtration with respect to which our processes do admit a (reverse) martingale property.
To elaborate, let
$$\calD_t^X = \sigma(X_1, X_2, \dots, X_t), \quad
  \calD_s^Y = \sigma(Y_1, Y_2, \dots, Y_s),\quad
  t,s=1,2,\dots$$
  denote the canonical filtrations
  associated with each sequence
  of samples. 
Our bounds have implicitly assumed
that at any pair of times $(t,s)$, 
the practitioner has access to the information
encoded by the data-generating filtration
\begin{equation}
\label{eq:data_gen_filtration}
    \calD_{ts} = \calD_t^X \bigvee \calD_s^Y,
\quad t,s=1,2,\dots
\end{equation}
The process $(M_{ts})$ is naturally adapted 
to $(\calD_{ts})$, but we are not aware of it 
satisfying a martingale-type property with respect
to this filtration in general\footnote{Nevertheless, under some conditions, it
can be deduced from Theorem~\ref{thm:submart} that there exists a bivariate canonical filtration 
$(\calF_{ts})$ and a process $(\widetilde M_{ts})$, which has the same
distribution as $(M_{ts})$, such that $(\widetilde M_{ts})$ is a reverse submartingale with respect
to $(\calF_{ts})$ rather than $(\calE_{ts})$---see Theorem B of~\cite{rzeszut2020}.}. 
It is, however, also adapted to the exchangeable filtration
$\calE_{ts} = \calE_t^X \bigvee \calE_s^Y$, but unlike before, $(M_{ts})$ is also a reverse submartingale with respect to $(\calE_{ts})$. 
Our paper reinforces the somewhat underappreciated view in sequential analysis that filtrations should not be viewed as being ``inherent'' to the problem, or as tedious formalism for ensuring measurability, but instead viewed as design tools---a nonstandard choice of filtration can yield a powerful design tool.

\textbf{Validity at Stopping Times.}
To better understand the underlying role
of the filtration $(\calD_{ts})$, 
we shall now prove that the results of 
Theorem~\ref{thm:main} can equivalently be stated 
as bounds which hold at arbitrary stopping times
with respect to $(\calD_{ts})$. We focus on the two-sample
case in what follows.

In order to define a notion of stopping time which is 
suitable for our purposes, define the set
$$\widebar\bbN^2 = \bbN^2 \cup \{(t,\infty): t \geq 1\}
\cup\{(\infty,s):s \geq 1\} \cup \{(\infty,\infty)\},$$
for some symbols $(t,\infty), (\infty,s), (\infty,\infty)$.
We endow $\widebar\bbN^2$ with the natural
partial order, given by that of $\bbN^2$
(described in Section~\ref{sec:background_partial_order}), together 
with the following additional relations:
$(t,s) \leq (t',\infty)$ whenever $t \leq t'$ and $s \in \bbN$;
$(t,s) \leq (\infty,s')$ whenever $s \leq s'$ and $t \in \bbN$;
$u \leq (\infty,\infty)$ for all $u \in \widebar\bbN^2$.
A map $\eta: \Omega \to \widebar\bbN^2$
is said to be a stopping time with respect to a filtration
$(\calF_{ts})$ 
if $\{\eta = (t,s)\} \in \calF_{ts}$
for all $(t,s) \in \bbN^2$.  

Intuitively, the event $\{\eta = (t,s)\}$ 
indicates that the data collection from each of 
$P$ and $Q$ was terminated
at times $(t,s)$, whereas the event $\{\eta=(t,\infty)\}$
indicates that data was collected from $P$ until time $t$,
but indefinitely so from $Q$. Likewise, the event 
$\{\eta=(\infty,\infty)\}$ indicates that neither 
of the two data
collections were halted. 
With these definitions in place, 
we arrive at the following general equivalence.
\begin{proposition} 
\label{prop:stopping}
Let $(A_{ts})_{t,s=1}^\infty$ be a sequence of events
adapted to a forward filtration $(\calF_{ts})_{t,s =1}^\infty$.
%Alternative: adapted to \calD_{ts} \wedge \calE_{ts}$
Define for all $t,s \geq 1$,
\begin{equation} 
\label{eq:Ainf_stopping}
A_{t\infty} = \limsup_{s \to \infty} A_{ts}, \quad
  A_{\infty s}= \limsup_{t \to \infty} A_{ts}, \quad 
  A_{\infty\infty} = \left(\limsup_{t \to \infty} A_{t\infty} \right)
  \cup \left(\limsup_{s \to \infty} A_{\infty s}\right).
 \end{equation}
Then, for all $\delta \in (0,1)$, 
the following statements are equivalent.
\begin{enumerate}
\item[(i)] $\bbP\left(\bigcup_{t,s=1}^\infty A_{ts}\right) \leq \delta$.
% \item[(ii)] For any stopping time 
% $\eta=(\tau,\sigma)$ with respect to $(\calE_{ts})$, we have $\bbP(A_{\tau\sigma}) \leq \delta$.
\item[(ii)] For any stopping time
$(\tau,\sigma)$ with respect to $(\calF_{ts})$, %$(\calD_{ts})$ 
we have $\bbP(A_{\tau\sigma}) \leq \delta$.
\item[(iii)] For any random time 
$(T,S)$, not necessarily a stopping time, 
we have $\bbP(A_{TS}) \leq \delta.$
\end{enumerate}
\end{proposition}
The proof of Proposition~\ref{prop:stopping}
is given in Appendix~\ref{app:proofs_subsec3_3}. 
Analogues of Proposition~\ref{prop:stopping} 
for one-sample processes have previously been given
by \cite{howard2021}, \cite{ramdas2020}, and
\cite{zhao2016}, so our result is an extension of theirs to partially
ordered processes. 
In our setting, recall
that $M_{ts}$ is $(\calD_{ts} \bigwedge \calE_{ts})$-measurable.
While  
Proposition~\ref{prop:stopping}
could be reformulated
in reverse time, so that
$(\calF_{ts})$ can be
taken to be the modeling filtration
$(\calE_{ts})$, it is most interpretable
to take it to be the
data-generating filtration $(\calD_{ts})$. 
Doing so, under the assumptions of Theorem~\ref{thm:main}, 
leads for instance to the bound
\begin{equation} 
\label{eq:Mts_stopping_time}
\bbP\big\{ \Psi (P , Q)  \in C_{\tau\sigma}\big\} \geq 1-\delta
\quad \text{
for all stopping times } \eta=(\tau,\sigma) \text{ 
with respect to } (\calD_{ts}),
\end{equation}
where $C_{ts}$ denotes the two-sided 
interval~\eqref{eq:two_sample_interval}, understood with 
conventions for infinities which can be deduced from equation~\eqref{eq:Ainf_stopping}.

\textbf{Alternate Data-Generating Filtrations.}
Though we presumed the data-generating filtration 
\eqref{eq:data_gen_filtration} 
throughout our development, 
slightly tighter confidence sequences can be obtained
if the user has access to additional information.
For instance, our confidence sequences
hold uniformly over arbitrary pairs of time $(t,s)$, but
such flexibility is unnecessary if the practitioner
knows the order in which sample points from $P$ and $Q$
arrive. We illustrate two such examples below, focusing
on lower confidence sequences:
% \item Theorem~\ref{thm:main}(ii)
% makes no assumptions on the relative sizes of 
% the sample sizes $t$ and $s$. While doing so 
% offers a great deal of flexibility,
% condition \eqref{eq:stitching_fn_main}
% on $g$ can be tightened if
% there is a known order in which sample points 
% $(X_t)$ and $(Y_s)$ appear.
% For instance, if these sequences of samples
% become available to the practitioner
% in pairs $(X_t, Y_t)_{t=1}^\infty$, the following
% two-sample bound
% is slightly sharper than
% that of Theorem~\ref{thm:main}(ii),
% and can be proven by a straightforward
% adaptation,
% \begin{align*}
% \bbP&\left\{\exists t \geq 1: D(P_{t} \| Q_t) \geq D(P \| Q) +(\psi_{tt}^*)^{-1}\Big(\log \osfn(\log_2 t) + \log(1/\delta)\Big)\right\}
% \leq \delta.
% \end{align*}
  
\begin{enumerate} 
% \item \textbf{Separate Samples.} When the divergence $D$ satisfies 
% the triangle inequality, possibly conservative
% bounds for $M_{ts}$ can be obtained by separately
% bounding the processes $(M_{t0})$ and $(M_{0s})$, due to 
% the inequality $|M_{ts}| \leq M_{t0} + M_{0s}$.
% In this case, it is natural to view the two sequences
% $(X_t)$ and $(Y_s)$ as separate, so that the practitioner has access
% to two distinct data-generating filtrations,
% $$\calD_t^X = \calC_t^X, \quad 
%   \calD_s^Y = \calC_s^Y, \quad
%   t,s,=1, 2, \dots$$

% \item \textbf{Arbitrarily-Ordered Samples.} When the observations $X_t$
% and $Y_s$ arrive in an arbitrary and possibly data-dependent order,
% the data-generating filtration is given by
% $$\calD_{ts} = \sigma(\calC_t^X \cup \calC_s^Y), \quad t,s =1,2,\dots$$
% % \item \textbf{One Sample.} When the distribution $Q$ is, in fact, known, 
% % the sample $(Y_s)$ is redundant 
\item[(i)] \textbf{Paired Samples.} When the observations $X_t$ and $Y_t$ are presumed to arrive
at the same time, in pairs $(X_t, Y_t)$, the data-generating
filtration may be replaced by
$$\calD_t = \sigma(X_t, Y_t), \quad t=1,2,\dots$$
In this case, following along similar lines as before, 
the following two-sample bound may be established, 
and is tighter than that of Theorem~\ref{thm:main}(ii),
\begin{align}
\label{eq:paired_samples}
\bbP&\left\{\exists t \geq 1: \Psi(P_{t}, Q_t) \geq \Psi(P , Q) +(\psi_{tt}^*)^{-1}\big(\log \osfn(\log_2 t) + \log(1/\delta)\big)\right\}
\leq \delta.
\end{align}
Unlike Theorem \ref{thm:main}, we note that the bound \eqref{eq:paired_samples}
can be taken to hold without assuming that $(X_t)$ and $(Y_s)$ are independent
of each other.

\item[(ii)] \textbf{Samples Ordered by External Randomization.} As a generalization of the previous
point, assume the observations
$X_t$ and $Y_s$ arrive in a possibly random
order which is independent of the data. Specifically, 
let $(\iota_n)_{n\geq 1}$ denote a sequence of random variables
taking values in $\{0,1\}$, which are independent of $(X_t)$ and $(Y_s)$,
but possibly dependent on an external source of randomness $U$, say distributed
uniformly on $[0,1]$. 
Let $t(n) = \sum_{i=1}^n \iota_n$, and
$s(n) = n-t_n$ 
% $O_n = \iota_n X_{t(n)}  + (1-\iota_n) Y_{s(n)}$ 
so that at any time $n \geq 1$, 
the practitioner observes 
$\iota_n X_{t(n)}  + (1-\iota_n) Y_{s(n)}$.
In this case, one has access to the filtration
$\calI_n = \sigma(U, \iota_1, \iota_2, \dots, \iota_n)$, $n \geq 1$,
which determines the order in which the sample points $X_t, Y_s$
arrive, as well as to the data-generating filtration
\begin{align} 
%\label{eq:efron_filtration1} 
\widebar \calD_n = \widebar \calD_{t(n)}^X\bigvee 
\widebar \calD_{s(n)}^Y, \quad n\geq 1,
\end{align} 
where $\widebar \calD_{t(n)}^X$
and $\widebar \calD_{s(n)}^Y$ are defined similarly as follows:
\begin{align} 
%\label{eq:efron_filtration2}
\widebar \calD_{t(n)}^X = \{A \in \calF: A \cap \{t(n)=t\} \in \widebar\calD_t^X, \ \forall t \geq 1\},
\quad \text{where}\quad 
\widebar \calD_{t}^X = \sigma(U, X_1, \dots, X_t).
\end{align}
Note that we could have assumed that the sequence $(\iota_n)$ is 
fully deterministic, in exchange for simpler notation. 
However, there are many situations, like clinical trials, 
in which we may wish to use external randomization (encoded by $U$) 
to determine how to obtain the next data point; for example, \citet{efron1971} shows how to adaptively randomize participants while encouraging balance between $t(n)$ and $s(n)$.
Under this setting, it can be shown that
$(\Psi (P_{t(n)}, Q_{s(n)}))_{n \geq 1}$ is a reverse submartingale, and assuming for simplicity that
$\psi_n = \psi_{n0} = \psi_{0n}$, one may obtain the confidence sequence
$$\bbP\big\{ \exists n \geq 1: 
\Psi(P_{t(n)}, Q_{s(n)}) 
\geq \Psi(P,Q) + (\psi_{\bar n}^*)^{-1} \big( \log \ell(\log_2 n) + \log(1/\delta)\big)\big\} \leq \delta.$$
\end{enumerate}
In contrast to the above two settings, 
our confidence sequence $C_{ts}$ satisfies the guarantee~\eqref{eq:confseq},
which is uniform over pairs $(t,s) \in \bbN^2$, and therefore yields 
valid coverage even if the orderings $t(n)$ and $s(n)$ depend arbitrarily on the samples observed from $P$ and $Q$. 

\section{Explicit corollaries for common divergences}
\label{sec:applications}
We now specialize the confidence sequence $C_{ts}$ %Theorem~\ref{thm:main} 
to several examples of divergences including
IPMs (Sections \ref{sec:dkw}, \ref{sec:mmd}), 
optimal transport costs (Section~\ref{sec:wasserstein}),
$\varphi$-divergences (Section~\ref{sec:phi_divergence}),
and divergences smoothed by convolution (Section~\ref{sec:smooth}). 
Moving beyond divergences, 
we also derive time-uniform generalization error bounds 
for binary classification problems~(Section~\ref{sec:rademacher}), and confidence sequences 
for multivariate means~(Section~\ref{sec:multivariate_means}).
These special cases will
illustrate how our framework can 
be used to port existing
fixed-time concentration inequalities
to time-uniform ones, typically at the expense
of iterated logarithmic factors. In these cases, any improvements
to existing fixed-time concentration inequalities
would typically carry over to our time-uniform setting.
% In several cases, we also derive new tail bounds which
% were not previously available in the non-sequential setting. 
Though our focus is on non-asymptotic bounds,
in Section~\ref{sec:lil}, we
also show how Theorem~\ref{thm:main} can be used to derive a one-sided analogue
of the classical law of the iterated
logarithm, for several divergences between an empirical and true underlying
measure.
We defer all proofs to 
Appendix~\ref{app:proofs_sec5}.

\subsection{Kolmogorov-Smirnov Distance}
\label{sec:dkw}
Theorem~\ref{thm:main} leads to a sequential 
analogue of the classical Dvoretzky-Kiefer-Wolfowitz (DKW) inequality~\citep{dvoretzky1956, massart1990}, based on distinct
techniques than those of~\cite{howard2019} and \cite{odalric-ambrym2020}.
Let $P$ be any distribution 
over $\bbR$ 
with cumulative distribution function 
(CDF) $F$. 
Let $F_t(x) = (1/ t) \sum_{i=1}^n I(X_i \leq x)$
denote the empirical CDF of $F$.
\begin{corollary} 
\label{cor:dkw}
For any $\delta \in (0,1)$, 
$$\bbP\left(\exists t \geq 1: \norm{F_t-F}_\infty
\geq \sqrt{\frac \pi{t}} +  2\sqrt{\frac 2 t \Big[ \log \ell(\log_2 t) + \log(1/\delta)\Big]}\right) \leq \delta.$$
 \end{corollary}
Notice that Corollary~\ref{cor:dkw} involves the term 
$\log(1/\delta)$, as opposed to the term $\log(2/\delta)$
which appears in the classical DKW inequality. This is due to the one-sidedness of the bounds in Theorem~\ref{thm:main}. 
The price to pay is the additional additive term $ \sqrt{ \pi/{t}}$, which is an upper bound on the expectation
of $\|F_{\lceil t/2\rceil}-F\|_\infty$. 

Corollary \ref{cor:dkw} and Proposition~\ref{prop:ucb} 
give rise to a sequential analogue
of the celebrated Kolmogorov-Smirnov
two-sample test. In what follows, 
$Q$ denotes a second distribution over
$\bbR$ with CDF $G$, 
and empirical CDF $G_s(y)=(1/s)\sum_{i=1}^s I(Y_i \leq y)$. Recall also the sequence  
$(\kappa_{ts})$ defined in 
equation~\eqref{eq:kappa_ts}.
\begin{corollary}
\label{cor:ks_test}
Given $\delta \in (0,1)$,
set $\gamma_{ts} = \sqrt{\pi/t} + \sqrt{\pi/s} + 2\sqrt{\frac {2ts} {t+s} \big[ \log \tsfn(\log_2 t + \log_2 s) + \log(2/\delta)\big]}$. Then, 
$$\bbP\big(\forall t,s \geq 1: -\gamma_{ts}\leq \norm{F-G}_\infty -
\norm{F_t-G_s}_\infty \leq \kappa_{ts}\big)\geq 1-\delta.$$
% $$\bbP\left(\exists t,s \geq 1:
% \norm{F_t-G_s}_\infty \geq 
% \norm{F-G}_\infty 
% + \gamma_{ts} \right) \leq \delta.$$
In particular, the sequential Kolmogorov-Smirnov test
which rejects the null hypothesis $H_0:P=Q$
when $\norm{F_t-G_s}_\infty > \gamma_{ts}$
has type-I error controlled at $\delta/2$.
\end{corollary}

We now turn our attention to another popular IPM that is based on reproducing kernels.

\subsection{Maximum Mean Discrepancy, V-Statistics, 
and U-Statistics}
\label{sec:mmd}
The kernel Maximum Mean Discrepancy (MMD) is an IPM measuring the distance between embeddings
of distributions in a reproducing kernel Hilbert space (RKHS). 
We provide a brief definitions in what follows, and refer the reader to \cite{scholkopf2018} for further details.  
Let $K:\bbR^d \times \bbR^d \to \bbR_+$ be a Mercer kernel, that is a 
symmetric and continuous function such that for any finite set of points
$x_1, \dots, x_n \in \bbR^d$, the matrix $(K(x_i,x_j))_{i,j=1}^n$ is positive semidefinite.
The RKHS $\calH$ corresponding to $K$ is the closure of the set
$$\calH_0 = \left\{\sum_{i=1}^{k} \alpha_i K(\cdot,x_i): \alpha_i \in \bbR, \ x_i \in \bbR^d, k \geq 1\right\},$$
endowed with the inner product and norm
$$\langle f, g\rangle_\calH = \sum_{i=1}^k \sum_{j=1}^{k'} \alpha_i \beta_j K( x_i, y_j),\quad \norm{f}_\calH = \sqrt{\langle f, f\rangle_\calH},$$
where $f = \sum_{i=1}^k \alpha_i  K(\cdot, x_i)$ and $g = \sum_{j=1}^{k'} \beta_j K(\cdot, y_j)$
denote the expansions of any two functions $f, g \in \calH_0$. 
The MMD is defined as the IPM over the unit ball $\calJ = \{f \in \calH: \norm f_\calH \leq 1\}$ of~$\calH$.
% ,
% that is,  for all distributions $P,Q$ on $\bbR^d$, 
% $$\mmd(P,Q) = \sup_{f \in \calH} \int fd(P-Q).$$
%satisfying $0 \leq K \leq B$, for some $B > 0$. 
%Let $\calH$ be a reproducing kernel Hilbert space (RKHS) with respect
%to a kernel $K:\bbR^d \times \bbR^d \to \bbR$ satisfying $0 \leq K \leq B$, for some $B > 0$. 
%where $\calG = \{f \in \calH: \norm f_\calH \leq 1\}$, and $\norm\cdot_\calH$
%denotes the norm induced by the inner product in $\calH$. 
% Let $P,Q$ be any two distributions over $\bbR^d$, and let $(X_t)_{t=1}^\infty$
% and $(Y_s)_{s=1}^\infty$ denote respective iid samples from $P$ and $Q$,
%Let
%$P_t = \frac 1 t \sum_{i=1}^t \delta_{X_i}$ and $Q_s = \frac 1 s \sum_{i=1}^s \delta_{Y_i}$
%denote their corresponding empirical measures, and note the representation
%Given sequences of i.i.d. observations $(X_t)_{t=1}^\infty$
%and $(Y_s)_{s=1}^\infty$ from $P$ and $Q$ respectively
% with corresponding empirical measures $P_t, Q_s$.
The plug-in estimator of the MMD admits the following representation, which makes its computation 
straightforward
\begin{equation} 
\label{eq:mmd_v_stat}
\mmd(P_t, Q_s) = \sqrt{\frac 1 {t^2} \sum_{i,j=1}^t K(X_i, X_j) + \frac 1 {s^2} \sum_{i,j=1}^s K(Y_i,Y_j) - 
								\frac 2 {st} \sum_{i=1}^t\sum_{j=1}^s K(X_i, Y_j)}.
\end{equation}								
In particular, if $s=t$ and the data
are available in pairs
$Z_t = (X_t, Y_t)$ for all $ t\geq 1$, 
as in Section~\ref{sec:choice_filtrations}(1), 
the above expression may be rewritten
as a square root of a second-order V-Statistic,
\begin{equation}
\label{eq:one_sample_mmd_v_stat}
\mmd(P_t, Q_t) = \sqrt{\frac 1 {t^2} \sum_{i,j=1}^t J(Z_i, Z_j)},
\end{equation}
where $J((x,y), (x',y')) = K(x,x') + K(y,y') - 2 K(x,y')$, for any $x,x',y,y' \in \bbR^d$. 
Assuming the kernel $K$ is bounded, we derive a sequential concentration bound for this statistic as follows.
%denote respective iid samples from $P$ and $Q$. Let
%$P_t = \frac 1 t \sum_{i=1}^t \delta_{X_i}$ and $Q_s = \frac 1 s \sum_{i=1}^s \delta_{Y_i}$
%denote their corresponding empirical measures, and note the representation
%We arrive at the following result.					
\begin{corollary}
\label{cor:mmd}
Let $P,Q \in \calP(\bbR^d)$. Assume that %for all $x,y \in \bbR^d$,
$\sup\{K(x,y):x,y \in \bbR^d\} \leq B <\infty$. 
For any $\delta \in (0,1)$, 
define 
$$\gamma_{ts} = 2 \sqrt{2 B} (t^{-\frac 1 2} + s^{-\frac 1 2})
 + 4\sqrt{\frac{B(t+s)}{ ts} \Big[ \log \tsfn(\log_2 t +  \log_2 s) + \log(2/\delta)\Big]},$$
and let $(\kappa_{ts})$ be the sequence
defined in equation~\eqref{eq:kappa_ts}. Then,  
\begin{align*}
\bbP\Big(\forall t,s \geq 1: -\gamma_{ts} \leq 
 \mmd(P,Q)-\mmd(P_t, Q_s) \leq  
 2\sqrt B\kappa_{ts} \Big) \geq 1-\delta.
 \end{align*}
 \end{corollary}
Assuming that the stitching
function $g$ is bounded above by a polynomial,
Corollary \ref{cor:mmd} provides a confidence
sequence for $\mmd(P, Q)$ scaling at the rate
$O(\sqrt{\log\log (t\vee s)/(t\vee s)})$. Up to the iterated logarithmic factor, we recover the fixed-time 
rate obtained by \cite{gretton2012} (Theorem 7),
which was shown to be minimax optimal by \cite{tolstikhin2016}. In the setting of equal sample sizes $t=s$ (cf.\ 
Section~\ref{sec:choice_filtrations}(1)),
it is well-known that the V-Statistic $\mmd^2(P_t, Q_t)$ has first-order degeneracy, so that
$\mmd^2(P_t, Q_t) = O_p(1/t)$ when $P=Q$
\citep{lee1990}.
On the other hand, the bound $|\mmd^2(P_t, Q_t) - \mmd^2(P,Q)| = O_p(1/\sqrt t)$ is tight
when $P,Q$ are fixed and $P \neq Q$. 
On the squared scale, the bound of Corollary \ref{cor:mmd} adapts to these distinct rates of convergence, since it implies that with probability
at least $1-\delta$,  
\begin{equation}
\label{eq:mmd_adaptive}
\forall t \geq 1, \quad
|\mmd^2(P_t, Q_t) - \mmd^2(P,Q) | = O\left(\mmd(P,Q) \sqrt{\frac{\log \log t}{t}} + \frac{\log \log t }{t}\right).
\end{equation}
Above, the right-hand side decays at the rate $O(\sqrt{(\log \log t)/t})$
in general, but improves to  
$O\left( (\log \log t)/t \right)$ 
when $P=Q$.
Similar considerations
are  discussed by \cite{gretton2012}.

While the plug-in estimator
$\mmd(P_t,Q_s)$ is typically upwards biased,
the \emph{squared} MMD also admits a widely-used unbiased 
estimator~\citep{gretton2012} 
obtained by replacing the V-Statistics in equations~\eqref{eq:mmd_v_stat} and~\eqref{eq:one_sample_mmd_v_stat} by U-Statistics. 
We derive confidence sequences 
for $\mmd^2(P,Q)$ based on this estimator
in Appendix \ref{app:proofs_mmd}. 
The bounds therein do not
adapt to the distinct rates of convergence
described above, however,
therefore we recommend those of Corollary~\ref{cor:mmd}
in practice. 
 
We conclude this section
with a more general discussion
of sequential inference
based on U-
and V-Statistics, for expectation functionals of the form
$$\Phi : \calP(\bbR^d)
 \to \bbR, \quad 
 \Phi(P) = \iint h(x,y)dP(x)dP(y),$$
where  $h:\bbR^d \times \bbR^d \to \bbR$ is
a symmetric function. 
 Following \citet{lee1990,serfling2009},
 the U- and V-Statistics corresponding
 to $\Phi$ are respectively defined by
 $$U_t = \bbE[h(X_1, X_2)|\calE_t^X] = 
 \frac {2}{t(t-1)}\sum_{1 \leq i < j \leq t} h(X_i,X_j), \quad
  V_t = \Phi(P_t) = 
  \frac 1 {t^2} \sum_{i,j=1}^t 
  h(X_i,X_j).$$
The above representation  immediately implies that 
$(U_t)$ is a reverse martingale 
whenever it is integrable---a
fact which can equivalently
be derived from the leave-one-out
property~\eqref{eq:loo},
which holds for $(U_t)$
with equality. 
Notice that this property holds 
irrespective of the kernel
$h$, so long as it is
symmetric. On the other hand,
the following can be said
about $(V_t)$.
\begin{proposition}
\label{prop:v_stats}
Assume that $h:\calX \times \calX \to \bbR_+$ is a  continuous, symmetric, 
positive definite kernel
over a compact set $\calX \subseteq \bbR^d$. 
Then $\sqrt{\Phi(\cdot)}$
is a convex functional, 
thus $(\sqrt{V_t})$ and $(V_t)$
are reverse submartingales
with respect to $(\calE_t^X)$.
\end{proposition}
The proof is in Appendix \ref{app:proofs_mmd}. We do not generally expect
the above result to hold 
for any symmetric kernel $h$,
because the functional $\Phi$
is akin to a quadratic form, 
which may be nonconvex
if its kernel is not positive
semidefinite. 
Under the above results, a straightforward
generalization of Theorem~\ref{thm:main}
can  be used to derive
two-sided confidence sequences for $\Phi(P)$ centered at $U_t$ 
for all symmetric~$h$,
or lower 
confidence sequences for~$\sqrt{\Phi(P)}$  and~$\Phi(P)$
based on~$(V_t)$ 
for all positive definite~$h$ 
(which may be coupled with 
upper confidence sequences
similarly
as in Section~\ref{sec:ucb}). 
While these considerations
make $(U_t)$-based
confidence sequences seem
attractive, we recall
that $(V_t)$-based
confidence sequences sometimes
have the advantage
of providing rate-optimal
inference even when $\Phi$
is degenerate. 
% We illustrated
% this fact in the special
% case of the MMD above, and we discuss it further in Appendix~\ref{app:proofs_mmd}.

\subsection{Optimal Transport Costs}
\label{sec:wasserstein}

Let
$c:\calX\times\calX\to\bbR_+$ be a nonnegative cost function,
and assume for simplicity that $c$ is bounded above over $\calX$ by $\Delta := \sup_{x,y\in \calX} c(x,y) < \infty$. We derive confidence sequences
for the optimal transport cost $\calT_c(P,Q)$, which will depend
on an upper bound $\alpha_{c,ts}$ for the bias
of the empirical plug-in estimator,
\begin{equation} 
\label{eq:wasserstein_alpha} 
\bbE \big[\calT_c(P_t, Q_s)\big] - \calT_c(P,Q)  \leq \alpha_{c,ts}.
\end{equation}
Such bounds have been 
derived in the literature for various choices
of the cost~$c$. For instance, 
one may
take~$\alpha_{c,ts}$ to scale
as~$(t\wedge s)^{-\alpha/d}$ for $d\geq 5$ when $c$ 
is an $\alpha$-H\"older smooth 
cost
for some $\alpha \in [1,2]$~\citep{manole2021a, niles-weed2019, chizat2020}, 
and as $(t\wedge s)^{-1/2}$ when
$\calX$ is  finite or 
one-dimensional, 
with suitable costs $c$~\citep{sommerfeld2018,forrow2018,munk1998, manole2019}. 
We refer to~\cite{hundrieser2022} for further examples.
The leading constants in these rates cannot always be made explicit, thus their use below will be of primarily theoretical interest. 
\begin{corollary} 
\label{cor:wasserstein}
Let $c:\calX\times \calX\to\bbR_+$
be a lower semi-continuous
cost function bounded above by $\Delta$, assume that 
equation \eqref{eq:wasserstein_alpha} is satisfied, and recall that $\bar t = \lceil t/2\rceil$. 
Furthermore, let $(\kappa_{ts})$ be the sequence
defined in equation~\eqref{eq:kappa_ts}.
Then, for any $\delta \in (0,1)$, 
\begin{align*} 
\bbP\bigg(\forall t,s\geq 1 :
-\alpha_{c,\bar t\bar s} - 2\Delta\sqrt{\frac{2(t+s)}{ts} \Big[ \log g(\log_2 t + \log_2 s) + \log(2/\delta)\Big]} &\\
&\mkern-135mu \leq  \calT_c(P,Q) -
\calT_c(P_t, Q_s) \leq 
\Delta \kappa_{ts}\bigg)  
  \geq 1-\delta.
 \end{align*}
\end{corollary} 
Corollary \ref{cor:wasserstein}
exhibits a confidence sequence for $\calT_c(P,Q)$ 
whose length
is typically dominated by the expected deviation bound
$\alpha_{c,\bar t\bar s}$. The potentially
severe dependence on dimensionality
in these rates can limit the applicability of
Corollary \ref{cor:wasserstein} in high-dimensional
problems. 
Section~\ref{sec:smooth} derives confidence sequences
for smoothed 1-Wasserstein distances, 
which admit significantly improved dimension dependence.

\subsection{$\varphi$-Divergences over finite sets}
\label{sec:phi_divergence}
Let $P$ be a probability distribution supported on a set $\calX = \{a_1, \dots, a_k\}$ of finite cardinality $k \geq 2$.
Set $p_j = P(\{a_j\})$ for all $j=1, \dots, k$. 
In this setting, we write the empirical measure as
$$P_t = \frac 1 t \sum_{i=1}^t \delta_{X_i} = \frac 1 t  \sum_{j=1}^k C_j \delta_{a_j},\quad \text{where}\quad
C_j = \sum_{i=1}^t I(X_i = a_j), \quad j=1, \dots, k.$$
The vector $(C_1, \dots, C_k)$ can be viewed
as a random sample from a multinomial experiment with $t$ trials and $k$ categories
with probabilities $(p_1, \dots, p_k)$. 
Concentration inequalities for $\varphi$-divergences $D_\varphi$ between the empirical measure and the 
finitely-supported true
distribution $P$ have received significant
attention in the offline setting. Here, 
we show how such results can be used together with   Theorem~\ref{thm:main} 
to obtain time-uniform bounds on $D_\varphi(P_t\|P)$. We focus on the 
 Kullback-Leibler divergence and Total Variation distance in what follows.

\textbf{Kullback-Leibler Divergence.}
 Tight upper bounds on the moment generating 
function of the scaled Kullback-Leibler
divergence $t\KL(P_t\|P)$
have recently been derived by 
\cite{guo2020} (see
also \cite{agrawal2020a}), who
prove 
\begin{equation} 
\bbE[\exp(\lambda t\KL(P_t\|P))] \leq G_{k,t}(\lambda)
:=\sum_{x_1, \dots, x_k} {t \choose {x_1, \dots, x_k}} \prod_{i=1}^k [\lambda x_i/t + (1-\lambda) p_i]^{x_i},
\end{equation}
for all $ \lambda \in [0,1]$.
\cite{guo2020} show that this upper
bound is nearly tight, in the
sense that it nearly matches the 
scaling in $k$ and $t$ of the moment
generating function of the limiting 
distribution of $t\KL(P_t\|P)$. 
Furthermore, the value of $G_{k,t}(\lambda)$
does not depend on $P$ (\cite{guo2020}, Proposition~1). 
Nevertheless, $G_{k,t}$ cannot
easily be used in Theorem~\ref{thm:main}, 
since the optimization problem $\sup_{\lambda \in [0,1]} \big\{\lambda u - G_{k,t}(\lambda)\big\},$
for any $u > 0,$
%\end{align}
is non-convex. 
\cite{guo2020} instead
derive several 
closed-form sequences $(\lambda_t)$
which approximately solve
this maximization problem.
We derive a sequential analogue
of their bounds in terms
of a generic choice of such a sequence.
% $(\lambda_t) \subseteq [0,1]$. 
The following result
is obtained by repeating
a similar stitching argument
as that of the proof of Theorem 
\ref{thm:main}.
\begin{proposition} 
\label{prop:kl_finite}
Let $\delta \in (0,1)$, and let 
$P$ be a distribution supported on a 
finite set of size $k \geq 2$. 
Let $(\lambda_t)_{t=1}^\infty \subseteq [0,1/2]$ be
a sequence of real numbers such that 
$$\gamma_t = \frac {2} {\lambda_{t}t} \log\left(\frac{ G_{k,\lfloor t/2\rfloor }(2\lambda_t) \osfn(\log_2 t)}{\delta}\right),
\quad t \geq 1,$$
is a nonincreasing sequence in $t$. Then, 
$\bbP\left\{\exists t \geq 1: \KL(P_t\|P) \geq \gamma_t \right\} \leq \delta.$
\end{proposition}
It can be seen
that the fitted probability
vector $(C_1/t, \dots, C_k/t)$ is precisely
the maximum likelihood estimator
of $(p_1, \dots, p_k)$, 
and that the
scaled Kullback-Leibler divergence 
$t\KL(P_t\|P)$ is a multiple
of the log-likelihood ratio of
$(p_1, \dots, p_k)$.
Proposition~\ref{prop:kl_finite}
therefore leads to a confidence
sequence for the probability vector
$p$ on the basis of the classical
likelihood ratio statistic. 

\textbf{Total Variation Distance.} We now similarly derive
time-uniform bounds for the discrete
Total Variation distance 
$\norm{P_t-P}_{\mathrm{TV}} :=
\frac 1 2 \sum_{j=1}^k \left|(C_j/t) - p_j\right|$. 
The following Corollary follows from Theorem~\ref{thm:main}
using elementary tail bounds for the Total Variation distance 
(see for instance \cite{berend2013}), together with an expectation bound due to~\cite{kamath2015}.
\begin{corollary}
\label{cor:tv}
For all $\delta \in (0,1)$, we have
\begin{align*}
\bbP \left\{ \exists t \geq 1 : \norm{P_t - P}_{\mathrm{TV}} \geq 
\sqrt{\frac{k}{\pi t}} + 2\left(\frac{2k}{t}\right)^{\frac 3 4} + 
  2\sqrt{\frac 2 t \Big[ \log \ell(\log_2 t) + \log(1/\delta) \Big] }\right\} \leq \delta.
\end{align*}
\end{corollary} 
Up to a polylogarithmic factor,
the bound of Corollary \ref{cor:tv} scales
at the parametric rate of convergence when the alphabet
size $k$ is fixed. 
For general distributions with uncountable support, 
such rates are not 
achievable under the Total Variation distance
due to the lack of absolute continuity of $P_t$ with respect to $P$.
The following subsection studies a notable exception, 
in which parametric rates are retained 
when the measures are smoothed
by convolution with a kernel admitting fixed bandwidth.

 \subsection{Smoothed divergences and 
 differential entropy}
 \label{sec:smooth}
Let $K:\bbR^d \to \bbR_+$ denote a smoothing kernel, that is, a nonnegative
and continuous function satisfying $\int_{\bbR^d}  K(x)dx = 1$. Given a bandwidth 
$\sigma > 0$, let 
$\calK_\sigma$ be the probability measure
admitting density $K_\sigma(x)=(1/\sigma^d)K(x/\sigma)$ with respect to the Lebesgue measure.
Let $D$ denote a convex divergence, and
define its smoothed counterpart by
$$D^\sigma : (P,Q) \longmapsto D(P \star \calK_\sigma \| Q \star \calK_\sigma).$$
It can be directly verified that $D^\sigma$
is itself a convex divergence, due to the linearity
of the convolution operator. 
Theorem~\ref{thm:main} can therefore
be used to derive a confidence sequence
for $D^\sigma(P\|Q)$ based on the
plug-in estimator $D^\sigma(P_t\|Q_s)$. 
We emphasize that this estimator is sensible
even if the original divergence $D$ requires
absolute continuity of the distributions being 
compared, as is the case for $\varphi$-divergences. 
In such cases, $D^\sigma$ forms a proxy of 
$D$ which can be estimated
using the empirical plug-in estimator. 
We refer to \cite{goldfeld2020b}
for upper bounds on $\bbE D^\sigma(P_t, P)$
under a wide range of divergences $D$. 
Smoothing by Gaussian convolution has also recently
been studied as a means of regularizing
optimal transport problems, and thereby reducing
the curse of dimensionality in estimating
Wasserstein distances. 
For instance, \cite{goldfeld2020c},
\cite{goldfeld2020} show that
that the empirical measure converges to $P$
at the parametric rate  under $W_1^\sigma$,
in expectation, contrasting 
the unavoidable $t^{-1/d}$ rate 
for this problem under $W_1$ itself
\citep{singh2019}.

Motivated by these two applications, we 
show in what follows how Theorem~\ref{thm:main} can
be used to derive confidence sequences
for the smoothed Total Variation distance, 
and for the smoothed 1-Wasserstein distance.
% , contrasting
% the bounds we previously obtained in Sections~\ref{sec:wasserstein}
% and~\ref{sec:phi_divergence}. 
The results which follow are obtained
by first  deriving upper bounds on the moment
generating functions of $D^\sigma(P_t\| P)$,
and second, invoking the upper
bounds of \cite{goldfeld2020b} on the 
expectation of this quantity. 
In order to appeal to their results,
we  assume in what follows
that $K(x) = e^{-\norm x_2^2/2}/\sqrt{2\pi}$ is taken to
be the standard Gaussian kernel.

We first recall a tail assumption which will be used
in the sequel. Given a metric $d$ on $\calX$, we say a measure $P \in \calP_1(\calX)$ satisfies
the $T_1(\tau^2)$ inequality with respect to $d$, 
for some $\tau > 0$, if 
$$\calT_d(\mu,P) \leq \sqrt{2\tau^2\KL(\mu\|P)}, \quad \text{for all }  \mu \in \calP_1(\calX).$$
Such inequalities are at
the heart of the transportation method for
deriving fixed-time
concentration inequalities---we 
refer to \cite{gozlan2010} for a survey. 
For our purposes,
transportation inequalities
are known to provide a natural tail assumption
on $P$ in order to guarantee sub-Gaussian
concentration of empirical Wasserstein distances:
\cite{niles-weed2019} (Theorem 6)
prove that $P$
satisfies $T_1(\tau^2)$
if and only if $W_1(P_t, P)$ is $(\tau^2/t)$-sub-Gaussian.
By extending their result to smoothed Wasserstein distances, 
we arrive at the following.
\begin{proposition}[Smoothed Divergences]
 \label{prop:smooth}
Let $\delta \in (0,1)$, $\sigma > 0$, and let
$P \in \calP(\bbR^d)$.
\begin{enumerate}
\item[(i)] (Total Variation Distance)
Assume $P$ is $\tau^2$-sub-Gaussian for some
$\tau > 0$.
Then,
$$\bbP\left(\exists t \geq 1: \norm{P_t-P}_{\mathrm{TV}}^\sigma 
 \geq \frac{c_d}{\sqrt t} + 4\sqrt{\frac 2 t \Big[\log \osfn(\log_2 t) + \log(1/\delta)\Big]}\right)
 \leq \delta,$$
 where   $c_d = \sqrt 2 \left(\frac 1 {\sqrt 2} + \frac {\tau}{\sigma}\right)^{\frac d 2} e^{\frac{3d}{16}}$.
\item [(ii)] ($W_1$ Distance) 
Assume $P$ satisfies the $T_1(\tau^2)$ inequality 
with respect to $\norm\cdot_2$, for some $\tau > 0$. Then, 
$W_1^\sigma(P_t, P)$ is $(\tau^2/t)$-sub-Gaussian, and for all $\delta \in (0,1)$,
\begin{equation}
\label{eq:smooth_wasserstein_bound}
\bbP\left(\exists t \geq 1: W_1^\sigma(P_t, P) \geq \frac{C_{d}}{\sqrt t} +2\sqrt{\frac{  \tau^2}{t} \Big[ \log \osfn(\log_2 t) + \log(1/\delta) \Big] }\right) \leq \delta,
\end{equation}
where $C_{d} = 2\sqrt{d\sigma^2}\left(\frac 1 {\sqrt 2} + \frac{\tau}{\sigma}\right)c_d$.%^{\frac d 2 +1} e^{\frac{3d}{16}}.$
\end{enumerate}
\end{proposition} 
Extensions of Proposition~\ref{prop:smooth}
to the two-sample setting are straightforward and 
omitted for brevity.  
Proposition~\ref{prop:smooth}(i)
yields a confidence sequence for the
smoothed Total Variation under a mere moment condition. Such a result
could not have been obtained by our framework
in the absence of smoothing, 
except in the special case of 
Section \ref{sec:phi_divergence}
where $P$ was assumed to have
countable support.
Proposition~\ref{prop:smooth}(ii)
contrasts our earlier
Corollary \ref{cor:wasserstein}, 
which implied a confidence sequence
for the 1-Wasserstein distance scaling at the rate $O(t^{-1/d})$ for $d \geq 3$. Smoothing removes the dimension
dependence from the rate itself when $\sigma$ is held fixed, 
although the constant $C_{d}$ continues
to grow exponentially in $d$.
It can be deduced from Theorem 1  of
\cite{goldfeld2020b}, combined
with equation \eqref{eq:entropy_wasserstein_bound},
below that exponential dimension dependence
in this constant is necessary, although the
optimal constant is not known.  
Any sharpening of these constants in future
work could directly be used to 
update the time-uniform bound in Proposition~\ref{prop:smooth}.

Inspired by \cite{weed2018}, we briefly close
this subsection by illustrating how
Proposition~\ref{prop:smooth} can further be used to obtain
sequential bounds for the smoothed differential entropy of $P$, 
$h(P\star \calK_\sigma) = -\int \log(P\star K_\sigma) dP\star \calK_\sigma$,
using the fact that it is Lipschitz with respect
to the $W_1$ metric~\citep{polyanskiy2016}. 
% where $\mu_2(P) = \bbE \norm {X_1}^2$
% and $\mu_2(Q) = \bbE \norm{Y_1}^2$. 
% This readily leads to the following result.
% Invoking Proposition~\ref{prop:smooth},
% we immediately obtain that with probability at least $1-\delta$,
% \begin{corollary}
% \label{cor:entropy}
% Let $P$ be a Borel probability
% measure supported on $[-1,1]^d$. Then, for all $\delta\in (0,1)$,
\begin{corollary}
\label{cor:entropy}
Let $\delta \in (0,1)$, $\sigma > 0$, and $P \in \calP ([-1,1]^d)$. Then, 
\begin{equation*} 
\bbP\left(\forall t \geq 1: |h(P_t \star \calK_\sigma) - h(P \star \calK_\sigma) | \leq 
\frac {3 } {\sigma^2}  \sqrt{\frac d t \Big[\log\ell(\log_2 t) + \log(4/\delta) \Big]}  + \frac{\sqrt d C_{d}}{\sqrt t\sigma^2}\right)\geq 1-\delta.
\end{equation*} 
\end{corollary}
Notice that the mapping $P \mapsto -h(P \star \calK_\sigma)$
is convex, therefore a confidence sequence
for the smoothed differential entropy
could also have been obtained by directly appealing
to Theorem~\ref{thm:main}, assuming a bound on the cumulant generating function 
of $h(P_t\star \calK_\sigma)$ were available.
Beyond this approach 
and the bound of Corollary~\ref{cor:entropy}, we are not aware of other existing sequential
concentration inequalities for this problem.

\subsection{Sequential generalization error
bounds for binary classification}
\label{sec:rademacher}
Theorem~\ref{thm:main} can
be used to derive  generalization error bounds
for classification or regression problems that are valid at  stopping times. 
We illustrate the special case of binary
classification. Let $\calX$ be a topological space, 
$P$ be a Borel probability
distribution over $\calX \times \{-1,1\}$,
and $(X_t, Y_t)_{t=1}^\infty$ a sequence of i.i.d.\ observations from $P$. Let $\calF$ be a collection of Borel-measurable
functions from $\calX$
to $\{-1,1\}$, and define the 
population and empirical classification risks by
$$R(f) = \bbP(f(X) \neq Y), \quad
\text{and} \quad R_t(f) = \frac 1 t \sum_{i=1}^t I(f(X_i) \neq Y_i),\quad f \in \calF.$$
High-probability bounds on the supremum of the empirical process 
$\sup_{f \in \calF} |R(f) - R_t(f)|$
for fixed times $t\geq 1$ are well-studied, 
and can lead to conservative confidence intervals for the generalization 
error $R(\hat f_t)$ of any data-dependent classifier $\hat f_t \in \calF$, such as an empirical risk minimizer, or an approximate one obtained by stochastic optimization. 
Such bounds necessarily depend on the complexity of $\calF$, as measured for instance by
population or empirical Rademacher complexities,
respectively defined by 
$$\calR_t(\calF) = \bbE_{\bepsilon_t,\bX_t} \left[ 
\sup_{f \in \calF} \frac 1 t \left|\sum_{i=1}^t \epsilon_i f(X_i)\right|\right],\quad 
\hat \calR_t(\calF) = \bbE_{\bepsilon_t} \left[ 
\sup_{f \in \calF} \frac 1 t\left| \sum_{i=1}^t \epsilon_i f(X_i)\right|\right].
$$
Here, we denote
$\bX_t=(X_1, \dots, X_t)$ and $\bepsilon_t = (\epsilon_1, \dots, \epsilon_t)$, where $(\epsilon_t)_{t=1}^\infty$ denotes a sequence
of i.i.d. Rademacher random variables (taking values $\pm1$ with equal probability).
We obtain the following bound which is 
uniform both over hypotheses $f \in \calF$ and over time $t\geq 1$.
\begin{corollary}
\label{cor:classification}
Let $\delta \in (0,1)$, and recall that 
$\bar t = \lfloor t/2\rfloor$ for all 
$t \geq 1$. 
\begin{enumerate} 
\item The population Rademacher complexity provides a time-uniform generalization error bound:
$$\bbP\left( \forall t \geq 1:
\sup_{f \in \calF}| R_t(f)-R(f)|
\leq  \calR_{\bar t}(\calF) + 
 2\sqrt{\frac 2 t \Big[ \log \ell(\log_2 t) + \log(1/\delta)\Big]}\right) \geq 1-\delta.$$
\item The empirical Rademacher complexity $(\hat\calR_t(\calF))$ is a reverse
submartingale with respect to $(\calE_t^X)$, 
and we have:
$$\bbP\Big( \forall t \geq 1:
 \calR_{\bar t}(\calF)
 \geq \hat \calR_t(\calF) - 
  4\sqrt{\frac 2 t \Big[ \log \ell(\log_2 t) + \log(1/\delta)\Big]}\Big) \geq 1-\delta.$$
 \end{enumerate}
\end{corollary} 
In particular, if $\tau$ is an arbitrary stopping time and $\widehat f_t$ is an arbitrary data-dependent classifier, then $\bbP\big(| R_\tau(\widehat f_\tau)-R(\widehat f_\tau)|
\leq  \calR_{\bar \tau}(\calF) + 
 2\sqrt{(2/\tau) [ \log \ell(\log_2 \tau) + \log(1/\delta)]}\big) \geq 1-\delta$. We are not aware of other such generalization bounds that hold at stopping times. 
 
 Corollary~\ref{cor:classification} is comparable
to the following
well-known fixed-time bound 
which can be deduced,
for instance, from the proof of Theorem 3.5
of~\cite{mohri2018}:
$$\bbP\left( \sup_{f \in \calF}
|R_t(f)-R(f)| \leq \calR_t(\calF)
 + \sqrt{\log(2/\delta)/2t}\right) \geq 1-\delta.$$
Once again, we observe that
our time-uniform bound only loses
iterated logarithmic factors and small universal constants
in comparison to the above display.
When the population Rademacher complexity $\calR_t(\calF)$
is unavailable in closed form, Corollary~\ref{cor:classification}(ii) 
may be used to provide a time-uniform 
lower bound on this quantity in terms of its
empirical counterpart. 
% implies that the population
% Rademacher complexity
% $\calR_t(\calF)$ may be 
% replaced by its empirical
% counterpart, up to 
% inflating 
% the error term
% $2\sqrt{(2/t) [ \log \ell(\log_2 t) + \log(1/\delta)]}$
% by a universal constant. 
We obtain
this result in 
Appendix~\ref{app:pf_rademacher} by noting
that $\hat\calR_t$
satisfies the leave-one-out
property in equation~\eqref{eq:loo}, 
although it 
cannot easily
be written as the evaluation
of a convex functional
at the empirical measure. We leave open the question
of providing upper confidence sequences on 
$\calR_t(\calF)$, which combined with Corollary~\ref{cor:classification}(i) would
lead to a fully empirical bound for 
the classification~risk.

The proof of Corollary~\ref{cor:classification}(i)
shows that analogous time-uniform concentration 
inequalities can be obtained for 
general suprema of empirical processes
over uniformly bounded function classes, 
up to modifying the expectation bound
$\calR_{\bar t}(\calF)$, which
yields time-uniform inference
for the risk of arbitrary estimators with respect to a bounded loss function,
in terms of their empirical~risk.

\subsection{Sequential estimation of multivariate means}
\label{sec:multivariate_means}
We next show how Theorem~\ref{thm:main}
can be used to derive confidence sequences for 
the mean $\mu$ of a multivariate distribution $P \in \calP(\bbR^d)$. 
In the special case $d=1$, our results
show how our reverse
submartingale techniques
can recover known
confidence
sequences for univariate 
sequential mean estimation
(summarized in Section~\ref{sec:background_confseq}), 
up to constant factors. 

Let $(X_t)_{t=1}^\infty$ be a sequence of i.i.d. random variables
with mean $\mu$, and let $\mu_t = (1/t) \sum_{i=1}^t X_i$. 
We state our bounds in terms of a general norm $\norm\cdot$ on $\bbR^d$, whose  dual norm  is denoted by
$\norm\cdot_\star = \sup_{\norm\lambda = 1} \langle \lambda, \cdot\rangle$. 
Assume further that there exists $\mlambda > 0$ 
and a convex function $\psi: [0,\mlambda) \to \bbR$ such that
\begin{equation}
\label{eq:tail_means}
\sup_{\nu \in \bbS^{d-1}_\star} \log\Big(\bbE_{X \sim P} \left[ \exp\left(\lambda \langle \nu, X-\mu\rangle \right)\right] \Big) \leq \psi(\lambda ), \quad \lambda \in [0,\mlambda),
\end{equation}
where $\bbS^{d-1}_\star = \{x \in \bbR^d: \norm x_\star = 1\}$. For instance, 
when $\psi(\lambda) = \lambda^2\sigma^2/2$, 
the above definition reduces to that of a
$(\sigma^2,\mlambda^{-1})$-sub-exponential random vector
given in~\cite{vershynin2018} when $\mlambda < \infty$, 
or of a $\sigma^2$-sub-Gaussian random vector when
$\mlambda = \infty$. 
Finally, for any $\gamma > 0$, let $N_\gamma$ denote the $\gamma$-covering number~\citep{vandervaart1996} 
of $\bbS^{d-1}_\star$ with respect to the norm $\|\cdot\|_\star$.
\begin{corollary}
\label{cor:multivariate_means}
Assume $P$ satisfies the tail assumption~\eqref{eq:tail_means}. Then, for all $\gamma \in [0,1)$,  $ \delta \in (0,1)$, 
$$\bbP\left\{ \exists t \geq 1 : \norm{\mu_t-\mu} \geq \frac 1 {1-\gamma} (\psi^*)^{-1}\left(\frac{\log\ell(\log_2 t) + \log(1/\delta)  + \log N_\gamma}{\lceil t/2\rceil}\right)\right\}\leq \delta.$$
\end{corollary}
% We further assume $\psi$ is CGF-like, in the sense
% that $\psi(0) = \nabla \psi(0) = \mu$, a property
% which can be seen to hold whenever $\psi$ is indeed the CGF of any mean-zero random
% variable. We also write $\psi_\mu(\lambda) = \psi(\lambda) + \lambda \mu$, 
% which is an upper bound on the CGF of the non-centered random variable $X$. 
Above, $\gamma=1/2$ is a reasonable default value.
We first illustrate 
the result of 
Corollary~\ref{cor:multivariate_means}
in the special case when $d=1$
and $P$ is 1-sub-Gaussian.
If the norm $\norm\cdot$ is taken to be the absolute
value, notice that one may choose 
$\gamma = 0$ and $N_\gamma = 2$, thus 
Corollary~\ref{cor:multivariate_means}
 implies
\begin{equation}
\label{eq:stitching_rev_martingale_means}
\forall t \geq 1 : |\mu_t - \mu| \leq    
2 \sqrt{\frac 1 t \Big[\log\ell(\log_2 t) + \log(2/\delta) \Big]},\quad
\text{with probability at least } 1-\delta.
\end{equation}
Equation~\eqref{eq:stitching_rev_martingale_means} is comparable to state-of-the-art
confidence sequences for univariate means, 
summarized in Section~\ref{sec:background_confseq}. 
For example, 
Theorem~1 of~\cite{howard2021}
provides a one-sided
bound for means
of 1-sub-Gaussian random
variables, which together
with a union bound leads to the two-sided
confidence sequence
% $$|\mu_t - \mu| \leq 
% k_1 \sqrt{t \left( s \log \log t + \log\frac{\zeta(s)}{\delta}\right)}$$
\begin{equation}
\label{eq:howard_stitching_means} 
\forall t \geq 1: |\mu_t - \mu| \leq 
k_1 \sqrt{\frac 1 t \Big[\log \ell(\log_2 t) + \log(2/\delta)\Big]},\quad
\text{with probability at least } 1-\delta,
\end{equation}
where $ 
k_1 = \frac{2^{1/4} + 2^{-1/4}}{\sqrt 2} \approx 1.8.$ 
% Equation~\eqref{eq:howard_stitching_means} holds under the condition~$\sum_{k=0}^\infty 1/\ell(k) \leq 1$, which is slightly stronger
% than our condition~\eqref{eq:stitching_fn_main}.
It can be seen from the preceding two displays
that our confidence sequence 
is wider by a mere factor of $ 2 /1.8 \approx 1.1$
compared to that
of~\cite{howard2021}. In fact, Corollary~\ref{cor:multivariate_means} is a  special case of a more general result that can be deduced
from Theorem~\ref{thm:main_better} 
in Appendix~\ref{app:proofs_prop_exp_bounds}, for which there are tuning parameters that we did not optimize here, so the factor 1.1 could presumably be lowered further. More importantly though, our result applies
more generally to means of multivariate distributions, for which we do not
know of any other confidence sequences
in the literature beyond those of~\cite{abbasi-yadkori2011}. The latter paper only has sub-Gaussian  bounds decaying at the $\sqrt{\log t/t}$ rate, instead of our $\sqrt{\log\log t/t}$~rate. 

As a multivariate example, suppose
now that $P$ is $(\sigma^2, \alpha)$-sub-exponential, 
so that $\psi(\lambda) = \lambda^2 \sigma^2/2$ 
with $\mlambda = 1/\alpha$. 
When $\norm\cdot$ is
taken to be the Euclidean norm 
$\norm\cdot_2$,
% using standard expressions for
% the convex conjugate of the sub-exponential
% CGF and the bound $N_\gamma \leq 5^d$, 
one may derive the bound
%  Derivation, don't delete.
% {\color{red}
% $$\psi(\lambda) = \frac{\lambda^2\sigma^2}{2}, 
% \quad \lambda \in [0,1/\alpha).$$
% $$\psi^*(u)  = \begin{cases} 
% \frac{u^2}{2\sigma^2}, & 0 \leq u < \frac{\sigma^2}{\alpha} \\
% \frac u \alpha - \frac {\sigma^2}{2\alpha^2}, &
% \frac{\sigma^2}{\alpha} \leq u 
% \end{cases},\quad 
% (\psi^*)^{-1}(x)  = \begin{cases} 
% \sqrt{2\sigma^2x}, & 0 \leq x < \frac{\sigma^2}{2\alpha^2} \\
% x\alpha + \frac{\sigma^2}{2\alpha}, &
% \frac{\sigma^2}{2\alpha^2} \leq x.
% \end{cases}
% $$}
$$\forall t \geq 1 : \norm{\mu_t-\mu}_2 \leq 2\begin{cases} 
\sqrt{2\sigma^2\gamma_t}, & t:0 \leq \gamma_t < \frac{\sigma^2}{2\alpha^2}, \\
\gamma_t\alpha + \frac{\sigma^2}{2\alpha}, &
t: \frac{\sigma^2}{2\alpha^2} \leq \gamma_t,
\end{cases},~~\text{with probability at least }  1-\delta,$$
% \begin{align*}
%     \exists t \geq 1 : \norm{\mu_t-\mu}_2 \leq 2 
% &\sqrt{\frac {\sigma^2} t \Big[\log\ell(\log_2 t) + \log(1/\delta) + d\log 5\Big]} \\ 
% &\qquad\qquad \vee \left\{ \frac {4\alpha} t \Big[\log\ell(\log_2 t) + \log(1/\delta) + d\log 5\Big] + \frac{\sigma^2}{2\alpha}\right\}I(\alpha > 0),
% \end{align*}
where $\gamma_t = (2/t) [ \log\ell(\log_2 t) + \log(1/\delta) + d\log 5].$
In particular, we recover the 
optimal dependence on both $d$ and $t$
from the fixed-time setting, up to iterated
logarithmic factors.

\subsection{An upper law of the iterated logarithm
for sub-Gaussian divergences}
\label{sec:lil}
% Thus far we have derived concentration inequalities
% for a wide range of divergences, 
% primarily
% motivated by their use as confidence sequences.
We now show how our finite-sample
results
can be used to derive an
asymptotic statement which 
mirrors the classical
(upper) law of the iterated logarithm (LIL; \cite{stout1970})
for sums of i.i.d. random variables. 
\begin{corollary}
\label{cor:lil}
Let $D$ be a convex divergence such that
$D(P_t\|P)$ is $(\sigma^2/t)$-sub-Gaussian
for all $t \geq 1$
and some $\sigma > 0$. Assume $\bbE D(P_t\|P)=o(\sqrt{(\log\log t)/t})$. 
Then, 
$$\limsup_{t \to \infty} \frac{t D(P_t\|P)}
{\sqrt{2 t \sigma^2 \log\log t}} \leq  1, \quad \text{a.s.}$$
\end{corollary}
Corollary
\ref{cor:lil} establishes an upper LIL
for convex divergences admitting the same constant as
the classical LIL, which
states that for any sequence of
mean-zero i.i.d. random variables
$(X_t)_{t=1}^\infty$
admitting finite variance $\sigma^2 > 0$,  
$$\limsup_{t \to \infty} 
\frac{1}{\sqrt{2t\sigma^2\log\log t}}
\sum_{i=1}^t X_i=1, \quad \text{a.s.}$$
Obtaining a matching lower
bound in Corollary \ref{cor:lil}
would, for instance, necessitate 
anti-concentration bounds
on the process $D(P_t\|P)$, and is
therefore beyond the scope of this work. 
The sub-Gaussianity assumption can also likely
be weakened, but given again that our purpose was
not asymptotics, we leave this for future work.
Though  results analogous to Corollary \ref{cor:lil}
have possibly appeared in past literature for various
divergences, we are only
aware of the LILs for the Kolmogorov-Smirnov statistic
derived by \cite{smirnov1944},
for the 1-Wasserstein distance in dimension $d=1$~\citep{delbarrio1999},
and for certain von Mises differentiable
functionals \citep{serfling2009}.

Adaptations of the proofs of 
Corollaries~\ref{cor:dkw},
\ref{cor:mmd}, \ref{cor:tv}, and~\ref{prop:smooth} 
respectively imply that 
Corollary \ref{cor:lil} 
holds when $D$ is taken
to be the Kolmogorov-Smirnov
distance, the Maximum Mean Discrepancy
with bounded kernel, 
the Total Variation distance  
for distributions supported
on a finite set, or the 
Total Variation and 1-Wasserstein distances
smoothed by Gaussian convolution under suitable
tail assumptions on $P$. 
The conditions
of Corollary 
\ref{cor:lil} 
can also be verified for 
the transportation cost 
$W_p^p$,  for any $p \geq 1$ satisfying 
$p > d/2$,
under suitable moment assumptions on $P$
\citep{fournier2015,niles-weed2019}.
To the best of our knowledge, this functional is not known to be von Mises differentiable without 
further assumptions on the connectedness of the support or 
absolute continuity of $P$~(cf. \cite{goldfeld2022a} and \citeauthor{manole2021} 
(\citeyear{manole2021}, Section 5.2)).

\section{Summary}
\label{sec:conclusion}

 Existing approaches to anytime-valid sequential inference, both classical ideas \citep{darling1967} and modern nonparametric treatments~\citep{howard2020}, are centered around the identification or design of certain nonnegative supermartingales, coupled with a maximal inequality due to~\cite{ville1939}. A primary technical contribution of this current work is to recognize that these techniques are unfortunately not well-suited to dealing with general convex divergences. Instead, we build on the elegant fact that the empirical measure is a measure-valued reverse martingale with respect to the exchangeable filtration, and utilize it to infer that every convex divergence is a  nonnegative reverse (partially ordered) submartingale. 
 Once this connection is made, the literature uncovers a little-known maximal inequality that is well-suited for our purposes---a time-reversed version of Ville's maximal inequality, together with its analogues for partially ordered processes. These tools act as batch-to-sequential devices: they allow us to port existing fixed-time Chernoff-type bounds that have been developed for a large variety of convex divergences into the sequential setting, at the loss of only iterated logarithm factors. As a result of our modular approach, if the fixed-time bounds have optimal rates in all problem-dependent parameters, then so do our sequential bounds, and if some of those bounds are improved in the future, then those improvements are immediately obtained in the sequential setting.

\subsection*{Acknowledgments}
TM would like to thank Ian Waudby-Smith for 
several conversations related to this work.
TM was supported in part by the 
Natural Sciences and Engineering Research Council of Canada, 
through a PGS D Scholarship. 
AR would like to acknowledge NSF DMS grant 1916320 and an Adobe faculty research award for funding work on this topic.

\bibliographystyle{apalike}
\bibliography{manuscript}

\clearpage 

 \addtocontents{toc}{\protect\setcounter{tocdepth}{1}}

\appendix
 
\section{Additional Lemmas}
\label{app:additional_lemmas}
We begin by recalling McDiarmid's inequality \citep{wainwright2019}.
We say that a map $G:\bbR^t \to \bbR$ satisfies the bounded
differences property with parameters 
$(L_1, \dots, L_t)\in \bbR^t_+$ if for every 
$x_1, \dots, x_t$, $x_1', \dots, x_t' \in \bbR$ and all $i=1, \dots, t$, 
$$|G(x_1, \dots, x_t) - G(x_1, \dots, x_{i-1}, x_i', x_{i+1}, \dots, x_t)| \leq L_i.$$
\begin{theorem}[McDiarmid's Inequality]
\label{thm:mcdiarmid}
Assume $G$ satisfies the bounded differences property with parameters
$(L_1, \dots, L_t)$ and that $X_1, \dots, X_t$ are independent random variables. Then, 
for all $u \geq 0,$
$$\bbP\Big( \big| G(X_1, \dots, X_t) - \bbE G(X_1, \dots, X_t)\big| \geq u \Big)
 \leq 2 \exp\left\{ - \frac{2u^2}{\sum_{i=1}^t L_i^2}\right\}.$$
\end{theorem}

In several of the proofs for Section \ref{sec:applications}, we will show that 
the processes $\Phi(P_t,Q)$ or $\Psi(P_t,Q_s)$  satisfy the bounded
differences property, when viewed as functions of the samples therein. The following standard
result will then imply that these processes are sub-Gaussian (see for instance \cite{rigollet2015}, Lemma 1.5, for a statement
with the exact constants used below).
\begin{lemma}
\label{lem:subg}
Let $P$ be a distribution over $\bbR$ such that for $X \sim P$ and $\sigma > 0$, $\bbE[X]=0$ and
$$\bbP(|X| > u) \leq 2e^{-u^2/2\sigma^2}, \quad u > 0.$$
Then, $P$ is $(8\sigma^2)$-sub-Gaussian.
\end{lemma}

%%%%%%%%%%%%%%%%%%%%%%%%%%%%%%%%%%%%%%%%%%%%%%%%%%%%%
% Proofs Section 2
%%%%%%%%%%%%%%%%%%%%%%%%%%%%%%%%%%%%%%%%%%%%%%%%%%%%
\section{Proofs from Section \ref{sec:background}}
\label{app:proofs_sec2} 

\subsection{Proofs of Theorem~\ref{thm:reverse_ville}}
\label{app:proof_ville}
Theorem~\ref{thm:reverse_ville}
was proven for instance by \cite{lee1990}, Theorem 3, p. 112.
Due to its centrality 
in our work, we provide 
two self-contained proofs of this result
below. 
The first proof follows directly
from Doob's submartingale inequality
(see equation \eqref{eq:forward_doob}).
The second proof is a restatement of Lee's
original proof, which we include for reference in the following subsection.

\textbf{First Proof (via Doob's submartingale inequality).}
For any integer $T \geq t_0$, 
define the process
$S_t = R_{T-t+t_0}$, for all $t_0 \leq t \leq T$,
as well as the forward filtration
$\calG_t = \calF_{T-t+t_0}$.
Since $(R_t)$ is a reverse submartingale, we have
$\bbE[R_t|\calF_{t+1}] \geq R_{t+1},$
whence for all $t_0+1 \leq t \leq T$,
$$\bbE[S_t|\calG_{t-1}] 
 = \bbE[R_{T-t+t_0}|\calF_{T-t+t_0+1}]
 \geq R_{T-t+t_0+1} = S_{t-1}.$$
 It follows that
 $(S_t)_{t=t_0}^T$ forms a 
 forward submartingale with
 respect to $(\calG_t)_{t=t_0}^T$.
 Applying Doob's submartingale inequality, 
 we therefore obtain
$$\bbP(\exists t_0\leq t \leq T: S_t \geq u)
\leq \frac{\bbE[S_T]}{u},$$
for all $u > 0$. 
Equivalently,
$$\bbP(\exists t_0\leq s \leq T: R_s \geq u)
\leq \frac{\bbE[R_{t_0}]}{u}.$$
Notice that the event within
the probability on the left-hand side 
of the above display is monotonically
increasing with $T$, converging to the
event $\{\exists s \geq t_0: R_s \geq u\}$.
Taking $T \to\infty$, we thus have
$$\bbP(\exists s\geq t_0: R_s \geq u)
\leq \frac{\bbE[R_{t_0}]}{u},$$
which proves the claim.
\hfill $\square$

\textbf{Second Proof (first principles).}
Let $T \geq t_0$, $u > 0$, and  define the disjoint sets
\[
A_{t} = \{R_{t} \geq u \} \cap \bigcap_{j = t+1}^{T} \{R_{j} < u\},\quad t=t_0,t_0+1, \dots, T,
\]
where $t$ represents the last time (in the range $t_0,t_0+1, \dots, T$) at which $R_t$ was larger than $u$.
We have,
\begin{align*}
\bbP(\exists t_0 \leq t \leq T: R_{t} \geq u) 
 = \sum_{t=t_0}^{T} \bbP\left( A_{t} \right)
 \leq  \frac 1 u \sum_{t=t_0}^{T}  \int_{A_{t}} R_{t} d\bbP
 \leq \frac 1 u \sum_{t=t_0}^{T}  \int_{A_{t}} \bbE[R_{t_0}|\calF_t] d\bbP,
\end{align*}
where the first inequality follows because $R_t > u$ on $A_t$, 
and the second inequality follows by the reverse submartingale property of 
$(R_t)$. Note that $R_j$ is $\calF_j$-measurable and hence also $\calF_t$-measurable for $t \leq j$ due to the reversed nature of the filtrations. 
Thus, $A_t \in \calF_t$, whence
$\int_{A_{t}} \bbE[R_{t_0}|\calF_t] d\bbP = \int_{A_{t}} R_{t_0} d\bbP$
and we obtain
\begin{align*}
\bbP(\exists t_0 \leq t \leq T: R_{t} \geq u) 
 \leq \frac 1 u \sum_{t=t_0}^{T}  \int_{A_{t}} R_{t_0} d\bbP  
 \leq  \frac 1 u \bbE[R_{t_0}],
\end{align*}
where the last step utilizes the nonnegativity of $R_t$ and the fact that the events $\{A_t\}_{t=t_0}^T$ are disjoint by construction.
The claim now follows as before by taking $T \to \infty$, noting that the right-hand side remains fixed
while the left-hand side is the probability of an increasing sequence of events
whose limit is $\displaystyle \left\{ \exists t \geq t_0: R_{t} \geq u\right\}$.
\qed

\subsection{Proof of Proposition 
\ref{prop:maximal_inequality_moments}}
\label{app:proof_prop_max_inequality_moments}
% We begin by discussing a g maximal inequalities for 
% partially ordered reverse submartingales,
% which will be used in the proof of
% Proposition~\ref{prop:maximal_exponential_bounds}(ii).
  
% Let $(\calF_{ts})_{t,s \geq 1}$ be a reverse filtration. 
% Throughout this section, we say that $(\calF_{ts})$
% satisfies the conditional independence property (CI)
% if for all $t,t',s,s' \geq 1$, 
% \begin{equation} 
% \label{eq:conditional_independence}
% \bbE \big\{ \bbE[ ~\cdot ~|~ \calF_{ts'}] ~\big|~\calF_{t's}\big\} = \bbE\{~\cdot ~|~\calF_{(t\vee t')(s\vee s')}\},
% \end{equation}
% or equivalently, that $\calF_{ts'}$ and $\calF_{t's}$ are conditionally
% independent given $\calF_{(t\vee t')(s\vee s')}$. The equivalence
% of these two statements has been noted, for instance, by 
% \cite{merzbach2003}. 

Let $(\calF_{ts})_{t,s \geq 1}$ be a reverse filtration. 
Under the conditional independence condition
\eqref{eq:conditional_independence}
on $(\calF_{ts})$, 
\cite{christofides1990} establish a maximal inequality
for partially ordered reverse submartingales, which we state below before proving
Proposition~\ref{prop:maximal_inequality_moments}. 
We will require the following notation. 
Let $1 \leq t_0 \leq T,1 \leq s_0 \leq S$, and let $\{C_{ts}:t,s \geq 0\}$ 
be a nondecreasing array of nonnegative numbers. 
Given a reverse submartingale $(R_{ts})_{t,s \geq 1}$
with respect to $(\calF_{ts})$,  
a general bound will be given on  $\bbP(A)$, where
for any $u > 0$,
\begin{equation} 
\label{eq:set_A}
A = \left\{\max_{\substack{t_0 \leq t\leq T\\s_0 \leq s \leq S}}  C_{ts}R_{ts} \geq u\right\}
= \bigcup_{t=t_0}^T \bigcup_{s=s_0}^S A_{ts}, ~~
\text{where} ~~ A_{ts} = \{ R_{ts}C_{ts} \geq u\} \cap \bigcap_{j=t+1}^T \bigcap_{k=s+1}^S \{R_{jk}C_{jk} < u\}.
\end{equation}
This decomposition into the sets $A_{ts}$  
is analogous to that in the proof of Ville's inequality for reverse submartingales (Appendix~\ref{app:proof_ville}). Unlike that result, however, the lack of total ordering
on $\bbN^2$ prevents
the sets $A_{ts}$ from being disjoint; for instance $A_{43}$ and $A_{34}$ could both potentially happen.  
\cite{christofides1990} instead form 
a partition $(B_{ts}^{(1)})_{t,s\geq 1}$
of $A$, defined recursively by the following algorithm.
\begin{align*} 
&\text{Let } D_0 = \emptyset, \  m := 1 \\
&\qquad\text{For } j = t_0 \text{ to } T \\
&\qquad\qquad \text{For } k = s_0 \text{ to } S \\
&\qquad\qquad\qquad B_{jk}^{(1)} := A_{jk} \setminus \bigcup_{l < m} D_l \\
&\qquad\qquad\qquad D_m := A_{jk}, \ m := m+1.
\end{align*}
A second partition $(B_{ts}^{(2)})_{t,s \geq 1}$
is further formed by changing
the order of the for-loops in the above display. 
Specifically,
\begin{align*}
B_{ts}^{(1)} = A_{ts} \setminus \left\{ \left( \bigcup_{j=t_0}^{t-1} \bigcup_{k=s_0}^S A_{jk} \right) \cup \left( \bigcup_{k=s_0}^{s-1} A_{tk} \right) \right\},\ \
B_{ts}^{(2)} = A_{ts} \setminus \left\{ \left( \bigcup_{k=s_0}^{s-1} \bigcup_{j=t_0}^T A_{jk} \right) \cup \left( \bigcup_{j=t_0}^{t-1} A_{js} \right) \right\},
\end{align*} 
with the convention that an empty union is equal to the empty set. Notice that for $j=1,2$, 
the sets $(B_{ts}^{(j)})_{t,s \geq 1}$ are mutually disjoint,
and $\bigcup_{t,s \geq 1} B_{ts}^{(j)} = A$. Further, 
unlike $(A_{ts})$, the sequence $(B_{ts}^{(j)})$
is not adapted to $(\calF_{ts})$, but instead satisfies
$B_{ts}^{(1)} \in \calF_{ts_0}$
and $B_{ts}^{(2)} \in \calF_{t_0s}$.
We are now in a position to state their bound.
\begin{lemma}[\cite{christofides1990}, Corollary 2.9]
\label{lem:maximal_partial_order}
Let $(R_{ts})$ be a nonnegative
reverse submartingale 
with respect to $(\calF_{ts})$,
and assume $(\calF_{ts})$ satisfies the conditional 
independence condition~\eqref{eq:conditional_independence}.
 Furthermore, 
let $\{C_{ts}:t,s \geq 0\}$ be a nondecreasing array of 
nonnegative numbers.
% such that 
% $C_{ts} = 0$ for all $t < t_0$ and $s < s_0$. 
Then, for all $u > 0$,
\begin{align*}
u \bbP\left\{\max_{\substack{t_0 \leq t \leq T \\
s_0 \leq s \leq S}} C_{ts} R_{ts} \geq u\right\}
\leq
 &\left\{\sum_{t=t_0}^T \sum_{s=s_0}^S (C_{ts} - C_{(t-1)s}) \bbE[R_{ts}]
 - \sum_{s=s_0}^S C_{t_0s} \int_{(\bigcup_{t=t_0}^T B_{ts}^{(1)})^c} R_{t_0s}d\bbP \right\} 
\\ &\wedge \left\{\sum_{t=t_0}^T \sum_{s=s_0}^S (C_{ts} - C_{t(s-1)}) \bbE[R_{ts}]
 - \sum_{t=t_0}^T C_{ts_0} \int_{(\bigcup_{s=s_0}^S B_{ts}^{(2)})^c} R_{ts_0}d\bbP \right\}.
\end{align*}
\end{lemma}
In the special case where 
$C_{ts} = I(t\geq t_0, s \geq s_0)$
for all $0 \leq t \leq T, 0 \leq s \leq S$,
 Lemma \ref{lem:maximal_partial_order}
 reduces to the following bound
%  which will be invoked in the proof
%  of Proposition~\ref{prop:maximal_inequality_moments},
\begin{equation} 
\label{eq:christofides_simplification}
\begin{aligned}
\bbP\left\{\max_{\substack{t_0 \leq t \leq T \\
s_0 \leq s \leq S}} R_{ts} \geq u\right\} 
&\leq
\frac 1 u  
\left[ \sum_{s=s_0}^S C_{t_0s} \bbE[R_{t_0s}]
 - \sum_{s=s_0}^S C_{t_0s} \int_{(\bigcup_{t=t_0}^T B_{ts}^{(1)})^c} R_{t_0s} d\bbP \right] \\
  &= \frac 1 u
\sum_{s=s_0}^S \int_{\bigcup_{t=t_0}^T B_{ts}^{(1)}} R_{t_0s} d\bbP.
\end{aligned}
 \end{equation}
This simplification of Lemma~\ref{lem:maximal_partial_order}
turns out to be simple to show, and we 
provide a self-contained proof before using
it to prove
Proposition~\ref{prop:maximal_inequality_moments}
below.

{\bf Proof of Inequality~\eqref{eq:christofides_simplification}.} We have for any $s_0 \leq s \leq S$,
\begin{align*}
u\bbP\left( \bigcup_{t=t_0}^T B_{t s}^{(1)} \right)
 &= u\sum_{t=t_0}^T \bbP(B_{ts}^{(1)}) \\
 &\leq \sum_{t=t_0}^T \int_{B_{ts}^{(1)}} R_{ts}d\bbP  \\
 &\leq \sum_{t=t_0}^T \int_{B_{ts}^{(1)}} \bbE [R_{t_0s} | \calF_{ts}] d\bbP & (\text{By the reverse submartingale property}) \\
 &= \sum_{t=t_0}^T\int_{B_{ts}^{(1)}} \bbE \Big\{ \bbE [R_{t_0s} | \calF_{t_0s}] ~|~ \calF_{t s_0}\Big\} d\bbP & \text{(By the CI condition)}\\
 &= \sum_{t=t_0}^T\int_{B_{ts}^{(1)}} \bbE[R_{t_0s} | \calF_{t_0s} ] d\bbP  & (\text{Since } B_{ts}^{(1)} \in \calF_{ts_0}) \\
 &= \sum_{t=t_0}^T\int_{B_{ts}^{(1)}} R_{t_0s} d\bbP
 = \int_{\bigcup_{t=t_0}^T B_{t s}^{(1)}} R_{t_0s} d\bbP.
  \end{align*} 
The claim follows by taking a summation over $s_0 \leq s \leq S$ on both sides.

Lemma \ref{lem:maximal_partial_order} leads to the following proof of Proposition 
\ref{prop:maximal_inequality_moments}, 
which generalizes Corollary  2.10 of \cite{christofides1990}.

{\bf Proof of Proposition 
\ref{prop:maximal_inequality_moments}.}
Let $T \geq t_0, S \geq s_0$. By inequality~\eqref{eq:christofides_simplification},
% Lemma \ref{lem:maximal_partial_order},
% \begin{align*}
% \bbP\left\{\sup_{\substack{t_0 \leq t \leq T \\
% s_0 \leq s \leq S}} R_{ts} \geq u\right\} 
% &=
% \bbP\left\{\sup_{\substack{t_0 \leq t \leq T \\
% s_0 \leq s \leq S}} C_{ts} R_{ts}^\alpha \geq u^\alpha\right\}\\
% &\leq
% \frac 1 {u^\alpha}  
% \left[ \sum_{s=s_0}^S C_{t_0s} \bbE[R_{ts}^\alpha]
%  - \sum_{s=s_0}^S C_{t_0,s} \int_{(\bigcup_{t=t_0}^T B_{ts}^{(1)})^c} R_{t_0s}^\alpha d\bbP \right] \\
%  \end{align*}
  \begin{align*} 
u^\alpha \bbP\left\{\sup_{\substack{t_0 \leq t \leq T \\
s_0 \leq s \leq S}} R_{ts} \geq u\right\} 
 &\leq 
  \sum_{s=s_0}^S \int_{\bigcup_{t=t_0}^T B_{ts}^{(1)}} R_{t_0s}^\alpha d\bbP  \\
%   \right\} 
% \wedge \left\{\sum_{t=t_0}^T  \int_{(\bigcup_{s=s_0}^S B_{ts}^{(2)})} Y_{ts_0}^\alpha d\bbP \right\}\Bigg]\\
 &\leq \sum_{s=s_0}^S \int_{\bigcup_{t=t_0}^T B_{ts}^{(1)}} \left(\max_{s_0\leq s \leq S} R_{t_0s}^\alpha\right) d\bbP  \\
%  \right\} 
% \wedge \left\{\sum_{t=t_0}^T  \int_{(\bigcup_{s=s_0}^S B_{ts}^{(2)})} 
% \left( \max_{t_0 \leq t \leq T} Y_{ts_0}^\alpha\right) d\bbP \right\}\Bigg]\\
 &= \int_A \left(\max_{s_0\leq s \leq S} R_{t_0s}^\alpha\right) d\bbP 
%  \right\} 
% \wedge \left\{\int_A
% \left( \max_{t_0 \leq t \leq T} Y_{ts_0}^\alpha\right) d\bbP \right\}\Bigg]
\\
 &\leq  \bbE\left( \max_{s_0\leq s \leq S} R_{t_0s}^\alpha\right)  \\
 &\leq \left(\frac \alpha{\alpha-1}\right)^\alpha \bbE[R_{t_0s_0}^\alpha],
\end{align*}
where the last inequality follows from 
\cite{doob1953}, Theorem 3.4, page 317. Taking $T,S \to \infty$ 
on both sides of the above display 
leads to the claim.\qed

\section{Proofs from Section \ref{sec:main}}
\label{app:proofs_sec4}

\subsection{Proofs from Subsection~\ref{sec:lcb}}
\label{app:proofs_prop_exp_bounds}
{\bf Proof of Proposition~\ref{prop:maximal_exponential_bounds}.}
To prove part (i), 
Theorem 
\ref{thm:submart} implies that $(N_t)$ forms a reverse submartingale
with respect to $(\calE_{t}^X)$. 
Furthermore,
the map $x\in\bbR \mapsto \exp\{\lambda x\}$
is convex and monotonic for any
fixed $\lambda \in [0,\mlambda)$, so Jensen's inequality implies that the process
$$L_{t}(\lambda) = \exp(\lambda   N_{t}), \quad t \geq 1,$$
is also a reverse submartingale with 
respect to $(\calE_{t}^X)$. 
By  
Theorem~\ref{thm:reverse_ville}, we 
 obtain for all $u > 0$,
\begin{align}
\nonumber 
\bbP\Big(\exists t \geq t_0: N_t \geq u\Big)
 &\leq \inf_{\lambda \in [0,\lambda_{\text{max}})} \bbP\Big( \exists t\geq t_0 : L_{t}(\lambda) \geq e^{\lambda u}\Big) 
 \nonumber 
 \leq \inf_{\lambda \in [0,\lambda_{\text{max}})} \bbE \left[\exp(-\lambda u)  L_{t_0}(\lambda) \right] \\   
 &\leq \inf_{\lambda \in [0,\lambda_{\text{max}})} \ \exp\big\{-\lambda u + \psi_{t_0}(\lambda)\big\}= \exp\big\{ -\psi^*_{t_0}(u)\big\},
 \label{eq:doob_ineq_pf}
\end{align}
as claimed. 
To prove
Proposition~\ref{prop:maximal_exponential_bounds}(ii),
recall that $\calE_{ts} = \calE_t^X \bigvee \calE_s^Y$ is a $\sigma$-algebra
generated by a union of independent $\sigma$-algebras. It follows that
the filtration $(\calE_{ts})$ satisfies the conditional
independence property \eqref{eq:conditional_independence},
by \cite{cairoli1975}, example (a), page 114 (see also \cite{christofides1990a}).
Furthermore, similarly as in part (i), the process
$$L_{ts}(\lambda,\alpha) = \exp(\lambda M_{ts}/\alpha),\quad t,s \geq 1,$$
is a 
partially ordered reverse submartingale
for any fixed choice of $\lambda \in [0,\lambda_{\text{max}})$
and $\alpha > 1$.
Notice also that $L_{ts}(\lambda,\alpha) \in L^\alpha(\bbP)$, thus we may apply Proposition~\ref{prop:maximal_inequality_moments}
to obtain
\begin{align*}
\bbP(\exists t \geq t_0,s \geq s_0: M_{ts} \geq u)
 &= \inf_{\lambda\in[0,\lambda_{\text{max}})} \bbP\Big(\exists t \geq t_0, s \geq s_0: L_{ts}(\lambda,\alpha)  \geq \exp(\lambda u/\alpha)\Big) \\
 &\leq  \left(\frac{\alpha}{\alpha-1}\right)^\alpha  \inf_{\lambda\in[0,\lambda_{\text{max}})} \exp(-\lambda u) \bbE\left[L_{t_0s_0}^{\alpha}(\lambda,\alpha)\right] \\
 &=  \left(\frac{\alpha}{\alpha-1}\right)^\alpha \inf_{\lambda\in[0,\lambda_{\text{max}})} \exp(-\lambda u) \bbE\left[\exp(\lambda M_{t_0s_0})\right] \\
 &\leq  \left(\frac{\alpha}{\alpha-1}\right)^\alpha \exp(-\psi_{t_0s_0}^*(u)).
\end{align*}
Taking the infimum over $\alpha > 1$ on both sides
of the above display leads to the claim.
\qed

We now turn to proving a more
general version of Theorem~\ref{thm:main}. 
% We will prove the following
% more general
% version of Theorem~\ref{thm:main}. 
Let
 $\eta,\xi > 1$ 
 be fixed constants which
 determine the sizes
 of the geometric epochs
 used in the proofs. 
 Furthermore, given
 $t, s \geq 1$, we use
 the shorthand notation
 $\bar t = \lceil t/\lceil\eta\rceil \rceil$
 and $\bar s = \lceil s/\lceil\xi\rceil
 \rceil$.
 \begin{theorem}
\label{thm:main_better}
Let $\Phi,\Psi$ be convex functionals, and let $\delta \in (0,1)$.
\begin{enumerate}
\item[(i)] (One-Sample) Assume $\psi_t^*$ is invertible for all $t\geq 1$, and 
if $\eta$ is not an integer, 
assume
$(\psi_t^*)^{-1}(\lambda)$  
is a decreasing (resp. increasing) sequence in $t$ 
(resp. $\lambda$). Set
$$\gamma_t = 
(\psi_{\bar t}^*)^{-1}\Big(\log \osfn(\log_\eta  t) + \log(1/\delta)\Big).$$
Assume further that $(\gamma_t)$ is a nonincreasing
sequence. 
% Furthermore, if 
% $\eta < 2$, assume $(\varphi_t^*)^{-1}(\lambda)$ is a decreasing function
% of $t \geq 1$ for any fixed $\lambda \in [0,\mlambda)$. 
Then,
\begin{align*}
\bbP&\big\{\exists t \geq 1: \Phi(P_t) \geq \Phi(P) + \gamma_t \big\} \leq \delta.
\end{align*}
\item[(ii)] (Two-Sample) Assume $\psi_{ts}^*$ is invertible for all $t,s\geq 1$, and if
$\eta$ (resp  $\xi$) is not an integer,
assume $\psi_{ts}^*(\lambda)$  
is decreasing in 
$t$ (resp. in $s$), and increasing in $\lambda$. Set
$$\gamma_{ts} = 
(\psi_{\bar t\bar s}^*)^{-1}\Big(\log g(\log_\eta t + \log_\xi s) + \log(1/\delta)\Big).$$
% Furthermore, if 
% $\eta < 2$ (resp. $\xi < 2$), 
% assume $\gamma_{ts}$ is a decreasing function
% of $s \geq 1$ (resp. $t \geq 1$). 
Assume further that $(\gamma_{ts})$
is a nonincreasing sequence in each of its indices.
Then,
\begin{align*}
\bbP&\left\{\exists t,s \geq 1: \Psi(P_t, Q_s) \geq \Psi(P, Q) + 
\gamma_{ts}\right\}
\leq \delta.
\end{align*}
\end{enumerate}

\end{theorem}

{\bf Proof of Theorems~\ref{thm:main}
and~\ref{thm:main_better}.}
The proofs of claims (i) and (ii) of Theorem~\ref{thm:main_better} 
are similar, thus we only prove (ii).
Theorem~\ref{thm:main}
will then follow by setting $\eta = \xi = 2$. 
Let $u_j = \lceil \eta^j \rceil$ and $v_k = \lceil \xi^k\rceil$ for all $j,k \in \bbN_0$. Since $\gamma_{ts}$ is decreasing
in $t$ and $s$, we have
\begin{align*}
\bbP &\left(\exists t,s \geq 1: M_{ts} \geq \gamma_{ts}\right)  \\ 
 &\leq \bbP\left(\bigcup_{j\in \bbN_0} \bigcup_{k \in \bbN_0} \Big\{\exists t \in \{u_j, \dots, u_{j+1}\}, s \in \{v_k, \dots, v_{k+1}\}: M_{ts}   \geq  \gamma_{ts}\Big\}\right) \\
 &\leq \bbP\left(\bigcup_{j\in \bbN_0} \bigcup_{k \in \bbN_0} \Big\{\exists t \in \{u_j, \dots, u_{j+1}\}, s \in \{v_k, \dots, v_{k+1}\}: M_{ts}   \geq  \gamma_{u_{j+1}v_{k+1}}\Big\}\right).
 \end{align*}
 Now, $u_{j+1}\leq u_j \lceil \eta \rceil$
(resp. $v_{k+1} \leq v_k \lceil \xi \rceil$),
with equality if $\eta$ (resp. $\xi$)
is an integer. 
Therefore, 
by definition of $\bar t,\bar s$, 
and by the fact that
$(\psi_{ts}^*)^{-1}$ is decreasing 
in $t$ (resp. $s$) when
$\eta$ (resp. $\xi)$ is not an integer, 
we have
\begin{align*}
\gamma_{u_{j+1}v_{k+1}}
 &\geq (\psi_{u_jv_k}^*)^{-1}\Big(\log g\big(\log_\eta(u_{j+1}) + \log_\xi(v_{k+1})\big) + \log(1/\delta)\Big).
\end{align*}
Since $(\psi^*_{ts})^{-1}(\lambda)$ is increasing in $\lambda$ when
$\eta$ or $\xi$ are not integers, we deduce 
\begin{align*}
\gamma_{u_{j+1}v_{k+1}}
 &\geq (\psi_{u_jv_k}^*)^{-1}\Big(\log g(j+k+2) + \log(1/\delta)\Big).
\end{align*}
 Applying a union bound together
 with Proposition~\ref{prop:maximal_exponential_bounds}
 then leads to
 \begin{align*} 
\bbP   \left(\exists t,s\geq 1: M_{ts} \geq \gamma_{ts}\right)     
 &\leq e\sum_{j=0}^\infty \sum_{k=0}^\infty \exp\left\{-\psi_{u_jv_k}^*\Big(
  (\psi_{u_jv_k}^*)^{-1}\big(\log g(j+k+2) + \log(1/\delta)\big)\Big) \right\}\\
 &\leq e\sum_{j=0}^\infty \sum_{k=0}^\infty \exp\left\{-\big(\log g(j+k+2) + \log(1/\delta)\big)  \right\}\\
&= \delta \sum_{j=0}^\infty\sum_{k=0}^\infty\frac e {g(j+k+2)}  
 \leq \delta,\end{align*}

% Therefore, 
% \begin{equation} 
% \label{eq:pf_main_summable_bound}
% \exp\big\{-\psi_{u_jv_k}^*(\gamma_{u_{j+1}v_{k+1}})\big\}
%  \leq \exp\left\{
%  -\log g(j+k+2) + \log(1/\delta)\Big)\right\}
%  =  \frac \delta {g(j+k+2)},
%  \end{equation}
%  so that
%  \begin{align*}
%  \bbP &\left(\exists t \geq 1,s\geq 1: M_{ts} \geq \gamma_{ts}\right)
%  \leq \delta  \sum_{j,k=0}^\infty 
% \frac 1 {g(j+k+2)} 
% =\delta \sum_{j,k=1}^\infty \frac 1 {g(j+k)}\leq \delta,
%  \end{align*}
 as claimed.
%  
% $$\psi_{u_kv_j}^*(\gamma_{u_{k+1}v_{j+1}})
% \geq 
% \psi_{u_{k}v_{j}}^*(\gamma_{u_{k}v_{j}})
%  \geq   
% (\psi_{ts}^*)^{-1}\Big(\log g(\log_\eta t + \log_\xi s) + \log(1/\delta)\Big).$$
% \begin{align*}
%  &\leq \sum_{k=0}^\infty \sum_{j=0}^\infty \exp(-\psi_{u_kv_j}^*(\gamma_{u_{k+1}v_{j+1}})) \\
%  &= \delta \sum_{k=0}^\infty \sum_{j=0}^\infty \frac 1 {h(\log_\eta(u_k) + \log_\xi(v_j))} \\
%  &\leq \delta \sum_{j=0}^\infty \sum_{k=0}^\infty \frac 1 {h(k+j)}  \leq \delta.
% \end{align*}
% The claim follows.
% \hfill $\square$
\qed

{\bf Proof of Corollary \ref{cor:sub_gaussian_bound}.}
We only prove the claim for $(M_{ts})_{t,s=1}^\infty$, as the proof for $(N_t)_{t=1}^\infty$ is similar. 
Since $M_{ts}$ is $\sigma_{ts}^2$-sub-Gaussian, 
we have for all $\lambda \in \bbR_+$,
\begin{align*}
\bbE\left\{\exp\left(\lambda M_{ts}\right)\right\}
 &= \bbE \left\{ \exp\left[\lambda (M_{ts} - \bbE (M_{ts}))\right]\right\} \exp\left\{ \lambda \bbE (M_{ts})\right\}  
 \leq \exp\left\{\frac{\lambda^2 \sigma_{ts}^2}{2}\right\} \exp\left\{\lambda \bbE(M_{ts})\right\},
\end{align*}
whence, an upper bound on the CGF of $M_{ts}$ is given by
$\psi_{ts}(\lambda) = \frac{\lambda^2 \sigma_{ts}^2}{2} + \lambda  \bbE(M_{ts})$. Thus,
for any $x \geq \bbE(M_{ts})$ and any $\gamma \in \bbR_+$,
$$\psi_{ts}^*(x) = \frac{(x -  \bbE(M_{ts}))^2\sigma_{ts}^{-2}}{2},\quad  (\psi_{ts}^*)^{-1}(\gamma) = \bbE(M_{ts}) + \sqrt{2\gamma \sigma_{ts}^2}.$$
The claim now follows from Theorem~\ref{thm:main}.
\hfill $\square$

\subsection{Proofs from Subsection~\ref{sec:ucb}}
\label{app:proofs_subsec3_2}
We state and prove the following stronger version of Corollary~\ref{cor:ucb_bounded}.
\begin{proposition}
\label{prop:ucb}
Assume that the processes $(N_t)$
and $(M_{ts})$ satisfy the conditions of 
Proposition~\ref{prop:affine_rev_martingale}, and fix $\delta \in (0,1)$. 
Assume further that there exists
$\mlambda > 0$ and convex functions
$\psi_{\Phi}, \psi_{\Psi},
\phi_{\Psi}$ such that for all $\lambda \in [0,\mlambda)$,
$$
\sup_{f \in \calF_\Phi} 
\log\left\{\bbE\left[e^{\lambda (f(X)-\bbE f(X))}\right]\right\} \leq  \psi_{\Phi}(\lambda),~~
$$
$$
\sup_{f \in \calF_\Psi}\log\left\{\bbE\left[e^{\lambda (f(X)-\bbE f(X))}\right]\right\} \leq \psi_{\Psi}(\lambda) ,~~
\sup_{g \in \calG_\Psi} 
\log\left\{\bbE\left[e^{\lambda (g(Y)-\bbE g(Y))}\right]\right\} \leq \phi_{\Psi}(\lambda).
$$
Assume further that the Legendre-Fenchel transforms
$\psi_\Phi^*,\psi_\Psi^*,\phi_\Psi^*$ are invertible,
with nondecreasing inverses.
\begin{enumerate} 
 \item[(i)] (One-Sample) We have
\begin{align*}
\bbP&\left\{\exists t \geq 1: \Phi(P_t) \leq  \Phi(P) - (\psi_\Phi^*)^{-1}\left( \frac 2 t \Big[ \log \ell (\log_2 t) + \log(2/\delta)\Big] \right) \right\} \leq \delta/2.
\end{align*}
\item[(ii)] (Two-Sample)  Define
$$\kappa_t^X = (\psi_\Psi^*)^{-1}\left( \frac 2 t \left[ \log \ell (\log_2 t) + \log\left(\frac 4 \delta\right)\right] \right), ~~
  \kappa_s^Y = (\phi_\Psi^*)^{-1}\left( \frac 2 s \left[ \log \ell (\log_2 s) + \log\left(\frac 4 \delta\right)\right] \right).$$
Then, we have
\begin{align}
\label{eq:two_sample_ucb}
\bbP\big(\exists t,s \geq 1: \Psi(P_t, Q_s) \leq  \Psi(P, Q)- \kappa_t^X - \kappa_s^Y \big)  \leq \delta/2.
\end{align}
\end{enumerate}
\end{proposition}
% For example, when $\psi_\Psi(\lambda)$ and $\phi_\Psi(\lambda)$ 
% are taken to scale as $\lambda^2\sigma^2/2$ for all $\lambda \geq 0$
% and some $\sigma > 0$, implying sub-Gaussian tails, 
% it can be verified that $\kappa_t^X = O(\sqrt{\log\ell(\log_2 t)/t})$
% and $\kappa_s^Y =O(\sqrt{\log\ell(\log_2 s)/s})$. 
% The proof appears in Appendix~\ref{app:proofs_subsec3_2}
When it exists, the cumulant generating function 
of any mean-zero random variable $Z$
scales quadratically near zero---specifically, it is easy to check that $\lim_{\lambda \to 0} \log (\mathbb{E} \exp(\lambda Z)) / (\lambda^2/2) = \text{Var}(Z)$. Thus, the inverse
of its Legendre-Fenchel transform typically scales
as the square root function near zero. 
The upper confidence sequences in Proposition~\ref{prop:ucb} thus typically
scale at the parametric rate up to a necessary iterated logarithmic factor.

{\bf Proof of Proposition~\ref{prop:ucb}.}
To prove part (i), 
define for this proof only, 
$$\eta_t=  (\psi_\Phi^*)^{-1}\left(\frac 2 t \Big[\log\osfn
(\log_2 t) + \log(2/\delta)\Big]\right).$$
Furthermore, let $u_j = 2^j$ for all $j \in \bbN_0$.
By Proposition~\ref{prop:affine_rev_martingale}, $(R_t)$ minorizes $(N_t)$ thus
\begin{align} 
\label{eq:pf_ucb_step}
\nonumber 
\bbP&(\exists t \geq 1: \Phi(P_t) \leq 
\Phi(P)-\eta_t)
 \leq 
\bbP\left(\exists t \geq 1: -R_t \geq 
\eta_t\right)\\
 \nonumber 
 &\leq \bbP\left(\bigcup_{j\in \bbN_0}   \Big\{\exists t \in \{u_j, \dots, u_{j+1}\}: -R_t   \geq 
\eta_t\Big\}\right) \\
&\leq \bbP\left(\bigcup_{j\in \bbN_0}   \left\{\exists t \in \{u_j, \dots, u_{j+1}\}: -R_t   \geq  (\psi_\Phi^*)^{-1}\left(\frac 2 {u_{j+1}} \left[\log \ell(j+1) + \log(2/\delta)\right]\right)\right\}\right),
\end{align}
where on the last line, we used the fact
that $(\psi_\Phi^*)^{-1}$ is nondecreasing.
Now, $(R_t)$ is a reverse martingale
with respect to $(\calE_{t}^X)$, whence
$(\exp(-\lambda R_t))_{t=1}^\infty$
is a reverse submartingale for any fixed 
$\lambda \in [0,t_0\mlambda)$. 
Applying  
Theorem~\ref{thm:reverse_ville} similarly
as in the proof of Proposition~\ref{prop:maximal_exponential_bounds}, we 
therefore obtain for all $u > 0$ and $t_0 \geq 1$,
\begin{align*}
\bbP\Big(\exists t \geq t_0: -R_t \geq u\Big)
 &= \inf_{\lambda \in [0,t_0\lambda_{\text{max}})} \bbP\Big( \exists t\geq t_0 : \exp(-\lambda R_{t}) \geq e^{\lambda u}\Big) \\ 
 &\leq \inf_{\lambda \in [0,t_0\lambda_{\text{max}})}\exp(-\lambda u) \bbE  \left[ \exp(- \lambda R_{t_0}) \right] \\ 
 &\leq \inf_{\lambda \in [0,t_0\lambda_{\text{max}})}\exp\big\{-\lambda u
 + t_0\psi_\Phi(\lambda/t_0)\big\} \\
 &\leq \inf_{\lambda \in [0,t_0\lambda_{\text{max}})} \ \exp\big\{-t_0[(\lambda/t_0) u - \psi_{\Phi}(\lambda/t_0)]\big\}\\
 &= \inf_{\lambda \in [0,\lambda_{\text{max}})} \ \exp\big\{-t_0[\lambda u - \psi_{\Phi}(\lambda)]\big\}
 = \exp\big\{ -t_0\psi_\Phi^* (u)\big\}.
\end{align*}
Returning to equation~\eqref{eq:pf_ucb_step},
we deduce
\begin{align*}
\bbP(\exists t \geq 1 &: \Phi(P_t) \leq 
\Phi(P)-\eta_t) 
 \leq \sum_{j=0}^\infty 
    \exp\left\{ -\frac{2u_{j}}{u_{j+1}} 
    \left[\log \ell(j+1) + \log(2/\delta)\right]\right\}\\
 &\leq \sum_{j=0}^\infty \exp\left\{ - 
 \left[\log \ell(j+1) + \log(2/\delta)\right]\right\} 
 =  \sum_{j=0}^\infty \frac {\delta/2} {\ell(j+1)}
 \leq \frac \delta 2.
 \end{align*}
The proof of claim (i) follows. The proof
of part (ii) follows by a similar
probability bound for each of $-R_t^X$
and $-R_s^Y$ at level $\delta/4$, combined with a union bound.
\qed

\textbf{Proof of Corollary~\ref{cor:ucb_bounded}.} 
By Hoeffding's Lemma, we may take %(\cite{wainwright2019}, Example 2.4),
$\psi_\Phi(\lambda) = \lambda^2 B^2/8$
for all $\lambda \in \bbR_+$,
thus, 
$(\psi_\Phi^*)^{-1}(\lambda) = %(\phi_\Psi^*)^{-1}(\lambda) = 
\sqrt{B^2 \lambda/2}$, and similarly for the
two-sample case. 
The claim follows directly from Proposition~\ref{prop:ucb}.
\qed

\subsection{Proofs from Subsection~\ref{sec:choice_filtrations}}
\label{app:proofs_subsec3_3}
{\bf Proof of Proposition~\ref{prop:stopping}.}
Assume first that 
$\bbP(\bigcup_{t,s=1}^\infty A_{ts})\leq \delta$. 
Let $\eta=(T,S)$ be any random time. Then
\begin{align*}
A_{TS}
 &= \left(\bigcup_{t=1}^\infty \bigcup_{s=1}^\infty 
A_{ts} \cap \{\eta=(t,s)\} \right)   \\
  &\quad\cup \left( \bigcup_{t=1}^\infty A_{t\infty}\cap\{\eta=(t, \infty)\}\right)
  \cup \left(\bigcup_{s=1}^\infty A_{\infty s}\cap\{\eta=(\infty, s)\}\right) 
 \cup \left( A_{\infty\infty} \cap \{\eta = (\infty,\infty)\}\right).
\end{align*}
Since  
$A_{\infty\infty},A_{t\infty},A_{\infty s}
\subseteq\bigcup_{t,s=1}^\infty A_{ts}$,
we deduce that $A_{TS} \subseteq  \bigcup_{t,s=1}^\infty A_{ts}$, implying
that $\bbP(A_{TS}) \leq \delta$. 
Thus (i) implies (iii). It is also clear
that (iii) implies (ii), thus it remains
to prove that (ii) implies (i). 
To this end, assume $\bbP(A_{\tau\sigma}) \leq \delta$
for any stopping time $(\tau,\sigma)$.
For any $\omega \in \Omega$, let 
$$I(\omega) = \left\{(t,s) \in \bbN^2: 
\omega \in A_{ts} \text{ and } \omega \not\in A_{t's'}, 
\forall (t',s') \in \bbN^2
,(t',s') < (t,s)\right\}.
$$
% $(\tau(\omega),\sigma(\omega))$
% denote any pair $(t,s) \in \widebar\bbN^2$ such that
% $\omega \in A_{ts}$ and $\omega \not\in A_{t's'}$
% for all $(t',s') \in \bbN^2$
% such that $t' < t$ and $s' < s$
% (note that $t$ and $s$ may equal the symbol
% $\infty$). 
% If no such pair exists, let $(\tau(\omega), \sigma(\omega))
% =(\infty,\infty)$. 
We may then
define 
$$(\tau(\omega), \sigma(\omega))
 = \begin{cases}
(\infty,\infty), & I(\omega) = \emptyset \\ 
\displaystyle \argmin_{(t,s) \in I(\omega)} t, & |I(\omega)| \geq 1  
 \end{cases}.$$
The minimizer in the above display is unique and unambiguous, because when $I(\omega)$
has cardinality greater or equal to 2, any of its distinct  
elements $(t,s)$ and $(t',s')$ must have $t\neq t'$ and $s \neq s'$ by construction; for example,  $(t,s)$ and $(t,s')$ cannot both be elements of $I(\omega)$ for $s\neq s'$.
Notice that $(\tau,\sigma)$ is a stopping
time with respect to $(\calF_{ts})$ because $A_{t's'} \in \calF_{ts}$ for all $(t',s')\leq (t,s)$. 
Furthermore, its definition
guarantees $\bigcup_{t,s=1}^\infty A_{ts} \subseteq A_{\tau\sigma}$.
We deduce by assumption that
$\bbP(\bigcup_{t,s=1}^\infty A_{ts}) \leq \delta,$ as claimed. 
\qed

We note that our precise definitions of the events 
$A_{t\infty}$ and $A_{\infty s}$, for $t,s \geq 1$, was not crucial in the preceding argument.

 %%%%%%%%%%%%%%%%%%%%% DKW

\section{Proofs from Section \ref{sec:applications}}
\label{app:proofs_sec5}

\subsection{Proofs from Subsection~\ref{sec:dkw}}
\label{app:proof_dkw}

{\bf Proof of Corollary \ref{cor:dkw}.}
%In what follows, we let $\norm\cdot_\infty$ denote the $\ell_\infty$ norm.
By the DKW inequality, for all $u > 0$, we have
$$\bbP \Big( \norm{F_t - F}_\infty \geq u\Big) \leq 2 e^{-2t u^2}.$$
Therefore, 
\begin{align*}
\bbE[\norm{F_t-F}_\infty]
= \int_0^\infty \bbP(\norm{F_t-F}_\infty \geq u) du 
\leq 2 \int_0^\infty e^{-2tu^2}du 
= \sqrt{\frac \pi {2t}}.
\end{align*}
Furthermore, it is 
a straightforward observation
that the map
$$(x_1, \dots, x_t) \in \bbR^t
\mapsto \sup_{x \in \bbR}\left|\frac 1 t \sum_{i=1}^t I(x_i \leq x) - F(x)\right|$$
satisfies the bounded differences
property with parameters
$L_1=\dots=L_t=1/t$. 
McDiarmid's inequality (Theorem~\ref{thm:mcdiarmid})
therefore implies the bound
$$\bbP\big( \big| \norm{F_t-F}_\infty-\bbE \norm{F_t-F}_\infty\big| \geq u \big) 
\leq 2 e^{-2tu^2}, \quad u > 0.$$
Lemma \ref{lem:subg} then implies that 
$\norm{F_t-F}_\infty$ is a $(2/t)$-sub-Gaussian random variable. 

Now, notice that $D_\calJ(P_t\| P) = \norm{F_t-F}_\infty$ is an IPM over the class $\calJ = \{(-\infty, x]: x \in \bbR\}$, and in particular $D_\calJ$ is convex by Lemma \ref{lem:convex}.
We may therefore apply Corollary \ref{cor:sub_gaussian_bound} to the process
$N_t =D_\calJ(P_t\| P)$, $t \geq 1$,
to obtain the claim.
\hfill $\square$

{\bf Proof of Corollary \ref{cor:ks_test}.}
By the triangle inequality, 
$$M_{ts}:=\norm{F_t-G_s}_\infty - \norm{F-G}_\infty
\leq \norm{G_s - G}_\infty + \norm{F_t - F}_\infty,$$
so that $\bbE M_{ts} \leq \sqrt{\pi/2t} + 
\sqrt{\pi/2s}$, 
by the same argument as in 
the proof of Corollary~\ref{cor:dkw}.
Furthermore, the map
$$(x_1, \dots, x_t, y_1, \dots, y_s)
 \mapsto \sup_{x \in \bbR} \left|
 \frac 1 t \sum_{i=1}^t I(x_i \leq x) - \frac 1 s 
 \sum_{i=1}^s I(y_i \leq x)\right|,$$
satisfies the bounded differences
property with parameters
$L_1=\dots=L_t=1/t$ and 
$L_{t+1}=\dots=L_{t+s}=1/s$. 
Therefore, McDiarmid's inequality implies
% $$\bbP\big(\big|\norm{F_t-G_s}_\infty
% - \bbE \norm{F_t-G_s}_\infty \big| \geq u\big)
% \leq 2\exp(-2tsu^2/(s+t)),\quad u>0$$
% so that 
$$\bbP\big(|M_{ts}
- \bbE M_{ts}| \geq u\big)
\leq 2\exp(-2tsu^2/(s+t)),\quad u>0.$$
% Therefore, by the DKW inequality, 
% $$\bbP(M_ts \geq u)
%  \leq \bbP(\norm{F_t-F}_\infty \geq u/2)
%   + \bbP(\norm{G_s-G}_\infty \geq u)
%   \leq 2\exp(-2tu^2) + 2\exp(-2su^2)
%   \leq 4\exp(-2(t\wedge s)u^2).$$
It now follows similarly
as in the proof of Corollary~\ref{cor:dkw} 
that 
$M_{ts}$ is sub-Gaussian with parameter
$2(t+s)/ts$. Applying
Corollary \ref{cor:sub_gaussian_bound}
leads to the bound
$\bbP(\exists t,s \geq 1: M_{ts} \geq \gamma_{ts})
\leq \delta/2$.
Furthermore, from Corollary~\ref{cor:ucb_bounded}, 
$\bbP(\exists t,s \geq 1: M_{ts} \leq -\kappa_{ts}) \leq \delta/2.$
Applying a union bound leads to the claim. 

To prove the validity of the 
test, notice simply that under
the null $F = G$, the aforementioned
bound reduces to 
$\bbP(\exists t,s \geq 1 :
\norm{F_t-G_s}_\infty \geq \gamma_{ts}) \leq \delta/2$.
\qed

%%%%%%%%%%%%%%%%% MMD

\subsection{Proofs from Subsection \ref{sec:mmd}}
\label{app:proofs_mmd}

\textbf{Proof of Corollary \ref{cor:mmd}}
The proof of Theorem 7 (equation (16)) of \cite{gretton2012}
yields the expectation bound
$$\bbE \big| \text{MMD}(P_t, Q_s) - \text{MMD}(P,Q)\big| \leq 2 \Big[(B/t)^{1/2} + (B/s)^{1/2}\Big].$$
Further, the following concentration bound follows
from equation~(15) of \cite{gretton2012}:
$$\bbP\big(\big|\text{MMD}(P_t, Q_s) - \bbE[\text{MMD}(P,Q)] \big| \geq u\big) \leq 2 \exp\left\{-\frac{tsu^2}{2B(t+s)}\right\}, \quad u > 0.$$
Therefore, Lemma \ref{lem:subg} implies that $\text{MMD}(P_t, Q_s)$ is $\frac{8B(t+s)}{ts}$-sub-Gaussian. 
Finally, $\mmd$ is an IPM by its definition, 
and is therefore convex by Lemma \ref{lem:convex}. 
Combining these facts with Corollary~\ref{cor:sub_gaussian_bound}, applied to the 
process $M_{ts} = \mmd(P_t,Q_s) - \mmd(P,Q)$,
leads to the bound
$$\bbP(\exists t,s \geq 1: \mmd(P_t,Q_s) - \mmd(P,Q) \geq  \gamma_{ts}) \leq \delta/2.$$
To obtain an upper confidence sequence, notice
that the set $\calJ = \{f \in \calH: \norm f_\calH \leq 1\}$ consists of functions taking values in the interval $[-\sqrt B, \sqrt B]$. Indeed, for all $f \in \calJ$,
$x \in \bbR^d$,
$$|f(x)| = |\langle f, K(x,\cdot)\rangle_{\calH}|
 \leq \norm f_\calH \norm{K(x,\cdot)}_\calH 
 \leq \sqrt{K(x,x)} \leq \sqrt B.$$
Applying Corollary~\ref{cor:ucb_bounded},
we thus have
$$\bbP(\exists t,s \geq 1: \mmd(P_t,Q_s) - \mmd(P,Q)  \leq - 2\sqrt{B}\kappa_{ts}) \leq \delta/2.$$
Applying a union bound leads to the claim.
\qed

Thus far we have studied the plug-in estimator $\mmd^2(P_t, Q_s)$, which is known to be a biased estimator of $\mmd^2(P,Q)$. 
The following unbiased estimator, which we state only in the case $t=s$, is widely-used and
is obtained by replacing the V-statistic in 
\eqref{eq:one_sample_mmd_v_stat} by the following 
U-Statistic
\begin{equation}
\label{eq:unbiased_mmd}
\hat M_t^2 = \frac 1 {t(t-1)} \sum_{i\neq j} \tilde J(Z_i, Z_j), 
\end{equation}
where $\tilde J((x,y), (x',y')) = K(x,x') + K(y,y') - K(x,y') - K(x',y)$.
The process $\hat M_t$ does not admit a simple characterization as a convex functional of the empirical measures
$P_t$ and $Q_t$, thus Theorem~\ref{thm:main} 
and Proposition~\ref{prop:ucb} cannot be directly applied. 
U-Statistics are, however, known to be reverse martingales, as discussed in Section~\ref{sec:mmd}, 
implying that $\hat M_t^2$ is a reverse martingale. 
While this does not imply that 
$\hat M_t-\mmd(P,Q)$ is a reverse submartingale, 
Theorem~\ref{thm:main} can be applied
directly to the mean-zero process 
$\hat M_t^2 - \mmd^2(P,Q)$.
\begin{proposition}
\label{prop:mmd_ustat}
Under the same conditions as Corollary \ref{cor:mmd}, we have for all $\delta \in (0,1)$,
$$\bbP \left( \exists t \geq 1 : \hat M_t^2 \geq \mmd^2(P,Q) + 16B\sqrt{\frac{1}{t-1}\Big[\log \osfn(\log_2 t)  + \log(1/\delta)\Big]}\right) \leq \delta.$$
\end{proposition}
% We illustrate this approach in Proposition 
% \ref{prop:mmd_ustat} 
% (Appendix \ref{app:proofs_sec5}), though
Unlike equation \eqref{eq:mmd_adaptive}, 
Proposition~\ref{prop:mmd_ustat}
does not lead to a confidence sequence
for $\mmd^2(P,Q)$ scaling at the
rate $O(\log\log t/t)$ when $P=Q$. We therefore recommend
the use of the plug-in estimator
$\mmd(P_t, Q_s)$ and
Corollary \ref{cor:mmd}
when a confidence sequence is needed in practice.

\textbf{Proof of Proposition~\ref{prop:mmd_ustat}.}
% Notice that the map
% $$G:(z_1, \dots, z_t) \mapsto \frac 1 {t(t-1)} \sum_{i\neq j} \tilde J(z_i, z_j)$$
% satisfies for all $z_1, \dots, z_n, z_1', \dots, z_n' \in \bbR^d$, and for all $k=1, \dots, n$,
% $$|G(z_1, \dots, z_t) - G(z_1, \dots, z_{k-1}, z_k', z_{k+1}, \dots, z_t)|
%  \leq \frac 1 {t(t-1)} \sum_{i\neq k}^t \Big| \tilde J(z_i,z_k) - \tilde J(z_i', z_k)\Big| 
%  \leq \frac{2B}{t}.
%  $$
% It follows that $G$ satisfies the bounded differences property, whence
% $$\bbP\Big( \big|G(Z_1, \dots, Z_t) - \bbE G(Z_1, \dots, Z_t)\big| \geq u\Big) \leq \exp(-tu^2/2B^2), \quad u > 0.$$
% Recalling equation \eqref{eq:unbiased_mmd} and the unbiasedness of $\hat M_t^2$, we arrive at the following bound,
% $$\bbP\Big( \hat M_t^2 - \mmd^2(P,Q) \geq u\Big) \leq \exp(-tu^2/2B^2), \quad u > 0.$$
By Theorem~10 of~\cite{gretton2012}
and Lemma \ref{lem:subg}, one can infer similarly
as in the proof of Corollary~\ref{cor:mmd} that
$\hat M_t^2$ is $(32B^2/\underline{t})$-sub-Gaussian, 
where $\underline t = \lfloor t/2\rfloor$.
The claim then follows by the same proof technique as Theorem 
\ref{thm:main}  and Corollary \ref{cor:sub_gaussian_bound}, 
using the fact that $\hat M_t^2$ is a 
reverse martingale, and using the inequality $\lceil\lfloor{t/2}\rfloor /2\rceil\geq (t-1)/4$.
\qed 

 We close this section with a proof of Proposition~\ref{prop:v_stats}. 

{\bf Proof of Proposition~\ref{prop:v_stats}.}
Let $\nu \in \calP(\calX)$ be any fixed reference measure. 
By Mercer's Theorem
(see for instance~\cite{christmann2008}, Theorem 4.49), 
a continuous, symmetric, and positive definite kernel $h$ admits the representation
$$h(x,y) = \sum_{i=0}^\infty \lambda_i \psi_i(x)\psi_i(y),
\quad x,y \in \calX,$$
where $(\lambda_i)_{i\geq 0} \subseteq \ell^2 \equiv \ell^2(\bbN_0)$
is the sequence of eigenvalues corresponding to 
the Hilbert-Schmidt
operator $f \in L^2(\nu) \mapsto \int h(\cdot,y)f(y)d\nu(y)$,  
and $(\psi_i)_{i \geq 0} \subseteq  L^2(\nu)$
is a corresponding sequence of eigenfunctions.  
It follows that one may write for any $\mu \in \calP(\calX)$,
$$\sqrt{\Phi(\mu)} 
 = \sqrt{\iint h(x,y)d\mu(x)d\mu(y)} 
 = \sqrt{\sum_{i=0}^\infty \lambda_i \left(\int \psi_i d\mu\right)^2} 
 = \left\|\left(\sqrt{\lambda_i} \int\psi_id\mu \right)_{i\geq 0}\right\|_{\ell^2},$$
The right-hand side of the above display is a composition
of the 
convex map $\norm{\cdot}_{\ell^2}$ with the
affine map $\mu \in V \mapsto \left(\sqrt{\lambda_i}\int \psi_i d\mu\right)_{i\geq 0} 
\in \ell^2$, where $V$ is the vector space of finite signed
measures on $(\calX, \bbB(\calX))$. 
It follows that the functional $\sqrt{\Phi(\cdot)}$
is itself convex. 
Since this functional is also nonnegative, 
and the square function is convex and increasing
on $\bbR_+$, it is straightforward to verify that
the functional $\Phi(\cdot)$ is likewise convex. 
The claim then follows from Theorem~\ref{thm:submart}. 
\qed

%%%%%%%%%%%%%%% Wasserstein

\subsection{Proofs from Subsection~\ref{sec:wasserstein}}
{\bf Proof of Corollary~\ref{cor:wasserstein}.}
We will make use of the Kantorovich duality
(cf.  Section~\ref{sec:background_divergences})
% (see Theorem 1.3 and Remark 1.13 of \cite{villani2003}). 
% \begin{lemma}[Kantorovich Duality]
% \label{lem:kantorovich}
% Let $\Phi$ denote the  set of pairs $(f,g)$ of
% bounded and continuous
% functions $f,g:\calX \to \bbR$ such that 
% $0 \leq f \leq D^p$, $-D^p \leq g \leq 0$, 
% and 
% $f(x) + g(y) \leq \norm{x-y}^p$
% for all $x,y \in \calX$. Then, 
% $$W_p^p(P,Q) = \sup_{(f,g) \in \Phi} \int f dP + \int g d Q.$$
% Furthermore, the supremum is achieved by some pair $(f,g) \in \Phi$. 
% \end{lemma}
% We will use Lemma \ref{lem:kantorovich}
to show that the map
$$G:(x_1, \dots, x_t, y_1, \dots, y_s)\in \calX^{t+s} \mapsto
\calT_c\left(\frac 1 t \sum_{i=1}^t \delta_{x_i}, \frac 1 s \sum_{j=1}^s \delta_{y_j}\right)
 = \sup_{(f,g)\in \calM_c} \frac 1 t \sum_{i=1}^t f(x_i) + \frac 1 s \sum_{j=1}^s g(y_j),$$
satisfies the bounded differences property. 
This generalizes the one-sample analogue
proven for instance by \cite{weed2019} (Proposition 20). 
Let $1 \leq k \leq t$, 
and $x_1, \tilde x_1, \dots, x_t, \tilde x_t \in \calX$
be such that $\tilde x_i = x_i$ for all $i \neq k$.
Furthermore, let $y_1, \dots, y_s \in \calX$. Let
$(f_0, g_0) \in \calM_c$ denote optimal
Kantorovich potentials satisfying
$$\calT_c\left(\frac 1 t \sum_{i=1}^t \delta_{x_i}, 
\frac 1 s \sum_{j=1}^s \delta_{y_j}\right)
 = \frac 1 t \sum_{i=1}^t f_0(x_i) + \frac 1 s \sum_{j=1}^s g_0(y_j).$$
Furthermore, recall from Section~\ref{sec:background_divergences}
that we may (and do) choose
$(f_0,g_0)$ such that $0\leq f_0 \leq \Delta$
and $-\Delta \leq g_0 \leq 0$. Then,
\begin{align*}
G(x_1, &\dots, x_t,y_1, \dots, y_s) - G(\tilde x_1, \dots, \tilde x_t,y_1, \dots, y_s) \\
 &\leq \frac 1 t \sum_{i=1}^t f_0(x_i) 
     + \frac 1 s \sum_{j=1}^s g_0(y_j) - 
       \frac 1 t \sum_{i=1}^t f_0(\tilde x_i) - 
       \frac 1 s \sum_{j=1}^s g_0(y_j)
\leq \frac 1 t [f_0(x_k) - f_0(\tilde x_k)] \leq \Delta/t.
\end{align*}
Repeating a symmetric argument, we obtain 
$$|G(x_1, \dots, x_t,y_1, \dots, y_s) - G(\tilde x_1, \dots, \tilde x_t,y_1, \dots, y_s)|\leq \Delta/t.$$
We similarly have that for all $\tilde y_1, \dots, \tilde y_s\in \calX$,
satisfying $\tilde y_i = y_i$ for all $i \neq k$,
$$|G(x_1, \dots, x_t,y_1, \dots, y_s) - G(x_1, \dots,  x_t,\tilde y_1, \dots, \tilde y_s)|\leq \Delta/s.$$
We deduce that $G$ satisfies the bounded differences
property with parameters $L_1=\dots=L_t = \Delta/t$
and $L_{t+1} = \dots L_{t+s} = \Delta/s$. 
McDiarmid's inequality then implies
$$\bbP\big(|\calT_c(P_t, Q_s) - \bbE \calT_c(P_t,Q_s)| \geq u\big)
\leq 2\exp\left\{-\frac{2tsu^2}{(t+s)\Delta^{2}}\right\},\quad u > 0.$$
It follows that $\calT_c(P_t,Q_s)$ is
$\frac{2\Delta^{2}(t+s)}{ts}$-sub-Gaussian by 
Lemma~\ref{lem:subg}. 
Since $\calT_c$ is convex, applying Corollary
\ref{cor:sub_gaussian_bound} yields
\begin{align*}
\begin{multlined}[0.85\textwidth]
\bbP\Bigg(\exists t,s \geq 1: 
\calT_c(P_t,Q_s) - \calT_c(P,Q) \geq \alpha_{c,ts}
 \\ + 2\Delta\sqrt{\frac{2(t+s)}{ts} \Big[ \log g(\log_2 t + \log_2 s) + \log(2/\delta)\Big]} \Bigg) \leq \delta/2.
 \end{multlined}
 \end{align*}
Furthermore, Corollary~\ref{cor:ucb_bounded}
and the Kantorovich duality immediately
lead to the bound
$$\bbP(\exists t,s \geq 1: \calT_c(P_t,Q_s)-\calT_c(P,Q) \leq -\Delta \kappa_{ts}) \leq \delta/2.$$
The claim follows.
 \qed

 %%%%%%%%%%%%%%%%%%%%% KL

\subsection{Proofs from Subsection~\ref{sec:phi_divergence}}
{\bf Proof of Proposition~\ref{prop:kl_finite}.}
We repeat a similar stitching argument as that of the proof of Theorem~\ref{thm:main}.
Let $N_t = \KL(P_t\|P)$, $t \geq 1$. 
$(N_t)$ forms a reverse submartingale
by Lemma \ref{lem:convex} and Theorem~\ref{thm:submart},
implying by Jensen's inequality that for
%any given $\lambda \in (0,1]$ and
any integer $t_1 \geq 1$,
$\big(\exp(t_1 \lambda_{t_1} N_t)\big)_{t=1}^\infty$
is a reverse submartingale.
Therefore, 
following along similar lines as 
the proof of Theorem~\ref{thm:main}, and applying
Theorem~\ref{thm:reverse_ville}, we have for all
$y > 0$ and all integers $t_0 \geq 1$,
\begin{align*}
\bbP\Big(\exists t \geq t_0: N_t   \geq y\Big)
 &= \bbP\Big( \exists t \geq t_0: \exp(\lambda_{t_1} t_1 N_t) \geq \exp(\lambda_{t_1} t_1y) \Big) \\
 &\leq \bbE[\exp(-yt_1\lambda_{t_1} + t_1 \lambda_{t_1}N_{t_0})] \leq \exp(-yt_1 \lambda_{t_1} )G_{k,t_0}(t_1\lambda_{t_1}/t_0).
\end{align*} 
Now, letting $u_j = 2^j$ for all integers $j \geq 0$, and
$\gamma_t = \frac {1} {\lambda_{t}t} \log\left(\delta^{-1} G_{k,\lfloor t/2\rfloor }(2\lambda_t) \osfn(\log_2 t)\right)$, we obtain
\begin{align*}
\bbP\left(\exists t \geq 1: N_t \geq \gamma_t\right) 
  &\leq \bbP\left(\bigcup_{j=0}^\infty \Big\{\exists t \in \{u_j, \dots, u_{j+1}\}: N_t  \geq  \gamma_t\Big\}\right) \\
  &\leq \bbP\left(\bigcup_{j=0}^\infty \Big\{\exists t \in \{u_j, \dots, u_{j+1}\}: N_t  \geq  \gamma_{u_{j+1}}\Big\}\right) \\
 &\leq \sum_{j=0}^\infty \exp\left\{-u_{j+1}\gamma_{u_{j+1}} \lambda_{u_{j+1}}\right\} G_{k,u_j}(2\lambda_{u_{j+1}}) 
 \leq  \sum_{j=0}^\infty \frac \delta {\ell(j+1)} \leq \delta,
\end{align*}
where on the final line, we used
the fact that $G_{k,t}(\lambda)$ increases with $t$
for all fixed $k \geq 2, \lambda \in [0,1]$  (\cite{guo2020}, Lemma~1). 
The claim follows. \qed

 %%%%%%%%%%%%%%%%%%%%% TV

{\bf Proof of Corollary \ref{cor:tv}.}
By \cite{berend2013}, Eq. (17), and references
therein, we have
$$\bbP\Big(\norm{P_t-P}_{\mathrm{TV}} - \bbE \big[\norm{P_t-P}_{\mathrm{TV}}\big] \geq u\Big) \leq 2 \exp(-2tu^2),\quad u > 0,$$
implying that $\norm{P_t-P}_{\mathrm{TV}}$ is $(2/t)$-sub-Gaussian. Furthermore,
Lemma~7 of~\cite{kamath2015} implies
$$\bbE[\norm{P_t-P}_{\mathrm{TV}}] \leq \sqrt{\frac{k}{2\pi t}} + 2\left(\frac{k}{t}\right)^{\frac 3 4}.$$
%
% Eq. (5) of \cite{berend2013}
%implies $\bbE[\norm{P_t-P}_{\mathrm{TV}}]\leq \sqrt{k/t}/4.$
Finally, $\norm\cdot_{\mathrm{TV}}$ forms
a convex divergence by 
Lemma \ref{lem:convex}. The claim now follows by Corollary~\ref{cor:sub_gaussian_bound}.\qed

 \subsection{Proofs from Subsection~\ref{sec:smooth}}
We shall make use of the following
result due to \cite{bobkov1999}.
\begin{lemma}[\cite{bobkov1999}, Theorem 1.3]
\label{lem:inf_conv}
Let $d$ denote a metric on $\calX$. Then, 
a measure $\mu \in \calP_1(\calX)$ satisfies the 
$T_1(\sigma^2)$ inequality
with respect to $d$ if and only if $f(X)$ is $\sigma^2$-sub-Gaussian
for all functions $f:\calX \to \bbR$ which are 1-Lipschitz with respect to $d$.
\end{lemma}

{\bf Proof of Proposition~\ref{prop:smooth}.}
To prove part (i), it is straightforward
 to show that the map
$$G: (x_1, \dots, x_t) \in \bbR^{t \times d} \longmapsto  \norm{\frac 1 t \sum_{i=1}^t (\delta_{x_i} - P) }_{\mathrm{TV}}^\sigma,$$
satisfies the bounded differences property.
Indeed, given $1 \leq j \leq t$, 
let $x_1, \tilde x_1, \dots, x_t, \tilde x_t \in \bbR^d$, such that $x_i = \tilde x_i$
for all $i \neq j$. Then,
the triangle inequality implies
\begin{align*} 
|G(x_1, \dots, x_t) - G(\tilde x_1, \dots, \tilde x_t)| 
%  &= \Bigg| \sup_{A \in \bbB(\bbR^d)} 
%  \left| \left(\frac 1 t \sum_{i=1}^t \delta_{x_i} \star \calK_\sigma\right)(A) - (P\star \calK_\sigma)(A)\right|-
%   \sup_{A \in \bbB(\bbR^d)} 
%  \left| \left(\frac 1 t \sum_{i=1}^t \delta_{\tilde x_i} \star \calK_\sigma\right)(A) - (P\star \calK_\sigma)(A)\right|~\Bigg|\\
 &\leq \sup_{A \in \bbB(\bbR^d)} 
  \left| \left(\frac 1 t \sum_{i=1}^t \delta_{x_i} \star \calK_\sigma\right)(A) - \left(\frac 1 t \sum_{i=1}^t \delta_{\tilde x_i} \star \calK_\sigma\right)(A)\right| \\
 &\leq \frac 1 t \sup_{A \in \bbB(\bbR^d)}
 \int_A  \left| K_\sigma(x-x_j) - K_\sigma(x-\tilde x_j)\right| dx\\
 &\leq \frac 1 t \sup_{A \in \bbB(\bbR^d)}
 \int_A \Big[  K_\sigma(x-x_j) +  K_\sigma(x-\tilde x_j) \Big] dx \leq 2/t.
\end{align*}
Therefore, by McDiarmid's Inequality (Theorem~\ref{thm:mcdiarmid}), we have 
$$\bbP\Big(\big| \norm{P_t-P}_{\mathrm{TV}}^\sigma 
- \bbE \norm{P_t-P}_{\mathrm{TV}}^\sigma \big| \geq u\Big) 
\leq 2\exp(-tu^2/2), \quad u > 0.$$
It follows from Lemma \ref{lem:subg} that 
$\norm{P_t-P}_{\mathrm{TV}}^\sigma$ is $8t$-sub-Gaussian. 
Furthermore, \cite{goldfeld2020b} show that
$\bbE \norm{P_t-P}_{\mathrm{TV}}^\sigma 
\leq c_{d} t^{-1/2}/\sqrt 2$.
Finally, the Total Variation distance
is convex by Lemma \ref{lem:convex}, 
thus $\norm\cdot_{\mathrm{TV}}^\sigma$
is also convex. 
The first claim now follows from Corollary \ref{cor:sub_gaussian_bound}.

To prove the second claim, we show similarly
as~\cite{niles-weed2019} that the map
$$G: (x_1, \dots, x_t) \in \bbR^{t\times d}\longmapsto tW_1^\sigma\left(\frac 1 t \sum_{i=1}^t \delta_{x_i}, P\right),$$
is Lipschitz with respect to the metric $c_t(x,y) := \sum_{i=1}^t \norm{x_i-y_i}_2$ on $\bbR^{d \times t}$,
where $x=(x_1, \dots, x_t),y=(y_1, \dots, y_t) \in \bbR^{d\times t}$.
Let $\calJ$ denote the set of 1-Lipschitz
functions on $\bbR^d$, and recall
that the $W_1$ distance coincides with
the IPM generated by $\calJ$,
by the Kantorovich-Rubinstein duality. 
%the triangle inequality for $W_1$ implies
We thus have, by the triangle inequality for $W_1$,
\begin{align*}
|G(x) - G(y)|
 &\leq tW_1\left( \left(\frac 1 t \sum_{i=1}^t \delta_{x_i}\right) \star \calK_\sigma, \left(\frac 1 t \sum_{i=1}^t \delta_{y_i}\right) \star \calK_\sigma\right) \\
 &= t\sup_{f\in\calJ} \int f d\left(\frac 1 t \sum_{i=1}^t (\delta_{x_i} \star \calK_{\sigma} - \delta_{y_i} \star \calK_{\sigma})\right)\\
 %& (\text{By equation \eqref{eq:w1_dual}})\\
% &= \sup_{\norm f_{\text{Lip}}\leq 1} \frac 1 t \sum_{i=1}^t \int ( f \star K_\sigma) d(\delta_{x_i}   - \delta_{y_i}  ) \\ 
 &= \sup_{f\in\calJ} \sum_{i=1}^t \big[ (f \star K_\sigma)(x_i) - (f\star K_\sigma)(y_i) \big] \\
 &= \sup_{f\in\calJ} \sum_{i=1}^t 
 \int [f(x_i-z) - f(y_i-z)]K_\sigma(z)dz \\
 &\leq \sum_{i=1}^t \norm{x_i-y_i}_2
 \int  K_\sigma(z)dz = c_t(x,y).
\end{align*} 
% $\norm{f \star K_\sigma}_{\text{Lip}} \leq \norm f_{\text{Lip}}$.
We deduce that $G$ is $1$-Lipschitz with respect to $c_t$. Furthermore,
by \cite{gozlan2010}, Proposition 1.9, the product 
measure $P^{\otimes t}$ satisfies the $T_1(t\tau^2)$
inequality over $\bbR^{d\times t}$ with respect to 
$c_t$. Therefore, $G(X_1, \dots, X_t)=tW_1^\sigma( P_t, P)$ is $(t\tau^2)$-sub-Gaussian by Lemma \ref{lem:inf_conv}, i.e. $W_1^\sigma(P_t,P)$
is $(\tau^2/t)$-sub-Gaussian. 
Furthermore, 
$\bbE W_1^\sigma(P_{t}, P) \leq C_{d} t^{-1/2}/\sqrt 2$ by \cite{goldfeld2020b}. 
% Finally, the 1-Wasserstein distance 
% is convex by Lemma \ref{lem:convex}, 
% thus $W_1^\sigma$ is also convex. 
Applying Corollary \ref{cor:sub_gaussian_bound} leads to the claim.
\qed 

{\bf Proof of Corollary~\ref{cor:entropy}.}
Since $P$ is supported in $[-1,1]^d$, Proposition 5 of \cite{polyanskiy2016} implies
\begin{equation} 
\label{eq:entropy_wasserstein_bound}
 |h(P_t\star \calK_\sigma) - h(P\star \calK_\sigma)| \leq \frac 1 {2\sigma^2} \left(
|\mu_t-\mu| + 2\sqrt d  W_1^\sigma(P_t,P)\right),
\end{equation}
where $\mu_t = \int xdP_t(x)$ and $\mu = \int xdP(x)$. 
Notice that $P$ is 1-sub-Gaussian by Hoeffding's Lemma, and thus also satisfies the $T_1(1)$ inequality (by 
Lemma~\ref{lem:inf_conv}). 
By Corollary~\ref{cor:multivariate_means} (see also the discussion thereafter), we have
\begin{equation*} 
\forall t \geq 1 : |\mu_t - \mu| \leq    
2 \sqrt{\frac 1 t \Big[\log\ell(\log_2 t) + \log(4/\delta) \Big]},\quad
\text{with probability at least } 1-\delta/2,
\end{equation*}
and by Corollary~\ref{prop:smooth},
$$
\forall t \geq 1: W_1^\sigma(P_t, P) \leq \frac{C_{d}}{\sqrt t} +2\sqrt{\frac{1}{t} \Big[ \log \osfn(\log_2 t) + \log(2/\delta) \Big] },\quad \text{with probability at least } 1-\delta/2.
$$
By a union bound and equation~\eqref{eq:entropy_wasserstein_bound}, it follows that with probability at least $1-\delta$, 
we have uniformly in $t \geq 1$,
\begin{align*} 
|h(P_t \star \calK_\sigma) - h(P \star \calK_\sigma) | 
 &\leq \frac 1 {2\sigma^2} \left\{ |\mu_t - \mu| + 2\sqrt d W_1^\sigma(P_t, P) \right\} \\
 &\leq \frac 1 {2\sigma^2} \left\{ (2+4\sqrt d) \sqrt{\frac 1 t \Big[\log\ell(\log_2 t) + \log(4/\delta) \Big]}  + \frac{2\sqrt d C_{d}}{\sqrt t}  \right\} \\
 &\leq \frac {3\sqrt d} {\sigma^2}  \sqrt{\frac 1 t \Big[\log\ell(\log_2 t) + \log(4/\delta) \Big]}  + \frac{\sqrt d C_{d}}{\sqrt t\sigma^2},
\end{align*} 
as claimed. 
\qed

\subsection{Proofs from Subsection~\ref{sec:rademacher}} 
\label{app:pf_rademacher}
{\bf Proof of Corollary~\ref{cor:classification}.}
Notice that
\begin{align*}
\sup_{f \in \calF} |R(f) - R_t(f)|
 = \sup_{f \in \calF} \left| \int I(f(x) \neq y) d(P-P_t)(x,y)\right| 
 = D_\calJ(P_t\|P),
\end{align*}
where $D_\calJ$ is the IPM generated by the class
$\calJ = \{(x,y) \mapsto I(f(x) \neq y): f \in \calF\}$.
Since the functions in $\calJ$ are uniformly bounded by 1, 
if follows by the same argument
as in the proof of Corollary~\ref{cor:dkw}
that
% it is a straightforward observation that the map
% $$((x_1, y_1), \dots, (x_n,y_n)) \mapsto D_\calJ\left( 
% \frac 1 t \sum_{i=1}^t \delta_{(x_i,y_i)}, P\right)$$
% satisfies the bounded differences property with parameter
% $1/t$. It follows 
that $D_\calJ(P_t\|P)$
is $(2/t)$-sub-Gaussian.  
Furthermore, a standard symmetrization
argument (see for instance equation~(4.18)
of \cite{wainwright2019}) implies 
$$\bbE[D_\calJ(P_t\|P)] \leq 2 \calR_t(\calJ) = \calR_t(\calF),$$
where the final equality follows from 
Lemma~3.4 of~\cite{mohri2018}. 
By Corollary~\ref{cor:sub_gaussian_bound}, we deduce
$$\bbP\left( \exists t \geq 1:
D_\calJ(P_t\|P) \geq \calR_{\bar t}(\calF)
+ 2\sqrt{\frac 2 t \Big[ \log \ell(\log_2 t) + \log(1/\delta)\Big]}\right) \leq \delta,$$
which readily implies the first claim.
To prove the second, abbreviate $\hat \calR_t(\calF)$
by $\hat \calR_t$, and let 
$$\hat\calR_t^i = \bbE_{\bepsilon_{t+1}}\left[\sup_{f \in \calF}\frac 1 t \left|\sum_{\substack{j=1 \\ j\neq i}}^{t+1} \epsilon_j f(X_j)\right|\right],
\quad i=1, \dots, t+1.$$
Then, 
\begin{align*}
\hat\calR_{t+1}
 &= \bbE_{\bepsilon_{t+1}}\left[\sup_{f \in \calF}\frac 1 {t+1}\left| \sum_{j=1}^{t+1} \epsilon_j f(X_j)\right|\right] 
= \bbE_{\bepsilon_{t+1}}\left[\sup_{f \in \calF}\frac 1 {t+1}\left| \sum_{i=1}^{t+1} \frac 1 t\sum_{\substack{j=1 \\ j\neq i}}^{t+1} \epsilon_j f(X_j)\right|\right]  \\
&\leq  \frac 1 {t+1} \sum_{i=1}^{t+1}  \bbE_{\bepsilon_{t+1}}\left[\sup_{f \in \calF}\frac 1 t\left| \sum_{\substack{j=1 \\ j\neq i}}^{t+1} \epsilon_j f(X_j)\right|\right] 
= \frac 1 {t+1} \sum_{i=1}^{t+1} \hat\calR_t^i,
\end{align*}
implying that $(\hat\calR_t)$ satisfies the leave-one-out
property in equation~\eqref{eq:loo}. It follows
from Proposition~\ref{prop:loo}
that $(\hat\calR_t)$ is a reverse submartingale
with respect to the exchangeable filtration
$(\calE_t^X)$. 
Thus, Theorem~\ref{thm:main} and 
Corollary~\ref{cor:sub_gaussian_bound} can be 
invoked with $(N_t)$ replaced by $(\hat\calR_t)$. 
Furthermore, it can again be deduced as before
that $\hat {\calR}_t$ is $(8/t)$-sub-Gaussian,
and has mean $\calR_t$. Corollary~\ref{cor:sub_gaussian_bound}
thus leads to the claim. 
\qed 

\subsection{Proofs from Subsection~\ref{sec:multivariate_means}}
Let $\calA_\gamma$ be a $\gamma$-cover of $\bbS^{d-1}_\star$ of size $N_\gamma$. 
By a straightforward covering argument, 
notice that for any $\nu \in \bbS^{d-1}_\star$, 
there exists $\nu_0 \in \calA_\gamma$ such that $\norm{\nu-\nu_0}_* \leq \gamma$
thus
$$\nu^\top X = (\nu-\nu_0)^\top X+ \nu_0^\top X  \leq \gamma \norm X + \nu_0^\top X ,$$
whence
$$\norm X = \sup_{\nu\in \bbS^{d-1}_\star} \nu^\top X \leq \gamma \norm X  + \max_{\nu \in \calA_\gamma} \nu^\top X,$$
implying that $\norm X \leq \frac 1 {1-\gamma} \max_{ \nu \in \calA_\gamma} \nu^\top X.$
We deduce that for all $\lambda\geq 0$, 
\begin{align*}
\bbE\left[ \exp \left( \lambda \norm {\mu_t - \mu}\right) \right]
 &\leq \bbE\left[ \exp \left( \frac \lambda {1-\gamma} \max_{\nu \in \calA_\gamma} \nu^\top(\mu_t - \mu)\right) \right]\\
 &\leq \sum_{\nu \in \calA_\gamma}\bbE\left[ \exp \left( \frac \lambda {1-\gamma}   \nu^\top(\mu_t - \mu)\right) \right]\\ 
 &= \sum_{\nu \in \calA_\gamma} \left(\bbE\left[ \exp \left( \frac \lambda {t(1-\gamma)}   \nu^\top(X - \mu)\right) \right]\right)^t\\ 
 &\leq N_\gamma \exp \left(t \psi\left(\frac{\lambda}{t(1-\gamma)} \right)\right).
%  \leq  \left(\sum_{\nu_0 \in \calA} \bbE\left[ \exp \left( \frac \lambda {n(1-\gamma)} \nu_0^\top(X - \mu)\right) \right]\right)^n \\ 
%   \leq  \left(|\calA| \bbE\left[ \exp \left( \frac \lambda {n(1-\gamma)} \nu_0^\top(X - \mu)\right) \right]\right)^n,
  \end{align*}
  where we extend the definition of the convex function 
$\psi$ to $\bbR_+$ by setting $\psi(\lambda) = \infty$
for all $\lambda \geq \mlambda$. 
We deduce that an upper bound on the cumulant generating function of $\norm{\mu_t-\mu}$ is given by
$$\widebar\psi_t(\lambda) = \log N_\gamma + t \psi\left(\frac{\lambda}{t(1-\gamma)} \right), 
\quad \lambda \geq  0,$$
It is readily seen that for all $x \in \bbR$ and all
$\lambda \geq 0$,
$$\widebar \psi_t^*(x) = -\log N_\gamma + t \psi^*\big((1-\gamma)x\big), \quad 
  \big(\widebar \psi_t^*\big)^{-1}(\lambda) = 
  \frac 1 {1-\gamma} (\psi^*)^{-1}\left(\frac{\lambda  + \log N_\gamma}{t}\right).$$
Finally, notice that the functional $\Phi(Q) = \norm{\int x dQ(x) - \mu}$ is convex, 
thus we may apply Theorem~\ref{thm:main} to deduce that for all $\delta \in (0,1)$, 
$$\bbP\left\{ \exists t \geq 1 : \norm{\mu_t-\mu} \geq \frac 1 {1-\gamma} (\widebar  \psi^*)^{-1}\left(\frac{\log\ell(\log_2t) + \log(1/\delta) + \log N_\gamma}{\lceil t/2\rceil }\right)\right\}\leq \delta.$$
The claim follows.  
\qed

\subsection{Proofs from Subsection~\ref{sec:lil}}
{\bf Proof of Corollary \ref{cor:lil}.}
Let $\eta,\alpha > 1$, and set
$u_k = \lceil \eta^k \rceil$
for all $ k\geq 0$. Define $\ell(k) = (1\vee k^\alpha) \zeta(\alpha)$,
where $\zeta(\alpha) = \sum_{k=1}^\infty \frac 1 {k^\alpha}$ and $\alpha > 1$. From the proof
of Theorem
\ref{thm:main_better},
for the process $N_t = D(P_t\|P)$ 
in the special case
of sub-Gaussian tails
$\psi_t(\lambda) = \lambda \bbE (N_t) + \lambda^2\sigma^2/2t$,
it can be seen that {$\bbP(A_k) \leq \delta/\osfn(k+1)$},
where, for all $k=0,1,\dots$, 
$$A_k = \left\{\exists u_k \leq t \leq u_{k+1}:
N_t > \bbE \left(N_{\lceil t/\lceil\eta\rceil\rceil}\right)+ \sqrt{\frac{2\sigma^2}{\lceil t/\lceil\eta\rceil\rceil }\Big[\log\osfn(\log_\eta t) + \log(1/\delta)\Big]}\right\}.$$
Thus, by definition of $\ell$ and by
the first Borel-Cantelli 
Lemma, we have $\bbP(\limsup_{k\to\infty} A_k) = 0$.
Therefore,
$$\bbP\left\{ 
N_t \leq
\bbE \left(N_{\lceil t/\lceil\eta\rceil\rceil}\right)  + 
\sqrt{\frac {2\sigma^2} {\lceil{t/\lceil\eta\rceil \rceil}}\Big[\log \osfn(\log_\eta t) + \log(1/\delta)\Big]} \ \text{eventually}\right\} =1.$$
Note that for all $t \geq \eta$,
\begin{align*}
\log \ell(\log_\eta t)
 = \alpha \log\log_\eta t  + \log\zeta(\alpha)
 = \alpha\log\log t - \alpha \log\log\eta + \log \zeta(\alpha).
 \end{align*} 
Therefore,
we have almost surely,
\begin{align*}
\limsup_{t \to \infty}
\frac{D(P_t\|P)}
{\bbE\left( N_{\lceil t/\lceil\eta\rceil\rceil}\right) + \sqrt{\frac {2 \sigma^2} {\lceil{t/\lceil\eta\rceil\rceil}}
\Big[\alpha\log \log t - \alpha\log\log \eta  +\log\zeta(\alpha)+ \log(1/\delta)\Big]}} \leq 1,
\end{align*} 
whence, by assumption on $\bbE N_t$
and by the fact that $\delta,\alpha,\eta$ are fixed,
we obtain
\begin{align*}
\limsup_{t \to \infty}
\frac{D(P_t\|P)}
{\sqrt{\frac {2 \sigma^2} {\lceil{t/\lceil\eta\rceil\rceil}}
 \alpha\log \log  t  }} \leq 1\quad \text{a.s.}
\end{align*} 
Since we may choose $\eta$ and $\alpha$
arbitrarily close to 1, the claim follows.
\qed

 \section{The failure of forward submartingales and canonical filtrations}
\label{sec:forward_submart}
Our main results in Section~\ref{sec:main}
derived confidence sequences
for convex divergences on the basis
of their empirical plug-in
estimators, which form
reverse submartingales with respect
to the exchangeable filtration. 
The latter are rarely used 
in sequential analysis. The aim
of this section is
to illustrate our difficulty in using more typical 
tools from sequential analysis.

Let $D$ be a convex divergence, and $(X_t)_{t=1}^\infty$
a sequence of i.i.d. observations from a distribution
$P \in \calP(\calX)$. We shall
focus on deriving time-uniform concentration inequalities
for the process 
\[
S_t = tD(P_t\|P), \quad t \geq 1.
\]
% As in Section~\ref{sec:main},
% we tacitly assume that
% $S_t$ is measurable for all $t \geq 1$,
% throughout the sequel. 
Specifically, given a level $\delta \in (0,1)$,
we derive sequences of nonnegative
real numbers $(u_{t})_{t=1}^\infty$
such that 
\begin{equation}
\label{eq:prelim_goal}
\bbP\Big(\exists t \geq 1: S_t \geq u_{t} \Big) \leq \delta.
\end{equation}
% for a given level $\delta \in (0,1)$, and for a nonnegative sequence $(u_{t,\delta})_{t=1}^\infty$. 
Assuming $D$ is additionally symmetric, satisfies the triangle
inequality, and that equation \eqref{eq:prelim_goal} 
can also be applied to the process $(sD(Q_s\|Q))_{s=1}^\infty$, 
one also obtains the bound
\begin{equation}
\label{eq:prelim_twosample}
\bbP\Big(\exists t,s \geq 1: |D(P_t\|Q_s) - D(P\|Q)| \geq u_{t}/t + u_{s}/s \Big) \leq 2\delta.
\end{equation}
As was discussed in Section~\ref{sec:main}, the passage from equation \eqref{eq:prelim_goal} to \eqref{eq:prelim_twosample}
may lead to confidence sequences 
with lengths of sub-optimal rate in general.

In Section \ref{sec:background_confseq},
we discussed two prominent
stitching constructions used in past
work to derive confidence
sequences for the common mean $\mu$
of $(X_t)$. 
Assuming 
$X_1-\mu$ admits a finite cumulant generating
function $\phi: [0,\mlambda) \to \bbR$
for some $\mlambda > 0$,
these approaches are either based on 
\begin{enumerate}
    \item[(A)] applying Ville's inequality to the  
forward nonnegative supermartingale $(L_{t}(\lambda))$, or
% $\big(\exp\{\lambda \sum_{i=1}^t (X_i-\mu) - t \phi(\lambda)\big]\big)_{t=0}^\infty$,
\item[(B)] applying Doob's submartingale inequality
to the forward submartingale $(U_t(\lambda))$.
\end{enumerate}
While concrete applications of
these approaches admit various nuances which lead 
to confidence sequences with distinct 
constants, under varying assumptions
on $\phi$, they both hinge upon
repeated applications of a maximal 
forward martingale
inequality over geometrically-increasing epochs of 
time, and lead to confidence sequences
of length typically scaling at the rate
$O\big(\sqrt{\log\log t/t}\big)$. 
We show in what follows that analogues
of constructions (A) and (B) can be obtained
for the process $(S_t)$, but may
be unsatisfactory in general. 

We begin with the following simple
result, which provides an analogue of the exponential
supermartingale in (A) for the process $(S_t)$.
In what follows, we let
$\calS_t = \sigma(X_1, \dots, X_t)$, 
$t \geq 1$,
denote the canonical filtration. 
\begin{lemma} 
\label{lem:convex_supermart}
Let $D$ denote a convex divergence, 
and set $S_t = tD(P_t\|P)$ for all 
$t \geq 1$. Assume there exists
$\mlambda > 0$ and a convex function $\psi$ such that 
\begin{equation} 
\label{eq:forward_supermart_tail}
\log\Big\{\bbE[\exp(\lambda D(\delta_{X_1}\|P))] \Big\} 
\leq \psi(\lambda),\quad  
\lambda \in [0,\mlambda).
\end{equation}
Then, for any $\lambda \in [0,\mlambda)$,
the process 
$$\calL_t(\lambda) = \exp\Big\{\lambda S_t - t \psi(\lambda)\Big\},\quad t \geq 1$$ 
is a nonnegative \underline{super}martingale with respect to $(\calS_t)_{t=1}^\infty$.
 \end{lemma}
Lemma \ref{lem:convex_supermart} is proven 
below. 
The process $(\calL_t(\lambda))$ is in direct
analogy with the supermartingale in 
(A), where the cumulant generating function $\phi$ of 
 $(X_1-\mu)$ is
replaced by that of $D(\delta_{X_1} \| P)$. 
 Unlike (A), however, the quantity
 $D(\delta_{X_1} \| P)$ has non-zero mean; for example if $D$ is the KL divergence and $P$ has finite support, then $\bbE[D(\delta_{X_1} \| P)]$ is the
 entropy of $P$.
 This in turn implies that the upper bound
 $\psi$ may typically
 decay linearly as $\lambda \downarrow 0$,
 as opposed to the typical quadratic
 rate of decay for cumulant
 generating functions \citep{howard2020}.
 It can be seen that standard stitching
 constructions applied to this process 
 would then only lead
 to a confidence sequence for $S_t/t = D(P_t\|P)$
 with non-vanishing length, which is
 overly conservative in most examples of interest.
 
 This intractibility of the process
 $(L_t(\lambda))$ arises 
 from the fact 
 that, unlike the setting
 (A) for sums of i.i.d. random variables,
 the term $t\psi(\lambda)$
can be an excessively  loose upper bound on 
 the true cumulant generating function
 of $tD(P_t\|P)$. While the form of this upper bound
 is pivotal for deriving the supermartingale property
 of $(L_t(\lambda))$, it would ideally
 be replaced by a tighter upper bound
 $\psi_t:[0,\mlambda) \to \bbR$, which is not
 necessarily linear in $t$, 
 satisfying 
\begin{equation}
\label{eq:forward_submart_tail}
\log \big\{ \bbE \big[ \exp(\lambda S_t) \big] \big\} \leq \psi_t(\lambda), \quad \lambda \in [0,\lambda_{\text{max}}).
\end{equation}
Without restrictions on the form of $\psi_t$, 
 however, it is unclear how to enforce a supermartingale property
 for the process $(\exp\{\lambda S_t - \psi_t(\lambda)\})$ for a general divergence $D$.
 
 These considerations motivate
 approach (B), though unlike Lemma \ref{lem:convex_supermart}, 
 however, we are not
 aware of a general result guaranteeing
 that $(\exp(\lambda S_t))$ is a forward
 submartingale for any convex
 divergence $D$. The following straightforward
 result presents a notable exception.
 \begin{lemma} 
\label{lem:ipm_submart}
Let $\calJ$ be any class of Borel-measurable functions from $\bbR^d$ to $\bbR$, let $D$ be the IPM
generated by $\calJ$, and recall $S_t = t D(P_t\| P)$. 
% Let $\calS_t = \sigma(\{S_0, \dots, S_t\})$, $t \geq 0$ denote the canonical
% filtration. 
Assume $S_t\in L^1(\bbP)$
for all $t \geq 1$. 
Then, $(S_t)_{t=1}^\infty$ is a forward nonnegative \underline{sub}martingale with respect to $(\calS_t)_{t=1}^\infty$.
In particular, for any $\lambda \in[0,\mlambda)$,
the process $\big(\exp(\lambda S_t)\big)_{t=1}^\infty$
is also a forward nonnegative submartingale with respect
to $(\calS_t)$ whenever it lies in $L^1(\bbP)$. 
\end{lemma}
We use Lemma \ref{lem:ipm_submart}
to derive a confidence sequence
for the process $(S_t)$, in the case
where $D$ is an IPM, or any other
convex divergence for which $(\exp(\lambda S_t))$
forms a forward martingale. 
Our result will depend on a stitching function
$\ell: [0,\infty) \to [1,\infty)$,
which satisfies $\sum_{k=0}^\infty 1/\ell(k) \leq 1$.
% and will dictate the shape of our confidence sequences which follow.
% A typical choice is the function
% $\ell(k) =  (k+1)^\alpha \zeta(\alpha)$,
% where $\alpha > 1$ and $\zeta(\alpha) = \sum_{k=1}^\infty (1/k^\alpha)$.
We also assume that the upper bound $\psi_t$ 
is convex. 
The main result of this section is stated as follows. 
\begin{theorem}
\label{thm:forward_submart}
Let $D$ be a convex divergence, and let $S_t=tD(P_t\|P)$
for all $t \geq 0$. Assume that for any 
$\lambda \in [0,\mlambda)$, the process $\big(\exp(\lambda S_t)\big)_{t=0}^\infty$ is a 
forward submartingale with respect to a filtration $(\calS_t)_{t=0}^\infty$, and satisfies inequality
\eqref{eq:forward_submart_tail}. 
Assume $\psi_t^*$ is invertible in $\lambda$ and
that the sequence
$$\gamma_t =  (\psi_{t}^*)^{-1} \left( \log \osfn \left(\log_2 t\right) + \log(2/\delta)\right),
\quad t \geq 1$$
%Define for all $t \geq 1$,
%$$\gamma_t = (\varphi_{\lceil t(\eta-1)\rceil}^*)^{-1} \left( \log h\left(\log_\eta t\right) + \log(2/\delta)\right).$$
is increasing in $t$. 
Then, 
for all $\delta \in (0,1)$, 
\begin{equation}
\label{eq:forward_submart_general_bound}
 \bbP\Big( \exists t \geq 1: S_t \geq 2\gamma_t \Big) \leq \delta.
\end{equation}
\end{theorem}  

Theorem~\ref{thm:forward_submart} is proved by dividing
time into geometrically
increasing epochs of the form $[2^j, 2^{j+1}-1]$,
$j \geq 0$,
with the goal of applying Doob's submartingale
inequality over each epoch, at level $\delta/\ell(j)$.
Applying a union bound over these epochs leads to the claim.
% It can be seen from the form of Doob's ienquality
% that doing so directly would be equivalent
% do dividing time into the epochs
% $[0,2^{j+1}-1]$, which we believe is loose
% due to their nestedness. 
Our proof makes use of the convexity of $D$, which
guarantees the following bound on the deviations of $(S_t)$,
\begin{equation}
\label{eq:key_inequality}
S_{t_1} - S_{t_0} \leq S_{t_0,t_1}:=  (t_1-t_0)  D\left(\frac 1 {t_1-t_0}\sum_{i=t_0+1}^{t_1} \delta_{X_i} \Bigg\| P\right),\quad \forall t_1 > t_0.
\end{equation}
Inequality \eqref{eq:key_inequality} is used in our proof
to relate
upcrossing probabilities over 
the interval $[2^j, 2^{j+1}-1]$
to those over its translation at zero, $[0, 2^{j+1}-1 - 2^j]$,
leading to improved constants in Theorem~\ref{thm:forward_submart}, and can be viewed
as a generalization of the proof technique employed
in Lemma 3 of \cite{jamieson2014}
for processes $(S_t)$ which are submartingales
rather than martingales.

While the bound of Theorem~\ref{thm:forward_submart}
is useful, it presents
two key shortcomings. On the one hand, 
it does not directly lead to a confidence sequence
for the quantity of interest $D(P\|Q)$,
as discussed at the beginning of this Appendix. 
% As was discussed at the beginning of this
% section, a confidence sequence for $D(P\|Q)$
% can be obtained from Theorem~\ref{thm:forward_submart}
% for divergences $D$ which satisfy
% the triangle inequality, but doing so
% may lead to suboptimal rates for certain
% divergences---we present such an example
% for Wasserstein distances
% in Section \ref{sec:wasserstein}.  
On the other hand, 
Theorem~\ref{thm:forward_submart}
can only be applied when $(\exp(\lambda S_t))$
forms a forward submartingale. 
We have shown in Lemma \ref{lem:ipm_submart} 
that $(\exp(\lambda S_t))$, 
is a submartingale when $D$ is an IPM, but we are not
aware of a generic result of this kind
for other convex divergences. 
In contrast, we showed in Section~\ref{sec:main} that
for any convex divergence, this process,
and the more general process $(M_{ts})$
in equation \eqref{eq:process_Mts},
form \textit{reverse}
submartingales with respect to the exchangeable
filtration, which turned out to be more natural to 
handle.

We close this Appendix by proving 
Lemmas~\ref{lem:convex_supermart}, 
\ref{lem:ipm_submart}, and 
Theorem~\ref{thm:forward_submart}.
\subsection{Proofs}
{\bf Proof of Lemma \ref{lem:convex_supermart}.}
$(\calL_t(\lambda))$ is clearly adapted to
$(\calS_t)$. Now, notice that 
%  we may write $P_t$ as the following convex combination
$P_t = \frac 1 t \delta_{X_t} + \frac{t-1}{t} P_{t-1}.$
Therefore, by convexity of $D$, 
$$D(P_t\|P) \leq  \frac{t-1}{t} D(P_{t-1}||P)
+ \frac 1 t D(\delta_{X_t}\|P) , \quad\text{thus,}
\quad S_t \leq S_{t-1} + D(\delta_{X_t}\|P).$$
It follows that for all $t \geq 1$,
\begin{align*}
\bbE[\calL_t(\lambda)|\calF_{t-1}] 
 &= \bbE \left\{ \exp\big(\lambda S_t - t\psi(\lambda)\big)~\big|~ \calS_{t-1}\right\}\\
 &\leq \bbE \left\{ \exp\big(\lambda S_{t-1} + \lambda D(\delta_{X_t}\|P) - t\psi(\lambda)\big)~\big|~ \calS_{t-1}\right\}\\
 &= \calL_{t-1}(\lambda) \bbE\left\{ \exp[\lambda D(\delta_{X_t}\|P) - \psi(\lambda)]\right\} 
 \leq \calL_{t-1}(\lambda),
\end{align*}
where the last inequality follows from the definition 
of $\psi$.
\qed

{\bf Proof of Lemma \ref{lem:ipm_submart}.}
The claim is straightforward. $(S_t)$ is clearly adapted to $(\calS_t)$. Furthermore, for all $t \geq 1$, notice that
\begin{align*}
\bbE[S_t|\calS_{t-1}] 
 &= \bbE\left[t \sup_{f \in \calJ} \int f d(P_t-P)\bigg| \calS_{t-1}\right] \\
 &\geq \sup_{f \in \calJ} \bbE\left[\sum_{i=1}^t [f(X_i) - \bbE f(X_i)]\bigg| \calS_{t-1}\right] 
 = \sup_{f \in \calJ} \sum_{i=1}^{t-1} [f(X_i) - \bbE f(X_i)]
 = S_{t-1}.
\end{align*}
The final claim follows from the convexity and monotonicity of the exponential function.
\qed

{\bf Proof of Theorem~\ref{thm:forward_submart}.}
The proof is inspired by 
\cite{jamieson2014} (Lemma 3).
By assumption, $(\exp(\lambda S_t))$  forms a submartingale with respect to the filtration $(\calS_t)$, for any $\lambda \in [0,\mlambda)$.
By Doob's submartingale inequality, we therefore have for any integer $T \geq 1$ and real number $y > 0$
the Cram\'er-Chernoff bound
\begin{align*} 
\bbP\Big( \exists t \leq T: S_t \geq y\Big) 
 &= \inf_{\lambda \in [0,\mlambda)} \bbP\Big( \exists t \leq T: \exp(\lambda S_t) \geq \exp(\lambda y) \Big) \\
 &\leq \inf_{\lambda \in [0,\mlambda)} \bbE\big[\exp\big(\lambda S_T - \lambda y\big)\big] \\
 &\leq \inf_{\lambda \in [0,\mlambda)} \bbE\big[\exp\big(\psi_T(\lambda) - \lambda y\big)\big] 
 = \exp\{-\psi_T^*(y)\}.
\end{align*} 
Now, set $u_k = 2^k$ for all $k \geq 0$.
Since $(\gamma_t)$ is an increasing sequence, we have
 \begin{align}
\nonumber &\bbP\left(\exists t \geq 1: S_t \geq 2\gamma_t\right) 
%\nonumber  &=\bbP\left(\bigcup_{k\in \bbN} \Big\{\exists t \in \{u_k+1, \dots, u_{k+1}\}: M_t - 2\bbE[M_t]\geq 2\gamma_t\Big\}\right)\\
\leq \bbP\left(\bigcup_{k=0}^\infty  \Big\{\exists t \in \{u_k, \dots, u_{k+1}-1\}: (S_t - S_{u_k}) + S_{u_k} \geq 2\gamma_{u_k}\Big\}\right)\\
% &(\text{Since } \gamma_t \text{ is increasing}) \\
%\nonumber  &\leq \bbP\left(\bigcup_{k\in \bbN} \Big\{\exists t \in \{u_k+1, \dots, u_{k+1}\}: M_{u_k,t} + M_{u_k} \geq 2\gamma_{u_k}\Big\}\right)
%   & (\text{By \eqref{eq::bound_fact1})}\\
\nonumber  &\leq \bbP\left(\bigcup_{k=0}^\infty \Big\{\exists t \in \{u_k, \dots, u_{k+1}-1\}: S_{u_k,t} + S_{u_k}
 \geq 2\gamma_{u_k} \Big\}\right),\qquad  \text{(by equation \eqref{eq:key_inequality})} \\
  &\leq \bbP\left(\bigcup_{k=0}^\infty \Big\{S_{u_k} \geq \gamma_{u_k}  \Big\}\right)+
   \bbP\left(\bigcup_{k=0}^\infty \Big\{\exists  t \in \{u_k, \dots, u_{k+1}-1\}: S_{u_k,t}   \geq \gamma_{u_k}\Big\} \right).
  \label{eq:pf_midway_expectations} 
\end{align}
We will upper bound the terms in \eqref{eq:pf_midway_expectations} separately. For the first term, by a union bound and another
standard application of
the Cram\'er-Chernoff technique,  
we have
\begin{align*}
\bbP\left(\bigcup_{k=0}^\infty \Big\{S_{u_k}  \geq \gamma_{u_k}  \Big\}\right) 
 \leq \sum_{k=0}^\infty \bbP\Big(S_{u_k} \geq \gamma_{u_k} \Big)  
 \leq\sum_{k=0}^\infty \exp\big\{-\psi_{u_k}^*(\gamma_{u_k})\big\}
 =\frac \delta 2 \sum_{k=0}^\infty \frac 1 {\osfn(k)}
 &\leq \frac \delta 2.
\end{align*}
For the second term in equation \eqref{eq:pf_midway_expectations}, notice that $S_t = S_{0,t} \overset{d}{=} S_{r,t+r}$ for all $t,r \geq 1$. 
Therefore, recalling that $u_{k+1}=2u_k$ we have for all $k\in\bbN$,
\begin{align}
\nonumber \sum_{k=0}^\infty
\bbP\Big(\exists t \in \{u_k, \dots, u_{k+1}-1\}: S_{u_k,t} \geq \gamma_{u_k}\Big) 
\nonumber &= \sum_{k=0}^\infty 
\bbP\Big(\exists t \in \{u_k, \dots, u_{k+1}-1\}:  S_{0,t-u_k} \geq \gamma_{u_k}\Big) \\
\nonumber &= \sum_{k=0}^\infty 
\bbP\Big(\exists u \in \{1, \dots, u_{k+1}-u_k\}: S_{0,u}  \geq \gamma_{u_k}\Big) \\
\nonumber  &= \sum_{k=0}^\infty 
\bbP\Big(\exists u \in \{1, \dots, u_k\}: S_u  \geq \gamma_{u_k}\Big)\\
 \nonumber &\leq \sum_{k=0}^\infty \exp\left\{ -\psi^*_{ u_k}(\gamma_{u_k}) \right\}
 \nonumber = \frac \delta 2 \sum_{k=0}^\infty \frac 1 {\osfn(k)} 
 \nonumber \leq \frac \delta 2.
\end{align} 
The claim thus follows. \qed

\end{document}